\documentclass{article}
\usepackage{epsfig}
\usepackage{hyperref}

\title{The Eckmann-Hilton argument and higher operads}

\author{M.A. Batanin\protect \footnote{The author holds the Scott Russell Johnson Fellowship in
the Centre of Australian Category Theory at Macquarie University.}\\ Macquarie University,  NSW 2109, Australia \\
e-mail: mbatanin@math.mq.edu.au}

\date{March 15, 2006}

\newtheorem{theorem}{\bf Theorem}[section]
\newtheorem{defin}{\bf Definition}[section]
\newtheorem{pro}{\bf Proposition}[section]
\newtheorem{lem}{\bf Lemma}[section]
\newtheorem{corol}{\bf Corollary}[pro]
\newtheorem{cor}{\bf Corollary}[theorem]

\newcommand{\ten}[1]{\mbox{$\diamond_{\scriptscriptstyle #1}$}}
\newcommand{\A}{\mbox{$\cal A$}}
\newcommand{\F}{\mbox{$\cal F$}}

\newcommand{\Q}{   \hfill  \ $\scriptstyle \clubsuit$ \ }

\newcommand{\h}{\mbox{$\bf h$}}
\renewcommand{\H}{\mbox{$\bf H$}}
\newcommand{\M}{\mbox{${\bf \tilde{m}}$}}
\newcommand{\ML}{\mbox{$\bf m$}}

\newcommand{\Proof}{\noindent {\bf Proof. \ }}
   
\renewcommand{\theequation}{\thesection.\arabic{equation}}
\makeatletter\@addtoreset{equation}{section}\makeatother

\newcommand{\Example}{\noindent \makebox[25mm]{{\bf
 Example  \hspace{-3mm}
\addtocounter{example}{1} 
\thesection.
\hspace{-3.5mm}
\theexample \ \  }}}

\newcommand{\Remark}{\noindent \makebox[23mm]{{\bf Remark
\hspace{-1mm}\addtocounter{remark}{1} 
\thesection.\theremark \ }}}
\newcounter{remark}[section]
\newcounter{example}[section]

\begin{document}

\maketitle

\begin{flushright}\it To the memory of my father. 
\end{flushright}

\begin{abstract}

The classical Eckmann-Hilton argument shows that  two 
monoid structures on a set, such that one is a
homomorphism for the other, coincide and, 
moreover, the resulting monoid is commutative. This
argument immediately gives a proof of the commutativity of the higher
 homotopy groups. A reformulation of this argument in
the language of higher categories is: suppose we have a one object, 
one arrow 2-category, then its $Hom$-set is
a commutative monoid. A similar argument due to A.Joyal and R.Street 
shows that a one object, one arrow tricategory is `the same' as a braided monoidal category.

In this paper we begin to investigate how one can extend this argument to arbitrary dimension. 
We provide a simple categorical scheme which allows us to formalise the Eckmann-Hilton type argument in terms of the calculation of left Kan extensions in an appropriate $2$-category. Then we apply this scheme to the case of $n$-operads in the author's sense and classical symmetric operads. 
We  demonstrate that there exists a functor of symmetrisation $Sym_n$ from  a certain subcategory of  $n$-operads  to the category of  symmetric operads
 such that the
category of one object, one arrow , . . . , one $(n-1)$-arrow algebras of $A$ is isomorphic  to the
category of  algebras
 of $Sym_n(A)$.    Under some mild conditions, we present
an explicit formula for $Sym_n(A)$ which involves taking the colimit over a remarkable categorical symmetric operad. 
 
 We will consider some applications of the methods developed to the theory of $n$-fold loop spaces  in the second paper of this series. 

\

1991 Math. Subj. Class.  18D05 , 18D50, 55P48     
\end{abstract}

\pagebreak

\tableofcontents

\section{How can symmetry emerge from nonsymmetry ?}\label{introduction}

Hopf and Alexandrov pointed out to \u{C}ech that his higher homotopy groups were commutative.  The proof follows
from the following statement which is known since \cite{EH} as  the  Eckmann-Hilton argument:  two 
monoid structures on a set such that one is a homomorphism for  the other coincide and, 
moreover, the resulting monoid is commutative.  A reformulation of this argument in
the language of higher categories is: suppose we have a one object, 
one arrow 2-category, then its $Hom$-set is
a commutative monoid. A higher dimensional generalization of this argument
 was provided by Joyal  and Street in
\cite{JS}. Essentially they proved that a 1-object, 1-arrow tricategory 
is a braided monoidal category and a one object,
one arrow, one 2-arrow tetracategory is a symmetric monoidal category.

Obviously we have here a pattern of some general higher categorical principle.   
Almost nothing, however, is known precisely except for   the above low dimensional examples and some higher
dimensional cases which can be reduced to the classical Eckmann-Hilton argument \cite{CransA}. Yet, there are
plenty of important conjectures which can be seen as  different manifestations of this principle. First of all
there is a bunch of hypotheses from Baez and Dolan
\cite{BD} about the so called `$k$-tuply monoidal'  $n$-categories, which are  $(n+k)$-categories with one
object, one arrow etc. up to
$(k-1)$.
 Basically these hypotheses state that these `$k$-tuply
monoidal'  $n$-categories are $n$-categorical analogues of $k$-fold loop spaces i.e. $n$-categories equipped with
an additional monoidal structure together with some sort of higher symmetry structures similar to 
the structure of a $k$-fold loop space. In particular, `$k$-tuply
monoidal' weak $\omega$-groupoids should model $k$-fold loop spaces. Many other hypotheses from \cite{BD} are
based on this analogy. 

Another problem, which involves the passage to  $k$-tuply monoidal $n$-categories, is the definition of higher
centers
\cite{BD,Crans,Simpson}. Closely related to this problem is 
 the Deligne conjecture from deformation theory \cite{K,KS} which tells us that there is an action of an
$E_2$-operad on the Hochschild  complex of an associative algebra. This conjecture is now proved by several
people. In higher dimensions the generalised Deligne conjecture was understood by Kontsevich \cite{K} as a
 problem of existence of some sort of homotopy  centre of any
$d$-algebra. This homotopy centre must have a structure of 
$(d+1)$-algebra. Here a $d$-algebra is an algebra of the  little $d$-cubes operad \cite{May}. To the best of our
knowledge this hypothesis is not proved yet in full generality, but there is  progress on it \cite{Tam}. 

In this paper we  consider a categorical basis for  the Eckman-Hilton argument in higher dimensions   using the
apparatus of higher-dimensional nonsymmetric operads  \cite{BatN}. They were introduced in \cite{BatN} 
for the purpose of defining weak $n$-categories for higher $n$. A weak $n$-category in
our sense is an algebra of a contractible (in a suitable combinatorial sense) $n$-operad.

Now consider the algebras of
an
$n$-operad
$A$ which have only one object, one arrow, ... , one $(k-1)$-arrow. The underlying $n$-globular
 object of such an
algebra can be identified with  an $(n-k)$-globular object and we can ask ourselves  what
sort of algebraic structure the action of $A$ induces on this $(n-k)$-globular object. 
Here we restrict ourselves  by considering
only
$k=n$. This provides  a great  simplification of the theory, yet clearly shows how  higher
symmetries can appear.  
We must say that we 
do not know the answer for arbitrary $k$. For this, perhaps, we need to develop  the theory of
symmetric higher operads, and some steps in this direction have already been taken in \cite{W}.

Returning to the case $k=n$  we  show that for an
$n$-operad $A$  one can construct a symmetric operad $Sym_n(A)$ (which in this case is just a classical
symmetric operad in the sense of May \cite{May} in a symmetric monodal category), called
symmetrisation  of
$A$, such that the category of one object, one arrow , ... , one $(n-1)$-arrow algebras of $A$ is isomorphic  to
the category of  algebras of
$Sym_n(A)$. Moreover,
 under mild conditions we present
an explicit formula for $Sym_n(A)$  involving the colimit over a remarkable categorical symmetric operad.

Fortunately, the restriction $n=k$ not only simplifies our techniques, but also makes almost unnecessary the use
of variable category theory from \cite{StV,StP}  which our paper \cite{BatN} used.
We can reformulate our theory of higher operads in a way that makes it very similar to the theory of
classical symmetric operads. So the reader who does not need to understand the full structure of a higher operad
may read the present paper without looking at \cite{BatN,StR,StP,BS}. In several places we do refer to some
constructions from
\cite{BatN,BS} but these references are not essential for understanding the main results. 

We now provide a brief description of each section.

In section \ref{nfoldsuspension}
we introduce the notion of the symmetrisation of an $n$-operad. This is the only section where we
seriously refer to the notion of monoidal globular category from \cite{BatN}. Nevertheless, we hope that the main
notion of  symmetrisation will be clear even without understanding all the details of the definition of
$n$-operad in a general monoidal globular category because  Proposition \ref{lambda} shows that the problem 
of finding a symmetrisation of an $n$-operad $A$ can be reduced to the case where $A$ is of a  special form,
which we call $(n-1)$-terminal. The latter is roughly speaking an operad which has strict $(n-1)$-categories 
as the algebras for its $(n-1)$-skeleton. The reader, therefore,  can start to read our paper from section \ref{symgroup}.

In section \ref{symgroup} we fix our terminology concerning symmetric operads  and obtain a useful combinatorial formula  needed later. 

In Section \ref{treesandmorphisms} we recall the definition of the
$\omega$-category of trees and of the category $\Omega_n$ \cite{BatN,BS,J} which is an $n$-dimensional
analogue of the   category $\Delta_{alg}$ of all finite ordinals and  plays an important role here. More generally, we
believe that the categories
$\Omega_n$ must be  one of the central objects of study in higher dimensional category theory, at least on the 
combinatorial side of the theory.  It appears that $\Omega_n$ contains all the information on the coherence
laws available in weak $n$-categories.

In Section \ref{noperads}
 we give a definition of $n$-operad in a symmetric monoidal category $V$, which is just
an $n$-operad in the monoidal globular category $\Sigma^n V$. This definition is much simpler than the
definition of general $n$-operad and is reminiscent of the classical definition of nonsymmetric operad. 

Section \ref{dessym} is devoted to the construction of a desymmetrisation functor $Des_n$ from symmetric operads to $n$-operads which 
incorporates the action of the symmetric groups.  We also show that the desymmetrisation functor does not alter the endomorphism operads. 
 Here we again refer to our paper \cite{BatN} for a construction of the endomorphism
$n$-operad. However, the reader, can accept our construction here as a definition of endomorphism $n$-operad, so
 again does not need to understand the  technical construction from \cite{BatN}. Our main activity for the
rest of the paper will be an explicit  construction of the symmetrisation functor $Sym_n$  left adjoint to $Des_n$.

In section \ref{InternalAlgebras} we  develop   a   general $2$-categorical method, which in the next sections will allow us to express the Eckmann-Hilton style arguments in terms of left Kan extensions in an appropriate $2$-category.   These techniques will be very useful in the sequel of this paper \cite{BFM}.

In Section \ref{Intop} we reap the first fruits of the theory developed in  Section  \ref{InternalAlgebras} by applying it to $n$-operads and symmetric operads. The results of this section   show that the symmetrisation functor $Sym_n$ exists.

In Section \ref{internaloperads} we consider internal $n$-operads inside categorical  symmetric operads and categorical $n$-operads and  prove that these theories can be represented by some categorical operads   
$\h^n$ and $\H^n .$ We provide   unpacked  definitions  of internal symmetric operads and internal $n$-operads and give some examples.

We continue to study internal operads in Section \ref{H^nh^n} and describe the  operad $\h^n$  in terms of generators and relations. 
We  show that our  theorem \ref{rep} is equivalent to the classical tree formalism for nonsymmetric and symmetric operads if $n=1$ or $n=\infty$ respectively \cite{MSS}. 
 
 In section \ref{IMon} we consider an  example of a categorical symmetric operad containing an internal operad, namely, the operad of $n$-fold monoidal categories of \cite{F}.  This example will be an important ingredient in  one of the  proofs  of a theorem which will relate our categorical constructions to the theory of $n$-fold loop spaces \cite{BFM}.  

Section \ref{freeoperad} has a technical character. We establish a useful formula for the free $n$-operad functor using the techniques developed in Section \ref{InternalAlgebras}.

Finally, in Section \ref{suspension}  we provide our symmetrisation  formula for  the $(n-1)$-terminal 
$n$-operad $A$ in a cocomplete symmetric   monoidal category $V$. The formula is
$$Sym_n(A)_k  \ \simeq \ \mbox{\rm co}\!\lim\limits_{\h^{n}_k} \tilde{A}_k$$
where $\tilde{A}$ is an operadic functor on $\h^n$ which appears from the universal property of $\h^n$.

We also show that in one important  case the symmetrisation functor commutes with the nerve functor, namely $$N(\h^n) \simeq Sym_n(N(\H^n)) .$$
Results like this will play an  important role in the homotopy theory of $n$-operads which we develop in the second part of this paper \cite{BFM}.  
We also  will connect  our symmetrisation formula with the geometry of the Fulton-Macpherson operad \cite{K} and coherence laws for $n$-fold loop spaces in  \cite{BFM} .

\

\noindent {\bf Acknowledgements.}  I would like to thank Ross Street for his constant encouragement and
stimulating discussions during my work on this paper. His intellectual and human contribution to this work cannot
be overestimated. I am   grateful to Alexei Davydov from the discussions with whom I formed many ideas which
constituted  the foundation of this work. I   wish to express my  gratitude to
C.Berger, S.Crans, J.Dolan, E.Getzler,  A.Joyal, S.Lack, M.Markl, P.May, J.McLure, J.Stasheff, D.Tamarkin and  M.Weber, each of whom
provided me with  an important detail of the mosaic  which finally led to the whole picture presented in this
paper. I am also grateful to my anonymous referee whose advice helped to    improve  considerably the presentation of this paper.

 Finally, I gratefully acknowledge the financial support of the Scott Russell Johnson Memorial  Foundation, Macquarie
University Research Commitee and Australian Research Council (grant \# DP0558372).

\section{General symmetrisation problem}\label{nfoldsuspension}

We introduce here the general notion  of symmetrisation  of an $n$-operad in an
augmented monoidal $n$-globular category. 

Let $M$ be an augmented monoidal $n$-globular category \cite{BatN}. Recall  that  part of the structure on $M$ are functors 
$s_k , t_k , z $ which make $M$ a reflexive graph in $Cat .$

  Let
$I$ be the unit object of $M_0$. Fix an integer $k>0$. Then we can
construct the following augmented monoidal $n$-globular category
$M^{(k)}.$ The category $M^{(k)}_l $ is the terminal category when $l< k$.
If $l\ge k$ then $M^{(k)}_l$ is the full subcategory of $M_l$ consisting of
objects $x$ with $$s_{k-1}x = t_{k-1}x = z^{k-1}I.$$ 
There is an obvious inclusion  
$$j:M^{(k)}\rightarrow M.$$
We also can form an augmented monoidal $(n-k)$-globular category
$\Sigma^{-k}M^{(k)}$ with 
$$(\Sigma^{-k}M^{(k)})_l = M^{(k)}_{({l+k})}$$
and obvious augmented monoidal $(n-k)$-globular structure.

Recall \cite{BatN} that a globular object of $M$ is a globular functor
from the terminal $n$-globular category $1$ to $M$. We will call a globular object
$$x:1\rightarrow M$$ {\it $(k-1)$-terminal} if $x$ can be factorised through
$j$. Analogously, {\it a morphism between two  $(k-1)$-terminal globular objects} is a natural transformation  which can
be factorised through
$j$. 

 Let us denote by $gl_n(M)$ and $gl_n^{(k)}(M)$  the categories of globular objects in 
$M$ and $(k-1)$-terminal globular objects in $M$ respectively. Then we have isomorphisms of categories
 $$ gl_n^{(k)}(M) \simeq gl_n(M^{(k)}) \simeq gl_{n-k}(\Sigma^{-k}M^{(k)}).$$

In the same way we can define  {\it $(k-1)$-terminal collections} in
$M$         
\cite{BatN} and {\it $(k-1)$-terminal $n$-operads} in $M$. Again the category of
$(k-1)$-terminal $n$-operads in $M$ is isomorphic to the category of $n$-operads
in
$M^{(k)}$ but is different from the category of $(n-k)$-operads in $\Sigma^{-k}M^{(k)}$.

Suppose now $A$ is an $n$-operad  in $M$ and colimits in $M$ commute
with the augmented monoidal structure \cite{BatN}. Then $A$ generates a
monad  
$\A$   on the category of $n$-globular objects $gl_n(M)$. The algebras of $A$ are, by definition, the algebras of the
monad $\A$.

More generally, let $\A$ be an arbitrary  monad on $gl_n(M)$.
 An algebra
$x$ of
$\A$ is called {\it $(k-1)$-terminal} provided its underlying globular object is
$(k-1)$-terminal. A morphism of $(k-1)$-terminal algebras is a morphism of underlying $(k-1)$-terminal 
objects which is also a morphism of $\A$-algebras. 

 Now let \
$Alg^{(k)}_{\A}$  \ be the category of $(k-1)$-terminal algebras of $\A$. 
We have a forgetful functor
$$U^{(k)}:Alg^{(k)}_{\A} \ \longrightarrow \ gl_n(M^{(k)})\simeq gl_{n-k}(\Sigma^{-k}M^{(k)}) .$$
\begin{defin}
If $U^{(k)}$ is monadic then we  call the corresponding monad on $ gl_{n-k}(\Sigma^{-k}M^{(k)})$ the
$k$-fold suspension of $\A$. \end{defin}
In the special  case $M = Span(Set)$ this definition was given by M.Weber
in his PhD thesis \cite{W}. He also proved that in this case the
suspension  exists for a large class of monads on globular sets. Observe
that  $gl_{\infty}(Span(Set)^{(k)})$ is equivalent to the
category of globular sets again. 

Suppose now that $\A$ is obtained from an $n$-operad $A$ in $M$.
Even if the $k$-fold suspension of $\A$ exists it
is often not true that  the suspension comes from an operad in $ \Sigma^{-k}M^{(k)}$. 
 To handle this situation we need a more general notion
of operad which is not available at this time.  M.Weber has a notion of
symmetric globular operad in the special case $M=Span(Set)$ which seems to be
a good candidate in this situation \cite{W}.  

However, there is one case where such a notion already exists. Indeed, if 
$k = n$ the globular category  $M^{(n)}$ has only one nontrivial 
category  $M^{(n)}_n =  \Sigma^{-n}M^{(n)} = V$. This category has to be braided
monoidal if $n=1$ and symmetric monoidal if $n>1 ;$ but we assume  that  $V$ is symmetric monoidal even if $n=1 .$ 

 The $n$-fold suspension of a monad
$\A$ on $gl_n(M)$ generates, therefore,  a monad on $V$. It
is now natural to ask whether this monad comes from a symmetric 
 operad in $V$.

\begin{defin}\label{symproblem} Let $A$ be an $n$-operad in $M$ such that
the $n$-fold suspension of $\A$ exists and comes from a symmetric
 operad $B$ on $V$. Then we call $B$ the 
symmetrisation  of $A$. The notation is $B= Sym_n(A) .$ \end{defin} 

\Remark  If $n=1$ and $V$ is a braided monoidal category we can give a similar definition with $B$ being a nonsymmetric operad in $V .$  
If the braiding in $V$ is actually a symmetry  we can show that $Sym_1(A)$ is a symmetrisation of the nonsymmetric operad $B$ in the classical sense. 

\

Now we will show that the problem of finding  a symmetrisation of
an $n$-operad  in  $M^n$ can often be subdivided into two steps.

Let $$t: gl_n(M^{(n)})\rightarrow gl_n(M^{})$$
be the  natural inclusion functor. Let $\tau$ be the other obvious inclusion 
$$\tau: O_{n}(M^{(n)})\rightarrow O_n(M^n)$$
where $O_n(C)$ means the category of $n$-operads in $C$.
Let $x$ be a  globular object of $M^{(n)}$. And suppose there
exist  endooperads $End(x)$ and $End(t(x))$ in $M^{(n)}$ and $M$ respectively \cite{BatN}. Then it is
not hard to check that $\tau(End(x)) \simeq End(t(x)).$

If now $x$ is an algebra of some $n$-operad $A$ in $M$ then we have an
operadic  morphism 
$$k:A\rightarrow End(t(x))\simeq \tau(End(x)) $$
Suppose in addition that $\tau$ has a left adjoint $\lambda$. This is
 true in the most interesting cases. Then we have that $k$ is uniquely
determined and determines an operadic map 
$$k':\lambda(A)\rightarrow End(x) .$$
Thus we have 
\begin{pro}\label{lambda} The category of $(n-1)$-terminal algebras of an $n$-operad $A$ in $M$ 
is isomorphic to the
category of algebras of the $n$-operad  $\lambda(A)$ in $M^{(n)}$. \end{pro}

Therefore, to define a symmetrisation of an $n$-operad $A$ we first  find an $(n-1)$-terminal $n$-operad $\lambda(A)$ and then calculate the symmetrisation of $\lambda(A) .$ 

If $V$ is a symmetric monoidal category, we can form the
augmented  monoidal $n$-globular category $L=\Sigma^n V$ where $L$ has $V$
in dimension $n$ and terminal categories in other dimensions. The
monoidal structure is given by $\otimes_i = \otimes$ where $\otimes$ is
tensor product in $V$.  For example,  $M^{(n)} = \Sigma^n(M^{(n)}_n)$.
 In the rest of the paper we will study the case $M=\Sigma^n V$.
We will show that many interesting phenomena appear already in this situation. The passage
from $A$ to $\lambda(A)$ will be studied elsewhere by a method similar to the method developed
in this paper.

\section{Symmetric operads}\label{symgroup}

 For a natural number $n$ we will denote by $[n]$ the ordinal
$$1 \ < \ 2 \ < \ \ldots \ < n .$$
In particular $[0]$ will denote the empty ordinal. Notice, that our notation is not classical. We  find it, however, more convenient for this  exposition.   

 A morphism from $[n]\rightarrow [k]$ 
is any function between underlying sets. It can be order preserving or not. 
It is clear that we then have a category. We denote this category by $\Omega^s$. Of course,
$\Omega^s$ is equivalent to the category of finite sets.    In particular, the symmetric group
$\Sigma_n$ is the group of automorphisms of $[n]$.  

Let $\sigma:[n]\rightarrow [k]$ be a morphism in $\Omega^s$ and 
let  $ 1\le i \le k .$ Then the preimage $\sigma^{-1}(i)$ has a   linear order induced from $[n] .$ Hence, there exists a unique object   $[n_i]\in \Omega^s$ and a unique order preserving bijection $[n_i]\rightarrow \sigma^{-1}(i) .$ We will call $[n_i]$ 
the {\it fiber} of $\sigma$ over $i$ and will denote it $\sigma^{-1}(i) $ or $[n_i].$   

Analogously, given a composite of morphisms  in $\Omega^s :$ 
\begin{equation}\label{nlk} [n]\stackrel{\sigma}{\longrightarrow} [l] \stackrel{\omega}{\longrightarrow} [k] \end{equation}
we will denote
 $\sigma_i$   the {\it $i$-th fiber} of $\sigma$; i.e. the
pullback 

{\unitlength=0.9mm

\begin{picture}(60,40)(-18,-3)

\put(10,25){\makebox(0,0){\mbox{$ \sigma^{-1}(\omega^{-1}(i))$}}}
\put(10,20){\vector(0,-1){10}}
\put(12,15){\shortstack{\mbox{$ $}}}

\put(22,25){\vector(1,0){10}}

\put(24,26){\shortstack{\mbox{$\sigma_i $}}}

\put(45,25){\makebox(0,0){\mbox{$ \omega^{-1}(i)$}}}
\put(45,20){\vector(0,-1){10}}

\put(75,25){\makebox(0,0){\mbox{$[1] $}}}
\put(75,20){\vector(0,-1){10}}

\put(57,25){\vector(1,0){10}}

\put(57,21){\shortstack{\mbox{$ $}}}

\put(57,28){\shortstack{\mbox{\small $ $}}}

\put(77,14){\shortstack{\mbox{\small $\xi_i $}}}

\put(10,5){\makebox(0,0){\mbox{$ [n]$}}}

\put(23,5){\vector(1,0){10}}

\put(27,6){\shortstack{\mbox{$\sigma $}}}

\put(45,5){\makebox(0,0){\mbox{$ [l]$}}}

\put(75,5){\makebox(0,0){\mbox{$[k] $}}}

\put(57,5){\vector(1,0){10}}
\put(80,4){\makebox(0,0){\mbox{$.$}}}

\put(60,6){\shortstack{\mbox{$\omega $}}}

\put(57,8){\shortstack{\mbox{\small $ $}}}

\end{picture}}

The following  is a slightly more functorial version of a classical definition of a symmetric operad \cite{May}. 

 Let $P $ be the subcategory of bijections in  $\Omega^s .$   A {\it right symmetric collection}  in a symmetric monoidal category $V$ is a functor $A:P^{op}\rightarrow V .$ The value of $A$ on an object $[n]$ will be denoted $A_n .$ 
\begin{defin}\label{defsymop} A right symmetric  operad in $V$ is a right symmetric collection $A$ equipped with the following additional structure:

- a morphism  $e:I\rightarrow A_1$

- for every order preserving  map $\sigma:[n]\rightarrow [k]$ in $\Omega^s$   a morphism
: 
$$\mu_{\sigma}: A_{k}\otimes(A_{n_1}\otimes ... \otimes
A_{n_k})\longrightarrow A_{n},
$$
where $[n_i] = \sigma^{-1}(i).$ 

They must satisfy the following identities:

- for any composite of order preserving morphisms in $\Omega^s$ $$[n]\stackrel{\sigma}{\longrightarrow} [l] \stackrel{\omega}{\longrightarrow} [k] ,$$
the following diagram commutes

{\unitlength=1mm

\begin{picture}(300,45)(2,0)

\put(20,35){\makebox(0,0){\mbox{$\scriptstyle A_k\otimes
A_{l_{\bullet}}\otimes A_{n_1^{\bullet}} \otimes  ...
\otimes 
 A_{n_i^{\bullet}}\otimes  ... \otimes A_{n_k^{\bullet}}   
$}}}
\put(20,31){\vector(0,-1){12}}

\put(94,31){\vector(0,-1){12}}

\put(88,35){\makebox(0,0){\mbox{$\scriptstyle A_k\otimes
A_{l_{1}}\otimes A_{n_1^{\bullet}} \otimes  ...
\otimes A_{l_{i}}\otimes
 A_{n_i^{\bullet}}\otimes  ... \otimes A_{l_{k}}\otimes
A_{n_k^{\bullet}}   
$ }}}

\put(50,37){\makebox(0,0){\mbox{$\scriptstyle \simeq $}}}
\put(45,35.5){\vector(1,0){11}}

\put(20,15){\makebox(0,0){\mbox{$\scriptstyle A_l\otimes 
A_{n_1^{\bullet}} \otimes  ...
\otimes 
 A_{n_i^{\bullet}}\otimes  ... \otimes A_{T_k^{\bullet}}
$}}}

\put(94,15){\makebox(0,0){\mbox{$\scriptstyle A_k\otimes 
A_{n_{\bullet}} 
$}}}

\put(60,5){\makebox(0,0){\mbox{$ \scriptstyle A_n 
$}}}

\put(35,11){\vector(4,-1){19}}

\put(85,11){\vector(-4,-1){19}}

\end{picture}}

\noindent 
Here $$A_{l_{\bullet}}= A_{l_1}\otimes ...
\otimes A_{l_k},$$  
$$A_{n_{i}^{\bullet}} = A_{n_i^1} \otimes ...\otimes A_{n_i^{m_i}}$$
and $$ A_{n_{\bullet} } =  A_{n_1}\otimes ...
\otimes A_{n_k};$$

- for an identity $\sigma = id : [n]\rightarrow [n]$ the diagram

{\unitlength=1mm
\begin{picture}(50,25)(30,2)

\put(97,20){\vector(-1,0){20}}

\put(60,17){\vector(0,-1){8}}

\put(60,20){\makebox(0,0){\mbox{\small$A_n\otimes 
A_{1}\otimes ... \otimes A_{1} 
$}}}

\put(114,20){\makebox(0,0){\mbox{\small$A_n\otimes 
{I}\otimes ... \otimes {I} 
$}}}

\put(60,5){\makebox(0,0){\mbox{\small$A_n 
$}}}

\put(105,15){\vector(-4,-1){30}}

\put(90,9){\makebox(0,0){\mbox{\small$id
$}}}

\end{picture}}

\noindent commutes;

- for the unique morphism $[n]\rightarrow [1]$ the diagram

{\unitlength=1mm
\begin{picture}(50,25)(30,2)

\put(87,20){\vector(-1,0){15}}

\put(60,17){\vector(0,-1){8}}

\put(60,20){\makebox(0,0){\mbox{\small$A_{1}\otimes 
A_n
$}}}

\put(98,20){\makebox(0,0){\mbox{\small$I \otimes
A_T
$}}}

\put(60,5){\makebox(0,0){\mbox{\small$A_n 
$}}}

\put(95,17){\vector(-3,-1){25}}

\put(84,11){\makebox(0,0){\mbox{\small$id
$}}}

\end{picture}}

\noindent commutes.

In addition the following 
  two equivariancy conditions must be satisfied:

\begin{enumerate}
\item

 For every commutative diagram in $\Omega^s $

{\unitlength=1mm

\begin{picture}(40,28)(-29,2)

\put(13,25){\makebox(0,0){\mbox{$ [n']$}}}
\put(13,21){\vector(0,-1){10}}
\put(10,15){\shortstack{\mbox{$\pi $}}}
\put(42,15){\shortstack{\mbox{$\rho $}}}

\put(22,25){\vector(1,0){10}}

\put(26,26){\shortstack{\mbox{$\sigma' $}}}

\put(41,25){\makebox(0,0){\mbox{$ [k']$}}}
\put(41,21){\vector(0,-1){10}}

\put(13,7){\makebox(0,0){\mbox{$ [n]$}}}

\put(23,7){\vector(1,0){10}}

\put(27,8){\shortstack{\mbox{$\sigma $}}}

\put(41,7){\makebox(0,0){\mbox{$[k]$}}}

\end{picture}}

whose vertical maps are bijections and whose horizontal
maps are order preserving the following diagram commutes:

{\unitlength=1mm

\begin{picture}(200,33)(-30,0)

\put(9,25){\makebox(0,0){\mbox{$A_{k'}\otimes(A_{n'_{\rho(1)}}
\otimes ... \otimes
A_{n'_{\rho(k)}}) $}}}
\put(10,10){\vector(0,1){10}}
\put(-5,15){\shortstack{\mbox{$\scriptstyle A(\rho)\otimes\tau(\rho) $}}}

\put(32,25){\vector(1,0){12}}

\put(35,27){\shortstack{\mbox{$\mu_{\sigma'} $}}}

\put(50,25){\makebox(0,0){\mbox{$ A_{n'}$}}}
\put(50,10){\vector(0,1){10}}

\put(10,5){\makebox(0,0){\mbox{$A_{k}\otimes(A_{n_1}
\otimes ... \otimes
A_{n_k})$}}}

\put(31,5){\vector(1,0){13}}

\put(35,7){\shortstack{\mbox{$\mu_{\sigma} $}}}

\put(50,5){\makebox(0,0){\mbox{$A_{n}$}}}

\put(55,15){\makebox(0,0){\mbox{$\scriptstyle A(\pi) $}}}

\put(60,4){\shortstack{\mbox{$, $}}}

\end{picture}}
\noindent where $\tau(\rho)$ is the symmetry in $V$  which corresponds to permutation $\rho .$ 

\item For every commutative diagram in $\Omega^s$ 

{\unitlength=1mm

\begin{picture}(40,28)(-29,2)

\put(13,25){\makebox(0,0){\mbox{$ [n'']$}}}
\put(13,21){\vector(0,-1){10}}
\put(10,15){\shortstack{\mbox{$\sigma $}}}
\put(42,15){\shortstack{\mbox{$\eta' $}}}

\put(22,25){\vector(1,0){10}}

\put(26,26){\shortstack{\mbox{$\sigma' $}}}

\put(41,25){\makebox(0,0){\mbox{$ [n']$}}}
\put(41,21){\vector(0,-1){10}}

\put(13,7){\makebox(0,0){\mbox{$ [n]$}}}

\put(23,7){\vector(1,0){10}}

\put(27,8.5){\shortstack{\mbox{$\eta $}}}

\put(41,7){\makebox(0,0){\mbox{$[k]$}}}

\end{picture}}

\noindent where $\sigma,\sigma'$ are bijections and $\eta,\eta'$
are order preserving maps, the following diagram commutes 

{\unitlength=1mm

\begin{picture}(200,53)(-30,0)

\put(10,25){\makebox(0,0){\mbox{$A_{k}\otimes (A_{n''_1}
\otimes ... \otimes
A_{n''_k}) $}}}
\put(10,8){\vector(0,1){12}}
\put(-18,13){\shortstack{\mbox{$\scriptstyle 1\otimes A(\sigma_1)\otimes\ldots\otimes A(\sigma_k) $}}}

\put(-18,34){\shortstack{\mbox{$\scriptstyle 1\otimes A(\sigma'_1)\otimes\ldots\otimes A(\sigma'_k) $}}}

\put(46,34){\shortstack{\mbox{$\scriptstyle  A(\sigma') $}}}
\put(46,13){\shortstack{\mbox{$\scriptstyle  A(\sigma) $}}}
\put(35,47){\makebox(0,0){\mbox{$ \mu_{\eta'}$}}}
\put(34,7){\makebox(0,0){\mbox{$ \mu_{\eta}$}}}

\put(45,25){\makebox(0,0){\mbox{$ A_{n''}$}}}
\put(45,8){\vector(0,1){12}}

\put(75,25){\makebox(0,0){\mbox{$ $}}}

\put(57,21){\shortstack{\mbox{$ $}}}

\put(10,5){\makebox(0,0){\mbox{$A_{k}\otimes (A_{n_1}
\otimes ... \otimes
A_{n_k})$}}}

\put(29,5){\vector(1,0){10}}

\put(27,6){\shortstack{\mbox{$ $}}}

\put(45,5){\makebox(0,0){\mbox{$A_{n}$}}}

\put(75,5){\makebox(0,0){\mbox{$ $}}}


\put(60,6){\shortstack{\mbox{$ $}}}


\put(10,45){\makebox(0,0){\mbox{$A_{k}\otimes(A_{n'_1}
\otimes ... \otimes
A_{n'_k}) $}}}

\put(45,45){\makebox(0,0){\mbox{$ A_{n'}$}}}

\put(29,45){\vector(1,0){10}}

\put(10,41){\vector(0,-1){12}}
\put(45,41){\vector(0,-1){12}}

\end{picture}}

\end{enumerate}
\end{defin}

Let us denote the category of operads in this sense by $SO_r(V) .$ Analogously, we   can construct the category of {\it left symmetric operads} $SO_l(V) $  by asking    a {\it left symmetric collection} to be a covariant functor on $P$ and inverting the corresponding  arrows in equivariancy diagrams.
Clearly, these two categories of operads are isomorphic.

 We will  define yet another category of symmetric operads $O^s(V) .$ 

\begin{defin} An $S$-operad is a  
  collection of objects $\{A_n\}, [n]\in \Omega^s ,$ equipped with :

 - a morphism  $e:I\rightarrow A_1$

- for every  map $\sigma:[n]\rightarrow [k]$ in $\Omega^s$   a morphism  
$$\mu_{\sigma}: A_{k}\otimes(A_{n_1}\otimes ... \otimes
A_{n_k})\longrightarrow A_{n},
$$
where $[n_i] = \sigma^{-1}(i).$ 
 
This structure must satisfy the associativity axiom from the Definition \ref{defsymop}  {\bf  with respect to all maps in} $\Omega^s $ and
 two other axioms concerning identity and trivial maps in $\Omega^s ,$ but no equivariance condition is imposed on $A .$ 
\end{defin}

\begin{pro} The categories   $SO_r(V) , SO_l(V)$ and $O^s(V)$ are isomorphic. \end{pro}

{\noindent Proof} We will construct a functor $$S :O^s(V)\rightarrow SO_r(V)$$ first.  Let $A$ be an object of $O^s(V) .$ 
We construct a symmetric collection $S(A)_n = A_n .$ 
Now we have to define the action of the symmetric groups on $S(A)$. Let 
$\sigma:[n]\rightarrow [n]$ be a permutation. Then  the composite 
$$S(A)_n = A_n \longrightarrow A_n \otimes I \otimes ... \otimes I {\longrightarrow} A_n\otimes A_1\otimes
...\otimes A_1
\stackrel{\mu_{\sigma}}{\longrightarrow}  A_n = S(A)_n$$
determines an endomorphism  $S(A)(\sigma)$. The reader may check as an exercise that $S(A)$ is a contravariant functor on $P .$  

The effect of an order preserving map on $S(A)$ is determined by the effect of this map on $A .$ The equivariance conditions follows easily from these definitions. 

Let us construct an inverse functor
$$(-)^s :SO_r(V)\rightarrow O^s(V)$$ 
On the level of collections, $(A)^s_n = A_n .$
To define $(A)^s$ on an arbitrary map from $\Omega^s$ we recall the following combinatorial fact. 
 Every morphism  $\sigma: [n]\rightarrow 
[k]$ in $\Omega^s$   has a unique factorisation

{\unitlength=1mm

\begin{picture}(200,25)(-17,7)

\put(32,25){\vector(1,0){13}}

\put(37,26){\shortstack{\mbox{$ \sigma $}}}

\put(27,25){\makebox(0,0){\mbox{$[n]$}}}
\put(30,22){\vector(1,-1){6}}
\put(26,16){\shortstack{\mbox{$ \pi(\sigma) $}}}

\put(75,25){\makebox(0,0){\mbox{$ $}}}

\put(50,25){\makebox(0,0){\mbox{$[k]$}}}

\put(41,16){\vector(1,1){6}}

\put(45,16){\shortstack{\mbox{$ \nu(\sigma) $}}}

\put(39,13){\makebox(0,0){\mbox{$[n']$}}}

\put(75,5){\makebox(0,0){\mbox{$ $}}}

\put(60,6){\shortstack{\mbox{$ $}}}

\end{picture}}

\noindent where
$\nu(\sigma)$ is order preserving, while $\pi(\sigma)$ is bijective and preserves order 
on the fibers of $\sigma_n$. We use this factorisation to define 
the effect on $\sigma$ of $(A)^s$ by requiring the  commutativity of the following diagram

{\unitlength=1mm

\begin{picture}(200,33)(-30,0)

\put(10,25){\makebox(0,0){\mbox{$A_{k}\otimes(A_{n_1}
\otimes ... \otimes
A_{n_k}) $}}}
\put(10,20){\vector(0,-1){10}}

\put(30,25){\vector(1,0){13}}

\put(33,26){\shortstack{\mbox{$\scriptstyle  A^s({\sigma}) $}}}

\put(48,25){\makebox(0,0){\mbox{$ A_n$}}}
\put(48,10){\vector(0,1){10}}
\put(50,15){\shortstack{\mbox{$\scriptstyle  A(\pi(\sigma)) $}}}

\put(75,25){\makebox(0,0){\mbox{$ $}}}

\put(57,21){\shortstack{\mbox{$ $}}}

\put(10,5){\makebox(0,0){\mbox{$ A_{k}\otimes(A_{n'_1}
\otimes ... \otimes
A_{n'_k})$}}}

\put(30,5){\vector(1,0){13}}

\put(31,6){\shortstack{\mbox{$\scriptstyle A(\nu(\sigma)) $}}}

\put(48,5){\makebox(0,0){\mbox{$A_{n'}$}}}

\put(75,5){\makebox(0,0){\mbox{$ $}}}

\put(60,6){\shortstack{\mbox{$ $}}}

\end{picture}}

\noindent where actually $[n'] = [n] ,$ $[n'_i] = [n_i]$ and the left vertical map is the identity since $\pi(\omega)$ is the identity on the fibers of $\sigma$ and $\nu(\sigma)$.  

Now consider the composite (\ref{nlk})
 of morphisms in $\Omega^s.$ It induces
the following factorisation diagram 

{\unitlength=1mm

\begin{picture}(180,60)(-10,-6)

\put(20,25){\makebox(0,0){\mbox{$[n]$}}}

\put(24,25){\vector(1,0){15}}
\put(24,30){\vector(1,1){15}}
\put(48,45){\vector(1,-1){15}}

\put(30,26){\shortstack{\small \mbox{$ \sigma $}}}
\put(53,26){\shortstack{\small \mbox{$ \omega $}}}

\put(18,36){\shortstack{\small \mbox{$\pi( \sigma\cdot
\omega)
$}}}

\put(58,36){\shortstack{\small \mbox{$\nu( \sigma\cdot
\omega)
$}}}

\put(61,16){\shortstack{\small \mbox{$\nu(
\omega)
$}}}

\put(18,16){\shortstack{\small \mbox{$\pi(
\sigma)
$}}}

\put(15,4){\shortstack{\small \mbox{$\pi(\nu(\sigma)\cdot
\pi(\omega))
$}}}
\put(49,4){\shortstack{\small \mbox{$\nu(\nu(\sigma)\cdot
\pi(\omega))
$}}}

\put(29,18){\shortstack{\small \mbox{$\nu(
\sigma)
$}}}

\put(49.5,18){\shortstack{\small \mbox{$\pi(
\omega)
$}}}

\put(24,20){\vector(1,-1){5}}
\put(47,20){\vector(1,-1){5}}
\put(35,7.5){\vector(1,-1){5}}

\put(45,3){\vector(1,1){5}}
\put(57,15){\vector(1,1){5}}
\put(34,16){\vector(1,1){5}}

\put(43,25){\makebox(0,0){\mbox{$[l]$}}}
\put(43,48){\makebox(0,0){\mbox{$[n''']$}}}
\put(43,0){\makebox(0,0){\mbox{$[n'']$}}}

\put(31.5,12.5){\makebox(0,0){\mbox{$[n']$}}}
\put(54.5,12.5){\makebox(0,0){\mbox{$[l']$}}}

\put(47,25){\vector(1,0){15}}

\put(66,25){\makebox(0,0){\mbox{$[k]$}}}

\end{picture}}

\noindent which in its turn generates the following huge diagram.

{\unitlength=1mm

\begin{picture}(300,130)(5,0)

\put(40,130){\makebox(0,0){\small\mbox{$ A_k
A_{l_{\bullet}} A_{n_1^{\bullet}} 
\ldots A_{n_k^{\bullet}}   
$}}}
\put(30,125){\vector(-1,-1){8}}

\put(56,130){\vector(1,0){8}}

\put(93,125){\vector(1,-1){8}}

\put(84,130){\makebox(0,0){\small\mbox{$ A_k
A_{l_{1}} A_{n_1^{\bullet}} \ldots
 A_{l_{k}}
A_{n_k^{\bullet}}   
$ }}}

\put(52,130){\makebox(0,0){\mbox{$ $}}}


\put(15,111){\makebox(0,0){\mbox{\small$A_k
A_{l'_{\bullet}} A_{n_1^{\bullet}} 
\ldots A_{n_k^{\bullet}}   
$}}}

\put(110,111){\makebox(0,0){\mbox{\small$ A_k
A_{l_{1}} A_{{n'}_1^{\bullet}} \ldots
 A_{l_{k}}
A_{{n'}_k^{\bullet}}   
$}}}

\put(13,105){\vector(-1,-2){7}}

\put(108,105){\vector(1,-2){7}}


\put(6,85){\makebox(0,0){\mbox{\small$
A_{l'} A_{n_1^{\bullet}} 
\ldots A_{n_k^{\bullet}}  
$}}}

\put(120,85){\makebox(0,0){\mbox{\small$A_k
 A_{{n'}_1} \ldots
A_{{n'}_k}   
$}}}

\put(6,58){\makebox(0,0){\mbox{\small$
A_{l} A_{n_1^{\bullet}} 
\ldots A_{n_k^{\bullet}}  
$}}}

\put(4,78){\vector(0,-1){12}}

\put(118,78){\vector(0,-1){12}}

\put(120,58){\makebox(0,0){\mbox{\small$A_k
 A_{{n}_1} \ldots
A_{{n}_k}   
$}}}

\put(16,31){\makebox(0,0){\mbox{\small$
A_{l} A_{{n'}_1^{\bullet}} 
\ldots A_{{n'}_k^{\bullet}}  
$}}}

\put(6,51){\vector(1,-2){7}}

\put(116,51){\vector(-1,-2){7}}

\put(110,31){\makebox(0,0){\mbox{\small$A_k
 A_{{n'''}_1} \ldots
A_{{n'''}_k}   
$}}}



\put(65,35){\makebox(0,0){\mbox{\small$A_{n''} 
$}}}



\put(44,18){\makebox(0,0){\mbox{\small$A_{n'} 
$}}}

\put(84,18){\makebox(0,0){\mbox{\small$A_{n'''} 
$}}}

\put(65,10){\makebox(0,0){\mbox{\small$A_n 
$}}}


\put(48,16){\vector(3,-1){12}}
\put(80,16){\vector(-3,-1){12}}

\put(26,26){\vector(2,-1){12}}
\put(99,26){\vector(-2,-1){12}}

\put(60,32){\vector(-1,-1){10}}
\put(65,30){\vector(0,-1){12}}

\put(65,55){\makebox(0,0){\small\mbox{$associativity
$}}}

\put(40,75){\makebox(0,0){\small\mbox{$ A_k
A_{{l'}_{\bullet}} A_{{n''}_1^{\bullet}} 
\ldots A_{{n''}_k^{\bullet}}   
$}}}
\put(45,70){\vector(-1,-4){5}}

\put(57,75){\vector(1,0){6}}

\put(84,70){\vector(1,-4){5}}

\put(84,75){\makebox(0,0){\small\mbox{$ A_k
A_{{l'}_{1}} A_{{n''}_1^{\bullet}} \ldots
 A_{{l'}_{k}}
A_{{n''}_k^{\bullet}}   
$ }}}


\put(36,45){\makebox(0,0){\mbox{\small$
A_{l'} A_{{n''}_1^{\bullet}} 
\ldots A_{{n''}_k^{\bullet}}  
$}}}

\put(44,42){\vector(2,-1){10}}

\put(85,42){\vector(-2,-1){10}}

\put(90,45){\makebox(0,0){\mbox{\small$A_k
 A_{{n''}_1} \ldots
A_{{n''}_k}   
$}}}

\put(23,66){\makebox(0,0){\mbox{\small$
A_{l'} A_{{n'}_1^{\bullet}} 
\ldots A_{{n'}_k^{\bullet}}  
$}}}
\put(35,42){\vector(-1,-1){8}}

\put(8,78){\vector(1,-1){8}}

\put(20,62){\vector(0,-1){22}}

\put(25,62){\vector(1,-1){12}}

\put(35,98){\makebox(0,0){\mbox{\small$A_k
A_{{l'}_{\bullet}} A_{{n'}_1^{\bullet}} 
\ldots A_{{n'}_k^{\bullet}}   
$}}}

\put(17,105){\vector(3,-1){10}}
\put(28,94){\vector(-1,-3){8}}
\put(36,94){\vector(1,-2){8}}

\put(42,30){\makebox(0,0){\small\mbox{$equivariance \ 1
$}}}

\put(82,30){\makebox(0,0){\small\mbox{$equivariance \ 2
$}}}

\put(100,90){\makebox(0,0){\small\mbox{$equivariance \ 1
$}}}


\put(59,111){\makebox(0,0){\mbox{\small$A_k
A_{l_{\bullet}} A_{{n'}_1^{\bullet}} 
\ldots A_{{n'}_k^{\bullet}}   
$}}}
\put(62,105){\vector(-1,-3){9}}
\put(45,125){\vector(1,-1){8}}
\put(47,105){\vector(-3,-1){10}}
\put(92,105){\vector(-1,-3){9}}
\put(78,111){\vector(1,0){8}}
\put(110,75){\vector(-1,-2){11}}
\put(108,53){\vector(-2,-1){8}}







\end{picture}}

\noindent In
this diagram the central region commutes because of associativity of
$A$ with respect to order preserving maps. Other regions commute
  by one of the equivariance conditions,   by naturality  or  
functoriality. The commutativity of this diagram means associativity
of $A^s$ with respect to arbitrary maps in $\Omega^s .$ 

It is also obvious that the functor $(-)^s$ is inverse to $S(-) .$

\Q

\

Recall  that the symmetric groups  form a symmetric operad $\Sigma$ in $Set$
sometimes called the permutation operad  in the literature. Let us describe this operad explicitly  as a right symmetric operad. 

 The collection $\Sigma_n$ consists of the bijections  from $[n]$ to $[n] .$
 Let $\Gamma$ be   multiplication in $\Sigma .$ 
 One can give the following explicit formula for $\Gamma :$
$$\Gamma(\sigma_k;\sigma_{n_1},\ldots,\sigma_{n_k}) = \Gamma(1_{[k]};\sigma_{n_1},\ldots,\sigma_{n_k})
\cdot \Gamma(\sigma_k;1_{[n_1]}, \ldots , 1_{[n_k]})$$
where $1_{n}$ means the identity bijection of $[n],$ 
$$\Gamma(1_{k};\sigma_{n_1},\ldots,\sigma_{n_k})= \sigma_{n_1}\oplus \ldots \oplus \sigma_{n_k}$$
and 
$$ \Gamma(\sigma;1_{n_1}, \ldots , 1_{n_k})(p)=
 \sum_{\sigma(k)< \sigma(i+1)} n_{k} +p  - \sum_{0\le l\le i}n_l ,$$
when $ n_0+ \ldots + n_i < p \le  n_0+ \ldots + n_{i+1} $ and we assume that $n_0=0 .$ 
In other words $\Gamma(\sigma;1_{n_1}, \ldots , 1_{n_k})$ permutes blocks $[n_1], \ldots , [n_k]$ in accord with the permutation $\sigma .$  

We can illustrate the multiplication $$\Gamma((132);(21),(12),(1))=
\Gamma((132);(21),1_{2},1_{1})$$ by  the following  picture to be read from top to bottom:

{\unitlength=0.7mm

\begin{picture}(60,53)(-28,-6)

\multiput(0,40)(2,0){44}{\line(1,0){1}}
\put(12,15){\shortstack{\mbox{$ $}}}
\multiput(0,20)(2,0){44}{\line(1,0){1}}

\multiput(0,0)(2,0){44}{\line(1,0){1}}
\multiput(4,40)(20,0){5}{\circle*{1}}
\multiput(4,20)(20,0){5}{\circle*{1}}
\multiput(4,0)(20,0){5}{\circle*{1}}
\put(4,40){\line(1,-1){20}}
\put(24,40){\line(-1,-1){20}}
\put(44,40){\line(0,-1){20}}
\put(64,40){\line(0,-1){20}}
\put(84,40){\line(0,-1){20}}

\put(4,20){\line(0,-1){20}}
\put(24,20){\line(0,-1){20}}
\put(44,20){\line(1,-1){20}}
\put(64,20){\line(1,-1){20}}
\put(84,20){\line(-2,-1){40}}

\end{picture}}

 
\begin{lem} \label{pisigma} For the  composite (\ref{nlk}) in $\Omega^s$  the following formula holds
$$ \pi(\sigma\cdot\omega)\cdot\Gamma(1; \pi(\sigma_1),..., \pi(\sigma_k)) 
 = 
\pi(\sigma)\cdot\Gamma(\pi(\omega);1_{\sigma^{-1}(1)},...,1_{\sigma^{-1}(l)})
.$$

\end{lem}

The idea of the lemma is presented in the  diagram

{\epsfxsize=250pt 
\makebox(300,300)[r]{\epsfbox{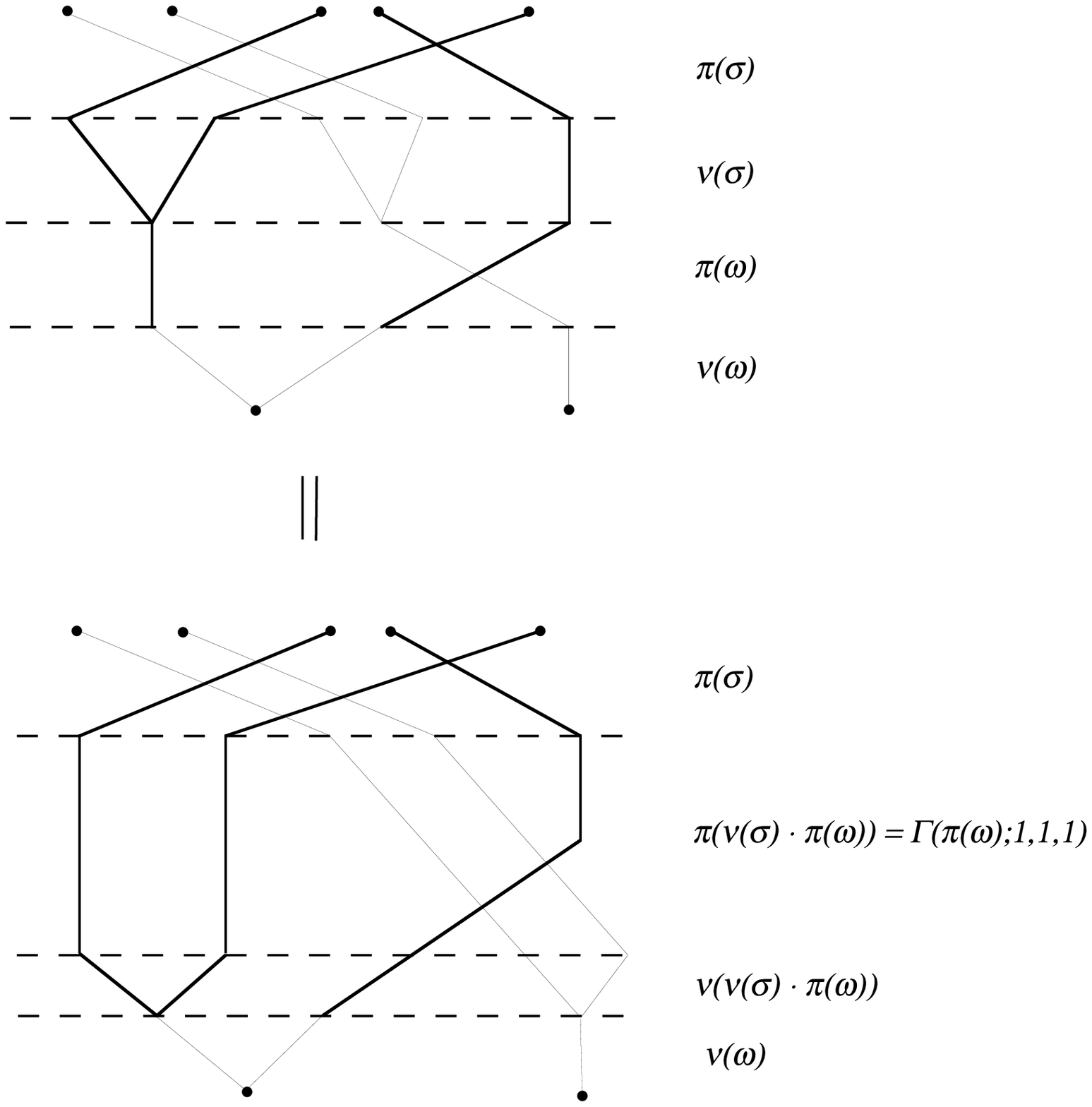}}}

\noindent where  the bottom diagram is just an appropriate deformation of the top one.

\

\Proof For a proof it is sufficient to 
consider the commutative   diagram for associativity generated by (\ref{nlk}) for the  $S$-operad  $\Sigma^s$   and then to calculate both sides of this diagram on identity permutations. 

\Q

\

From now on we accept as agreed that the term {\it symmetric operad} will mean   the {\it left symmetric operad} unless a different understanding   is not required explicitly. The reason for this agreement is practical: many  classical   operads are described as left symmetric operads. Also the description of multiplication in a left symmetric operad is often easier.

\section{Trees and their morphisms}\label{treesandmorphisms}

\begin{defin} A tree of height $n$ (or simply $n$-tree) is a chain of  order preserving maps of
ordinals
$$T=
[k_n]\stackrel{\rho_{n-1}}{\longrightarrow}[k_{n-1}]\stackrel{\rho_{n-2}}{\longrightarrow}...
\stackrel{\rho_0}{\longrightarrow} [1]
$$  
 \end{defin}
If $i\in [k_m]$ and there is no $j\in [k_{m+1}]$ such that $\rho_{m}(j)=i$ then we call $i$ {\it
a leaf of $T$ of height $i$}.   We will call the leaves of $T$
of  height 
$n$ {\it the tips} of 
$T$. If for an $n$-tree $T$ all its leaves are tips we call
such a tree {\it pruned}.

We illustrate the definition in a picture

{\epsfxsize=300pt 
\makebox(300,130)[r]{\epsfbox{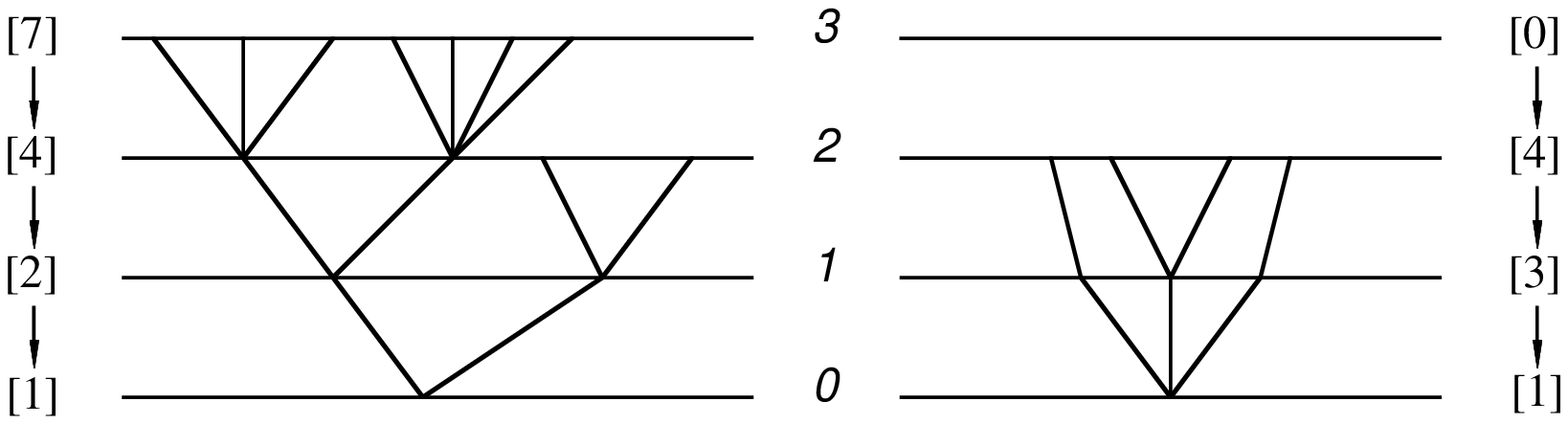}}}

The tree on the right side of the picture has the empty ordinal at the highest level. We will
call such trees {\it degenerate}. There is actually an operation on trees which we denote by
$z(-)$ which assigns to the $n$-tree
$[k_n]\stackrel{}{\rightarrow}[k_{n-1}]\stackrel{}{\rightarrow}...
\stackrel{}{\rightarrow} [1]$ the $(n+1)$-tree
$$[0]\longrightarrow
[k_n]\stackrel{}{\longrightarrow}[k_{n-1}]\stackrel{}{\longrightarrow}...
\stackrel{}{\longrightarrow} [1] .$$

Two other operations on trees are {\it truncation}
$$\partial([k_n]\stackrel{}{\rightarrow}[k_{n-1}]\stackrel{}{\rightarrow}...
\stackrel{}{\rightarrow} [1])=[k_{n-1}]\stackrel{}{\rightarrow}...
\stackrel{}{\rightarrow} [1] $$
and {\it suspension}
$$S([k_n]\stackrel{}{\rightarrow}[k_{n-1}]\stackrel{}{\rightarrow}...
\stackrel{}{\rightarrow} [1])=[k_n]\stackrel{}{\rightarrow}[k_{n-1}]\stackrel{}{\rightarrow}...
\stackrel{}{\rightarrow} [1]\stackrel{}{\rightarrow} [1] . $$

\begin{defin} A tree $T$ is called a $k$-fold suspension if it can be obtained from another tree by
applying  the operation of suspension $k$-times.  The suspension index  
 $susp(T)$ is the  maximum integer $k$ such that $T$ is a $k$-fold
suspension. \end{defin}
The only $n$-tree with suspension index equal to $n$ is the linear tree $$U_n = [1]\rightarrow
\ldots\rightarrow [1].$$

We now define the source and target of a tree $T$ to be equal to $\partial( T)$. So we have a globular
structure on the set of all trees. We actually have more. The  trees form an $\omega$-category
$Tr$ with the set of $n$-cells being equal to the set of the trees of height $n$. 
If two $n$-trees $S$ and $T$ have the same $k$-sources and $k$-targets (i.e.
$\partial^{n-k}T=\partial^{n-k}S$ ) then they can be composed, and the  composite will be
denoted by $S\otimes_k T$. Then  $z(T)$ is  the identity of
the $n$-cell $T$. Here are examples of the $2$-categorical operations on trees

{\epsfxsize=320pt 
\makebox(320,140)[r]{\epsfbox{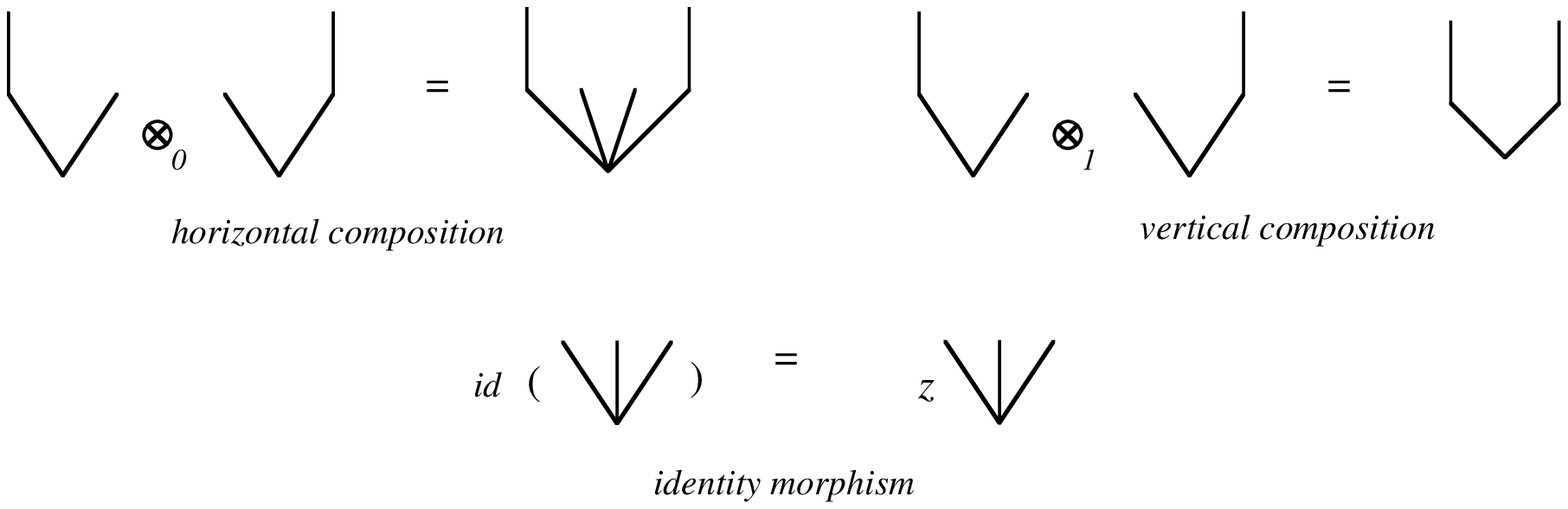}}}

The $\omega$-category $Tr$ is actually the free $\omega$-category generated by the terminal globular
set. Every $n$-tree can be considered as a special sort of $n$-pasting diagram called {\it
globular.} This construction was called the $\star$-construction in \cite{BatN}. Here are a couple
of examples.

 {\epsfxsize=230pt 
\makebox(280,150)[r]{\epsfbox{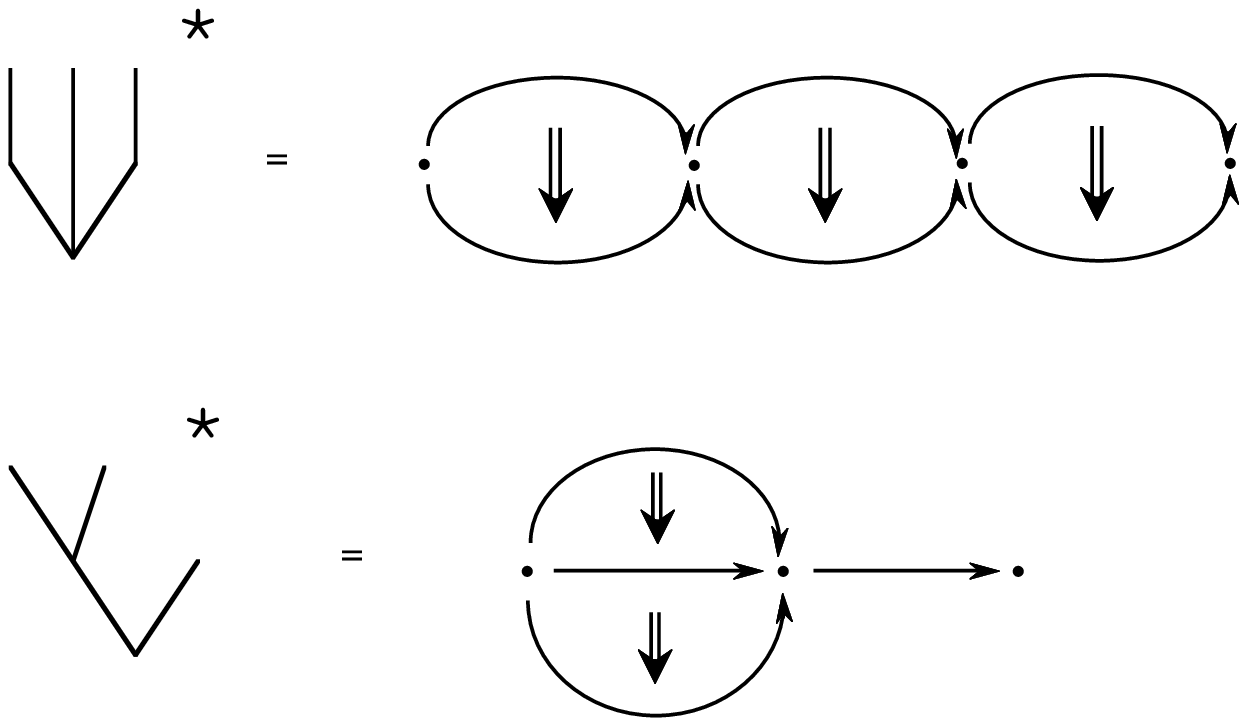}}}

\

 For a globular set $X$ one can then form the set $D(X)$
of all globular pasting diagrams labelled in $X$. This
is  the free $\omega$-category generated by $X$.
In this way we have a monad $(D,\mu,\epsilon)$ on the
category of globular sets, which plays a central role in
our work \cite{BatN}. 

In particular, $D(1)= Tr.$ We also can consider
$D(Tr) = D^2(1)$. It was observed in \cite[p.80-81]{BatN} that the $n$-cells of $D(Tr)$ which were called in \cite[p.80]{BatN} diagrams of $n$-stage trees,  can be identified
with the morphisms of another category of $n$-planar  trees (or the same as  open maps of $n$-disks) introduced by A.Joyal in \cite{J}. 
   This category was called
$\Omega_n$ in \cite[p.10]{BS}. 
It was found that the collection of categories $\Omega_n, n\ge 0,$ forms an $\omega$-category in $Cat$ and, moreover, it
is freely generated by an internal $\omega$-category (called globular monoid in \cite{BS}). So it is a higher dimensional analogue of the 
category $\Delta_{alg}=\Omega_1$ (which is of course the free monoidal category generated by a monoid \cite{ML}). A general theory of such universal objects is developed in   section \ref{InternalAlgebras} of our paper. We also would like  to mention that  C.Berger  also describes   maps in $\Omega_n$ in \cite[1.8-1.9]{BergerC} as dual 
to his   {\it cover} maps of trees.

The definition below is taken from \cite{J} but also presented  in
\cite{BS}[p.11].  
\begin{defin} The category $\Omega_n$ has as objects the trees of height $n$. The morphisms of $\Omega_n$ are
commutative diagrams    

{\unitlength=0.9mm

\begin{picture}(60,35)(-12,0)

\put(10,25){\makebox(0,0){\mbox{$ [k_n]$}}}
\put(10,20){\vector(0,-1){10}}
\put(12,15){\shortstack{\mbox{$ $}}}

\put(17,25){\vector(1,0){10}}
\put(17,5){\vector(1,0){10}}
\put(18,26){\shortstack{\small\mbox{$\rho_{n-1} $}}}
\put(18,6){\shortstack{\small\mbox{$\xi_{n-1} $}}}

\put(35,25){\makebox(0,0){\mbox{$ [k_{n-1}]$}}}
\put(35,20){\vector(0,-1){10}}
\put(35,5){\makebox(0,0){\mbox{$ [s_{n-1}]$}}}

\put(42,25){\vector(1,0){10}}
\put(42,5){\vector(1,0){10}}
\put(43,26){\shortstack{\small\mbox{$\rho_{n-2} $}}}
\put(43,6){\shortstack{\small\mbox{$\xi_{n-2} $}}}

\put(56,25){\shortstack{\mbox{$.  \   .  \  .  $}}}
\put(56,5){\shortstack{\mbox{$ . \ . \ . $}}}

\put(80,25){\makebox(0,0){\mbox{$[1] $}}}
\put(80,20){\vector(0,-1){10}}

\put(66,25){\vector(1,0){10}}
\put(66,5){\vector(1,0){10}}
\put(69,26){\shortstack{\small\mbox{$\rho_0 $}}}
\put(69,6){\shortstack{\small\mbox{$\xi_0 $}}}

\put(57,21){\shortstack{\mbox{$ $}}}

\put(57,28){\shortstack{\mbox{\small $ $}}}

\put(10,5){\makebox(0,0){\mbox{$ [s_n]$}}}

\put(80,5){\makebox(0,0){\mbox{$[1] $}}}

\put(57,8){\shortstack{\mbox{\small $ $}}}

\put(5,14){\shortstack{\mbox{\small $\sigma_n $}}}
\put(26.5,14){\shortstack{\mbox{\small $\sigma_{n-1} $}}}
\put(75,14){\shortstack{\mbox{\small $\sigma_0 $}}}

\end{picture}}

\noindent of sets and functions such that for all $i$ and all $j\in [k_{i-1}]$ the
restriction of
$\sigma_i$ to  $\rho^{-1}_{i-1}(j)$ preserves the natural order on it.
\end{defin}

Let $T$ be an $n$-tree and let $i$ be a leaf of height $m$ of $T$. Then $i$ determines a unique morphism 
$\xi_i: z^{n-m}U_m\rightarrow T$ in
$\Omega_n$ such that $\xi_{i}(1)=i$. We will often identify the leaf with this morphism. 

Let $\sigma:T\rightarrow S$ be a morphism in $\Omega_n$ and let $i$ be a leaf of $T$ . Then 
{\it the fiber  of $\sigma$ over $i$} is the following pullback in $\Omega_n$ 

{\unitlength=0.9mm

\begin{picture}(60,35)(-2,0)

\put(45,25){\makebox(0,0){\mbox{$ \sigma^{-1}(i)$}}}
\put(45,20){\vector(0,-1){12}}

\put(77,25){\makebox(0,0){\mbox{$z^{n-m}U_m $}}}
\put(75,20){\vector(0,-1){12}}

\put(54,25){\vector(1,0){13}}

\put(77,14){\shortstack{\mbox{\small $\xi_i $}}}

\put(45,5){\makebox(0,0){\mbox{$ T$}}}

\put(75,5){\makebox(0,0){\mbox{$S $}}}

\put(54,5){\vector(1,0){14}}

\put(60,6){\shortstack{\small \mbox{$\sigma $}}}

\put(57,8){\shortstack{\mbox{\small $ $}}}

\end{picture}}

\noindent which can be calculated as a levelwise pullback in $Set$.

Now, for a morphism $\sigma:T\rightarrow S$ one can construct a labelling of the
pasting scheme $S^{\star}$  in the $\omega$-category $Tr$ by associating to a
vertex $i$ from $S$ the fiber of $\sigma$ over  $i$. The result of the
pasting operation will be exactly $T$. We will use extensively this
correspondence between morphisms in $\Omega_n$ and pasting schemes in $Tr$.

Some trees will play a special role in our paper.  We will denote by $M_l^j$ the tree
$$\underbrace{U_n\otimes_l \ldots \otimes_l U_n}_{j-\mbox{\small times}} .$$
A picture for $M_l^j$ is as follows
 
{\epsfxsize=70pt 
\makebox(200,100)[r]{\epsfbox{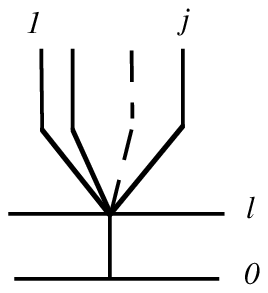}}}

 Now let $T$ be a tree with $susp(T)=l$. Then it is easy to see that we have a unique
representation
$$T=T_1\otimes_l \ldots \otimes_l T_j$$
where $susp(T_i)>l .$ In \cite{BatN} we called this representation the {\it canonical
decomposition} of
$T$. We also will refer to the canonical decomposition when talking about the
 morphism
$$T\longrightarrow M_l^j$$
it generates.

\section{ n-Operads in symmetric monoidal categories}\label{noperads}

It is clear from the definitions of the previous section that 
 the assignment to an $n$-tree 
$$S = [k_n]\stackrel{}{\rightarrow}[k_{n-1}]\stackrel{}{\rightarrow}...
\stackrel{}{\rightarrow} [1]$$
 of its ordinal of tips $[k_n]$ gives us a functor
\begin{equation}\label{[]}[-]:\Omega_n \rightarrow \Omega^s .\end{equation}
We also introduce the notation $|S|$ for the number of tips of the $n$-tree $S .$

Here and in all subsequent sections a  fiber of a morphism $\sigma:T\rightarrow S$ in $ \Omega_n$ 
will mean
only a fiber over  a tip of $S$. So every $\sigma:T\rightarrow S$ with $|S|= k$ determines a list
of trees
$T_1,\ldots,T_k$ being fibers over tips of $S$ ordered according to the order in $[S]$. From now on
we will always relate to $\sigma$ this list of trees in this order.

The definition below is a specialisation  of a general definition
of  $n$-operad in an augmented monoidal $n$-globular category $M$ given in \cite{BatN}. Let $(V,\otimes,I)$ be a
(strict) symmetric monoidal category. Put $M= \Sigma^n V$ , which
means that $M$ has terminal categories up to dimension $n-1$ and
$V$ in dimension $n$. The augmented monoidal structure is given by $\otimes_i = \otimes$ for all
$i$. Then an $n$-operad in $V$ will mean an $n$-operad in $\Sigma^n V$. Explicitly  it means
the following.

\begin{defin}\label{defnoper}  An $n$-operad in $V$ is 
a collection $A_T, \ T\in \Omega_n$, of objects of $V$ equipped with the following structure :

- a morphism $e: I \rightarrow  A_{U_n}$ (the unit);

- for every morphism $\sigma:T \rightarrow S$ in $\Omega_n ,$ 
a morphism 
$$m_{\sigma}:A_S\otimes A_{T_1}\otimes ... \otimes A_{T_k}
 \rightarrow A_T \mbox{\ \ (the multiplication}).$$

They must satisfy the following identities:

- for any composite $T\stackrel{\sigma}{\rightarrow} S \stackrel{\omega}{\rightarrow} R ,$
the associativity diagram

{\unitlength=1mm

\begin{picture}(300,45)(2,0)

\put(20,35){\makebox(0,0){\mbox{$\scriptstyle A_R\otimes
A_{S_{\bullet}}\otimes A_{T_1^{\bullet}} \otimes  ...
\otimes 
 A_{T_i^{\bullet}}\otimes  ... \otimes A_{T_k^{\bullet}}   
$}}}
\put(20,31){\vector(0,-1){12}}

\put(94,31){\vector(0,-1){12}}

\put(88,35){\makebox(0,0){\mbox{$\scriptstyle A_R\otimes
A_{S_{1}}\otimes A_{T_1^{\bullet}} \otimes  ...
\otimes A_{S_{i}}\otimes
 A_{T_i^{\bullet}}\otimes  ... \otimes A_{S_{k}}\otimes
A_{T_k^{\bullet}}   
$ }}}

\put(50,35){\makebox(0,0){\mbox{$\scriptstyle \simeq $}}}

\put(20,15){\makebox(0,0){\mbox{$\scriptstyle A_S\otimes 
A_{T_1^{\bullet}} \otimes  ...
\otimes 
 A_{T_i^{\bullet}}\otimes  ... \otimes A_{T_k^{\bullet}}
$}}}

\put(94,15){\makebox(0,0){\mbox{$\scriptstyle A_R\otimes 
A_{T_{\bullet}} 
$}}}

\put(60,5){\makebox(0,0){\mbox{$ \scriptstyle A_T 
$}}}

\put(35,11){\vector(4,-1){19}}

\put(85,11){\vector(-4,-1){19}}

\end{picture}}

\noindent commutes,
where $$A_{S_{\bullet}}= A_{S_1}\otimes ...
\otimes A_{S_k},$$  
$$A_{T_{i}^{\bullet}} = A_{T_i^1} \otimes ...\otimes A_{T_i^{m_i}}$$
and $$ A_{T_{\bullet} } =  A_{T_1}\otimes ...
\otimes A_{T_k};$$

- for an identity $\sigma = id : T\rightarrow T$ the diagram

{\unitlength=1mm
\begin{picture}(50,25)(30,2)

\put(97,20){\vector(-1,0){20}}

\put(60,17){\vector(0,-1){8}}

\put(60,20){\makebox(0,0){\mbox{\small$A_T\otimes 
A_{U_n}\otimes ... \otimes A_{U_n} 
$}}}

\put(114,20){\makebox(0,0){\mbox{\small$A_T\otimes 
{I}\otimes ... \otimes {I} 
$}}}

\put(60,5){\makebox(0,0){\mbox{\small$A_T 
$}}}

\put(105,15){\vector(-4,-1){30}}

\put(90,9){\makebox(0,0){\mbox{\small$id
$}}}

\end{picture}}

\noindent commutes;

- for the unique morphism $T\rightarrow U_n$ the diagram

{\unitlength=1mm
\begin{picture}(50,25)(30,2)

\put(87,20){\vector(-1,0){15}}

\put(60,17){\vector(0,-1){8}}

\put(60,20){\makebox(0,0){\mbox{\small$A_{U_n}\otimes 
A_T
$}}}

\put(98,20){\makebox(0,0){\mbox{\small$I \otimes
A_T
$}}}

\put(60,5){\makebox(0,0){\mbox{\small$A_T 
$}}}

\put(95,17){\vector(-3,-1){25}}

\put(84,11){\makebox(0,0){\mbox{\small$id
$}}}

\end{picture}}

\noindent commutes.

\end{defin}

The definition of morphism of $n$-operads is the obvious one, so we have a category of $n$-operads in $V$ which
we will denote by $O_n(V)$.

 We give  an example of a $2$-operad to provide the reader with a
feeling  of how these operads  look. Other examples will appear later in the
course of the paper. 

\
 
\Example    One can construct a $2$-operad $G$ in $Cat$ such that the algebras of
$G$ in
$ Cat$  are braided strict monoidal  categories. If we apply the functor $\tau$ to $G$ (i.e.
consider it as a $1$-terminal  operad in
$Span(Cat)$) then the algebras of $\tau(G)$ are Gray-categories \cite{GPS}. The
categories
$G_T$ are  chaotic groupoids with objects  corresponding to so called
$T$-shuffles. The nice geometrical pictures below  show these groupoids in some low
dimensions.

{\epsfxsize=260pt 
\makebox(270,290)[r]{\epsfbox{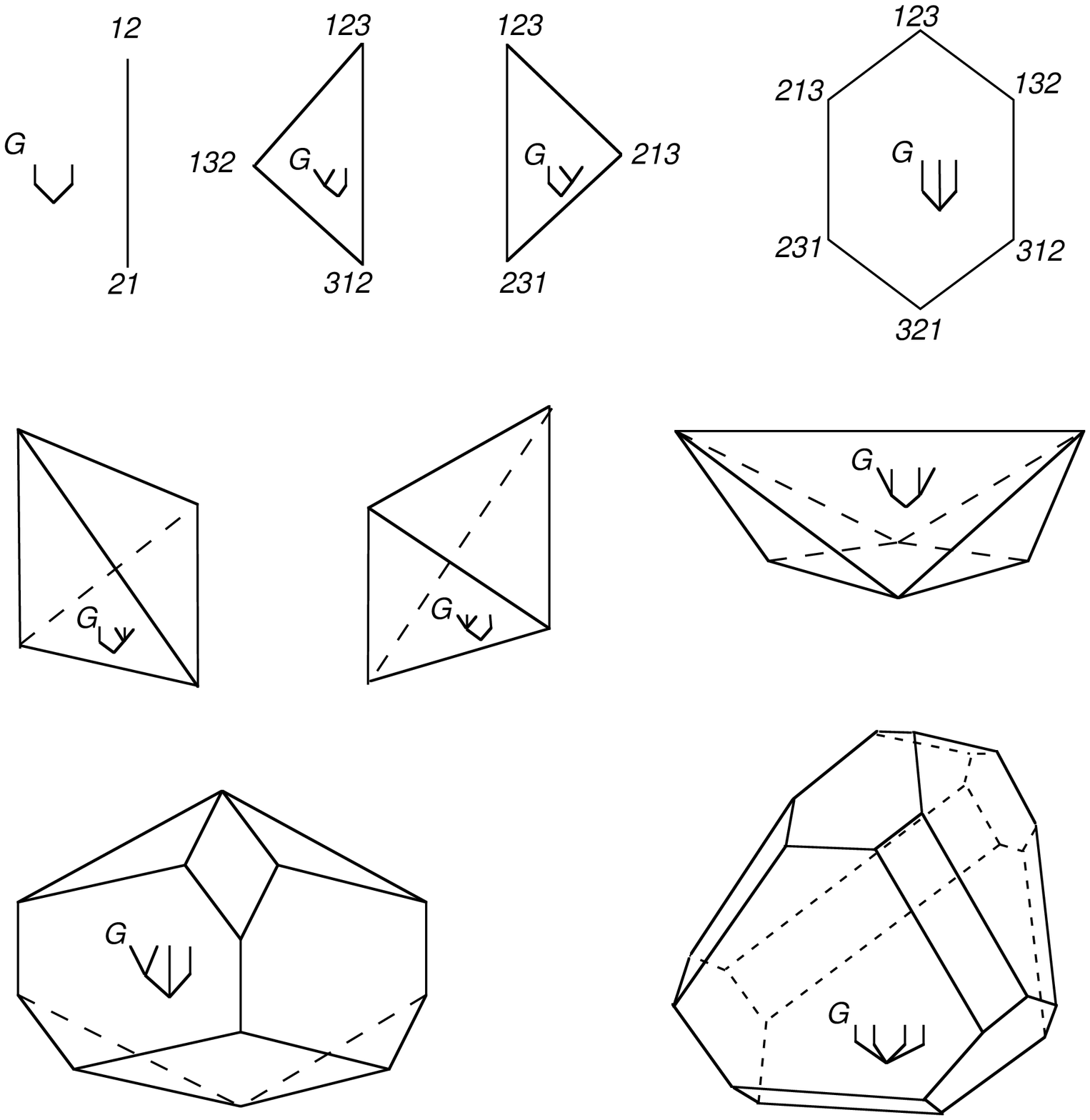}}}

In general, the groupoid $G_T$ is generated by the $1$-faces of the so-called
$T$-shuffle convex polytopes $P_T$ \cite{BatSh} which themselves form a
topological $2$-operad. The polytope $P_T$ is a point if $T$ has only one tip. The polytope
$P_{M_0^j}$ is the permutohedron
$P_j$, and  the polytope $P_{M_1^p\otimes_0 M_1^q}$ is the resultohedron $N_{pq}$
\cite{KV,GKZ}. We finish this example by presenting a picture for multiplication in
$P$ (or $G$ if you like). The reader might find this picture somewhat familiar. 

{\epsfxsize=280pt 
\makebox(290,220)[r]{\epsfbox{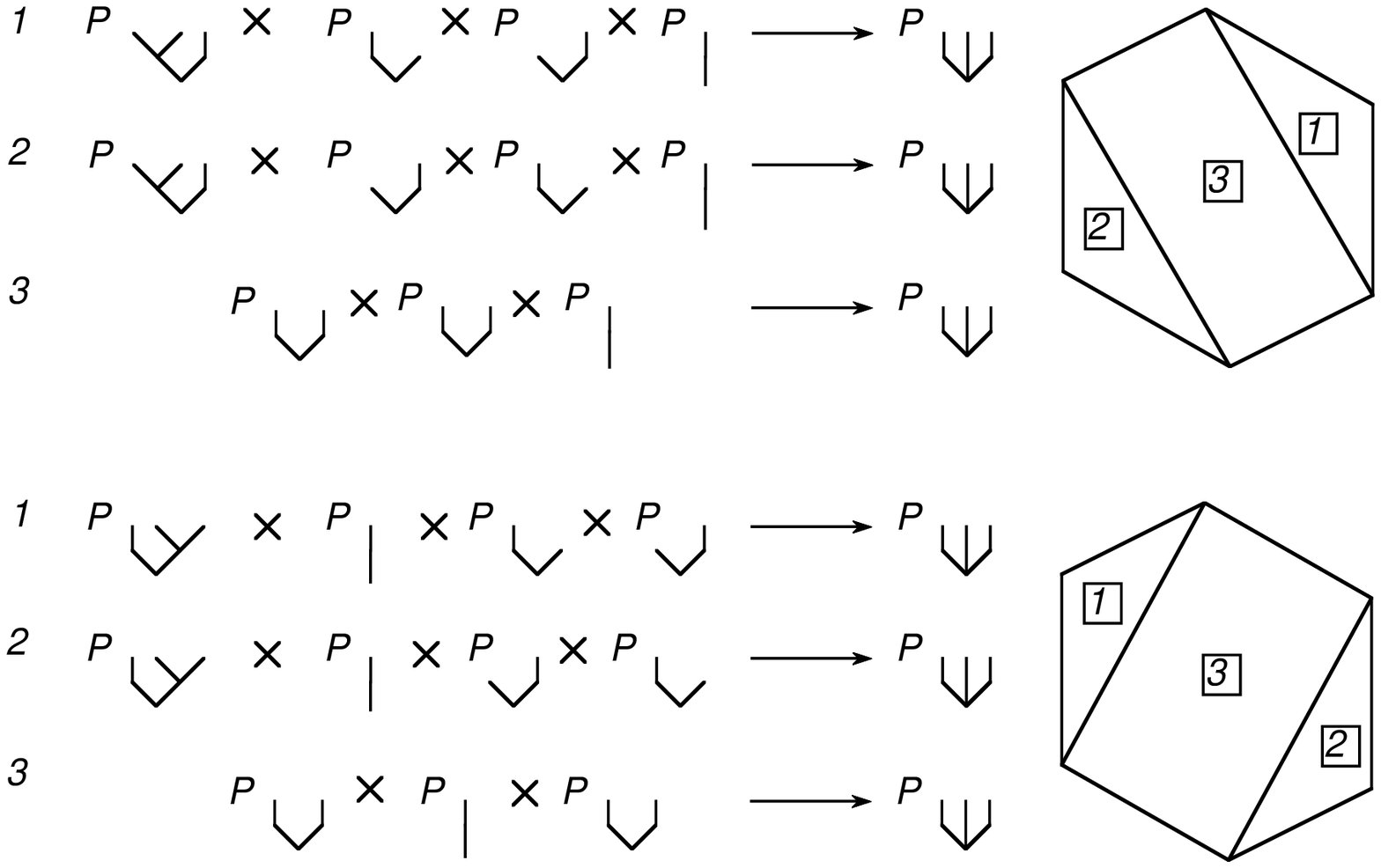}}}

There are similar pictures for higher dimensions. In general the multiplication
in $P$ produces some subdivisions of $P_T$ into products of shuffle polytopes of low 
dimensions. Some special cases of these subdivisions were discovered  in \cite{KV}.  Two new examples
are presented below.

{\epsfxsize=280pt 
\makebox(290,160)[r]{\epsfbox{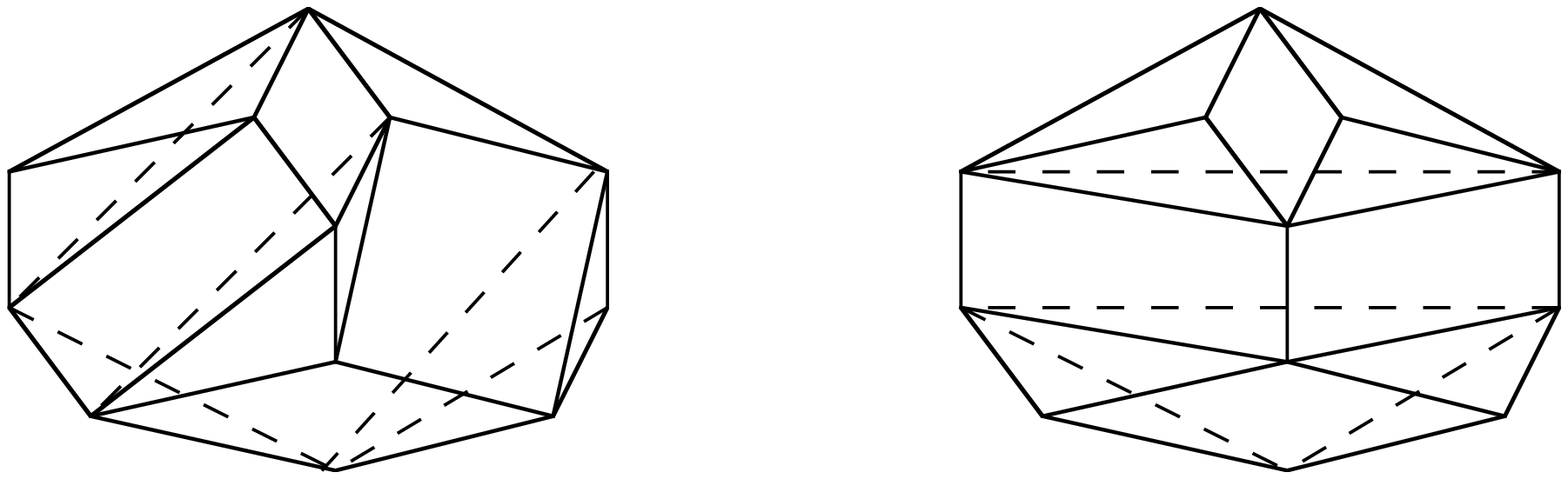}}}

\section{Desymmetrisation of symmetric operads}\label{dessym}

We define the desymmetrisation of an $S$-operad $A\in O^s(V)$ as a pullback along the functor (\ref{[]}). 

However, since we prefer to work with left symmetric operads the following explicit definition of the desymmetrisation functor will be used in the rest of this paper.

 Let $A$ be a left symmetric operad in $V$ with  multiplication $m$ and unit $e$. Then an
$n$-operad 
$Des_n(A)$ is defined by
$$Des_n(A)_T = A_{|T|},$$  
with unit morphism $$e:I\longrightarrow A_{U_n} = A_1$$
and multiplication
$$ m_{\sigma}: A_{|S|}\otimes A_{|T_1|}\otimes ... \otimes A_{|T_k|}\stackrel{m}{\longrightarrow}
A_{|T|}\stackrel{\pi(\sigma)^{-1}}
{\longrightarrow} A_{|T|}$$ 
for  $\sigma: T\rightarrow S$. 
 We, therefore,  have a functor
 $$Des_n:SO_l(V) \longrightarrow O_n (V).$$

\

Now we will consider how the desymmetrisation functor acts on endomorphism operads. 

Let us recall a construction from \cite{BatN}. If $M$ is a monoidal globular category
then a corollary of the coherence theorem for monoidal globular categories \cite{BatN}
says that $M$ is a pseudo algebra of the monad $D$ on the $2$-category of globular
categories. So we have an action $k: D(M)\rightarrow M$ and an isomorphism  in the
square

{\unitlength=1mm

\begin{picture}(60,33)

\put(45,25){\makebox(0,0){\mbox{$D^2(M)$}}}
\put(45,21){\vector(0,-1){10}}
\put(53,26){\shortstack{\mbox{$D(k) $}}}

\put(67,25){\makebox(0,0){\mbox{$D(M) $}}}
\put(67,21){\vector(0,-1){10}}

\put(52,25){\vector(1,0){9.5}}

\put(68,15){\shortstack{\mbox{$ k$}}}

\put(41,15){\shortstack{\mbox{$\mu $}}}

\put(45,7){\makebox(0,0){\mbox{$D(M)$}}}

\put(67,7){\makebox(0,0){\mbox{$M$}}}

\put(52,7){\vector(1,0){10}}

\put(55,8){\shortstack{\mbox{$ k$}}}

\put(58.6,12.4){\line(1,1){4}}
\put(58,13){\line(1,1){4}}

\put(58,12.4){\line(1,0){2}}
\put(58,12.4){\line(0,1){2}}

\end{picture}}

If $x: 1\rightarrow M$ is a globular object of $M$ then the composite 
$$t: D(1) \stackrel{D(x)}{-\!\!\!\longrightarrow} D(M)\stackrel{k}{\longrightarrow} M$$
can be considered as a globular version of the tensor power functor. The value of this functor
on a tree $T$ is denoted by $x^T$. Moreover, the square above gives us an isomorphism $\chi$:  
  
{\unitlength=1mm

\begin{picture}(60,33)

\put(45,25){\makebox(0,0){\mbox{$D^2(1)$}}}
\put(45,21){\vector(0,-1){10}}
\put(53,26){\shortstack{\mbox{$D(t) $}}}

\put(67,25){\makebox(0,0){\mbox{$D(M) $}}}
\put(67,21){\vector(0,-1){10}}

\put(52,25){\vector(1,0){9.5}}

\put(68,15){\shortstack{\mbox{$k$}}}

\put(41,15){\shortstack{\mbox{$\mu $}}}

\put(45,7){\makebox(0,0){\mbox{$D(1)$}}}

\put(67,7){\makebox(0,0){\mbox{$M$}}}

\put(52,7){\vector(1,0){10}}

\put(55,8){\shortstack{\mbox{$ t$}}}

\put(58.6,12.4){\line(1,1){4}}
\put(58,13){\line(1,1){4}}

\put(58,12.4){\line(1,0){2}}
\put(58,12.4){\line(0,1){2}}
\put(57,16){\shortstack{\mbox{$\chi$}}}

\end{picture}}

\noindent This isomorphism $\chi$  gives a canonical isomorphism
$$\chi: x^{T\otimes_l S} \rightarrow x^T \otimes_l x^S , $$ 
for example.

In the
special case of $M= \Sigma^n (V),$ we identify  globular objects of $\Sigma^n (V)$ with
objects of $V$ and follow the constructions from
\cite{BatN} to get the following    inductive description of the tensor
power functor and isomorphism $\chi$.

For the $k$-tree $T$, $k\le n$, and an object $x$ from $V$, let us define the
object 
$x^T$ in the following inductive way:

- if $k< n$, then $x^T = I$;

- if $k=n$
and $T= zT'$ then $x^T = I$;

- if $k=n$ and $T = U_n$, then $x^T = x$;

- now we use  induction on $susp(T)$ : suppose we already have defined 
$x^S$ for $S$ such that $susp(S)\ge k+1$, and let $T= T_1\otimes_k
\ldots \otimes_k T_j$ be a canonical decomposition of $T$. Then we
define
$$x^T = x^{T_1}\otimes \ldots \otimes x^{T_j}.$$

Clearly, $x^T$ is  isomorphic to $\underbrace{x \otimes \ldots
\otimes x}_{|T|-\mbox{\small times}}$ . 
Now we want to provide an explicit description of $\chi$. 

\begin{lem}\label{7.2} For  $\sigma:T\rightarrow S,$ the isomorphism 
$$\chi_{\sigma}: x^{T_1}\otimes \ldots \otimes x^{T_k} \rightarrow x^T $$
is induced by the permutation inverse of the permutation
$\pi(\sigma)$. \end{lem}
\noindent {\bf Proof.}
We will prove the lemma 
by induction. 

If $S=U_n$ then $\chi_{\sigma}$ is the identity morphism.  Suppose we already have
proved our lemma   for all $\sigma$'s with codomain being an $(l+1)$-fold
suspension.  As a first step we study $\chi_{\sigma}$ in the special case  $\sigma:
T\rightarrow M_l^k$.

Now we start another induction on $susp(T)$.
If $susp(T)>l$ then  $\sigma$ factorises through one of the tips, so the fibers are either
$T$ or degenerate trees and we get
$$\chi_{\sigma} = id: I\otimes \ldots \otimes x^T\otimes \ldots I \rightarrow x^T.$$

If $susp(T)=l$ then we get 
$$\chi_{\sigma} = id : x^{T_1}\otimes \ldots \otimes x^{T_k} \rightarrow x^T.$$ 
Suppose we already have proved our lemma  for all $T$ with $susp(T)>m$. 
Now suppose we have a $\sigma$ with  $susp(T)= m  < l$. In this case we have the canonical
decomposition 
$$T_i = T_i^1\otimes_m \ldots \otimes_m T_i^j ,  $$
where $j$ is the same for all $1\le i \le k \ .$
Then  $\chi_{\sigma}$ is equal to the composite
$$x^{T_1}\otimes \ldots \otimes x^{T_k} = (x^{T_1^1}\otimes \ldots \otimes x^{T_1^j})\otimes
\ldots \otimes (x^{T_k^1}\otimes \ldots \otimes x^{T_k^j})
\stackrel{\pi^{-1}}{\longrightarrow}$$
$$\stackrel{\pi^{-1}}{\longrightarrow} (x^{T_1^1}\otimes \ldots \otimes x^{T_k^1})\otimes
\ldots \otimes (x^{T_1^j}\otimes \ldots \otimes x^{T_k^j})
\stackrel{\scriptscriptstyle \chi_1\otimes\ldots\otimes\chi_j}{\longrightarrow}
$$
$$  
\stackrel{\scriptscriptstyle \chi_1\otimes\ldots\otimes\chi_j}{\longrightarrow}
x^{{T_1^1}\otimes_l \ldots \otimes_l {T_k^1}}\otimes \ldots \otimes x^{{T_1^j}\otimes_l
\ldots
\otimes_l {T_k^j}} = x^T$$
where $\pi$ is the corresponding permutation and $\chi_i$ is already constructed by the
inductive hypothesis as $susp({T_1^i}\otimes_l \ldots \otimes_l {T_k^i})>m. $ Again by
induction 
$\chi_i$ is induced by  the permutation $\pi(\phi_i)^{-1}$, where
$$\phi_i: {T_1^i}\otimes_l \ldots \otimes_l {T_k^i} \rightarrow M_l^k .$$
So $\chi_{\sigma}$ is induced by $\Gamma(\pi(\omega);\pi(\phi_1),\ldots,\pi(\phi_j))^{-1}.$

From the point of view of morphisms in $\Omega_n$ what we have used here is a decomposition 
of $\sigma$ into 
$$T\stackrel{\xi}{\longrightarrow} M_m^j\otimes_l\ldots \otimes_l  M_m^j
\stackrel{\omega}{\longrightarrow} M_l^k.$$
Then we have $\pi = \Gamma(\pi(\omega);1,\ldots,1)$.  By construction we have
$\pi(\xi_i)=1$ and by the inductive hypothesis,
 $\pi(\xi) = \Gamma(1;\pi(\phi_1),\ldots,\pi(\phi_j)).$
By  Lemma \ref{pisigma} $$\pi(\sigma) = \pi(\xi)\cdot\Gamma(\pi(\omega);1,\ldots,1) =\Gamma(1;\pi(\phi_1),\ldots,\pi(\phi_j))\cdot
\Gamma(\pi(\omega);1,\ldots,1)  = $$ $$
= \Gamma(\pi(\omega);\pi(\phi_1),\ldots,\pi(\phi_j)).$$
So we have completed our first induction. 

To complete the proof it remains to show the lemma when $S$ is an arbitrary tree with 
$susp(S)= l$. Then we have a  canonical decomposition $\omega:S\rightarrow M_l^j$ and we can
form the composite
$$T\stackrel{\sigma}{\longrightarrow} S\stackrel{\omega}{\longrightarrow} M_l^j.$$
Since $\omega$  is an order preserving map, by  Lemma \ref{pisigma} again, 
$$\pi(\omega)= \pi(\sigma\cdot\omega)\cdot\Gamma(1;\pi(\sigma_1),\ldots,\pi(\sigma_j)).$$
By the inductive hypothesis again we can assume that we already have proved our lemma for the
$\sigma_i$'s and, by the previous argument, for $\sigma\cdot\omega$ as well. So the result
follows. 

\

\Q

\

We now recall the construction of endomorphism operad from \cite{BatN} in the  special case of
an augmented monoidal globular category equal to $\Sigma^n V$, where  
 $V$ is a  closed symmetric monoidal category. 

Let $x$ be
an object of $V$. The endomorphism $n$-operad of $a$ is the following
$n$-operad $End_n(x)$ in $V$. For a tree $T$, 
$$End_n(x)_{_{\scriptstyle T}}  = V(x^T,x) \ ; $$
the unit of this operad is given by the identity
$$ I\rightarrow V(x^{U_n},x) = V(x,x). $$
of $x$. For a morphism $\sigma:T\rightarrow S ,$
the multiplication is given by
$$ V(x^S,x)\otimes V(x^{T_1},x)\otimes \ldots \otimes V(x^{T_k},x)  \rightarrow
 V(x^S,x)\otimes V(x^{T_1}\otimes \ldots \otimes x^{T_k},x^k) \rightarrow$$
$$\stackrel{1\otimes V(\chi^{-1}_{\sigma},x^k)}{-\!\!\!-\!\!\!-\!\!\!\longrightarrow} 
V(x^S,x)\otimes V(x^T,x^k)
\longrightarrow V(x^T,x).$$ 

We also can form the usual symmetric endooperad of $x$ in the  symmetric closed monoidal
category
$V$ \cite{May}. Let us denote this operad by $End(x)$.
Now we want to  compare $End_n(x)$  with the $n$-operad $Des_n(End(x))$. In $End(x)$ the action
of a bijection $\pi:[n]\rightarrow [n]$ is defined to be 
$$V(x^n,x)\stackrel{V(\pi,x)^{-1}}{-\!\!\!-\!\!\!-\!\!\!-\!\!\!\longrightarrow} V(x^n,x).$$   
So, for  $\sigma:T\rightarrow S ,$ we have in $Des_n (End(x))$ the multiplication
$$ V(x^{|S|},x)\otimes V(x^{|T_1|},x)\otimes \ldots \otimes V(x^{|T_k|},x) \rightarrow V(x^{|S|},x)\otimes
V(x^{|T_1|+ \ldots + {|T_k|}},x^k)  \rightarrow $$
$$ \longrightarrow
V(x^{|T|},x)\stackrel{V(\pi(\sigma),x)}{-\!\!\!-\!\!\!-\!\!\!-\!\!\!\longrightarrow}
V(x^{|T|},x).$$ The following proposition follows now from   Lemma \ref{7.2} and associativity of composition in $V .$  

\begin{pro} \label{endomorphism} For any object $x\in V$, there is a natural isomorphism of $n$-operads
$$End_n(x) \simeq Des_n(End(x)).$$ \end{pro}

\section{Internal algebras of cartesian monads}\label{InternalAlgebras}

Suppose $C$ is a category with finite limits. Recall that a monad $(T, \mu , \epsilon)$ on $C$ is called cartesian if $T$ preserves pullbacks and $\mu$ and $\epsilon$ are cartesian natural transformations in the sense that all naturality squares for these transformations are pullbacks. 

If  $(T, \mu , \epsilon)$ is a cartesian monad 
then it can be extended to a $2$-monad ${\bf T= (T, \mu , \epsilon)}$ on the $2$-category ${\bf Cat}(C)$ of internal categories in $C .$ Slightly abusing notation we will speak about categorical $T$-algebras,   having in mind $\bf T$-algebras, if it does not lead to confusion.  So we can speak of pseudo-$T$-algebras, strict morphisms or simply morphisms between $T$-algebras as well as strong or pseudo morphisms and lax-morphisms \cite{BPK,Lack}. The last notion requires some clarification because of the choice 
 of direction of the structure cells. So we give the following definition:
 
 \begin{defin} Let $A$ and $B$ be two categorical $T$-algebras. Then a lax-morphism  $(f,\phi): A\rightarrow B$ is a functor $f: A\rightarrow B$ together with a natural transformation

{\unitlength=1mm

\begin{picture}(60,33)
\put(45,25){\makebox(0,0){\mbox{$T(A)$}}}
\put(45,21){\vector(0,-1){10}}
\put(53,26){\shortstack{\mbox{$ $}}}

\put(70,25){\makebox(0,0){\mbox{$T(B) $}}}
\put(70,21){\vector(0,-1){10}}

\put(53,25){\vector(1,0){12}}

\put(45,7){\makebox(0,0){\mbox{$A$}}}

\put(70,7){\makebox(0,0){\mbox{$B$}}}

\put(53,7){\vector(1,0){12}}

\put(55.8,12.4){\line(1,1){4}}
\put(55,13){\line(1,1){4}}

\put(55,12.2){\line(1,0){2}}
\put(55,12.2){\line(0,1){2}}
\put(54,16){\shortstack{\mbox{$\phi$}}}

\end{picture}}

\noindent which must satisfy the usual coherence conditions \cite{BPK,Lack}.

\end{defin}

 Notice that our terminology  is different from \cite{BPK} but identical with \cite{Lack} :  we call pseudo morphism  what in \cite{BPK} is called morphism  between $T$-algebras. We introduce the following notations: $Alg_T$ is the category of algebras of $T$, while ${\bf CAlg}_{T}$ is the $2$-category of categorical $T$-algebras   and   strict categorical $T$-algebras  morphisms. Notice  also that   ${\bf CAlg}_{T}$ is isomorphic to  the $2$-category of internal categories in the category of $T$-algebras. 

 
We first make the following observation about algebras and pseudoalgebras in our settings\footnote{ I  would like  to thank S.Lack for explaining to me that the proof of Power's general coherence result  works in this situation. }.

\begin{theorem} Every  pseudo-$T$-algebra  is equivalent to a $T$-algebra in the $2$-category of categorical pseudo-$T$-algebras and their pseudo morphisms. \end{theorem}

\Proof  We can easily  adapt  the proof of the general coherence result from \cite{Power} since our monad $T$ has the property of preserving functors which are isomorphisms on objects.
 That is, if $f:A\rightarrow B$ is such that $f$ is an isomorphism on the objects of objects then the same is true for $Tf .$ 

\

\Q

\
 
In practice, the  pseudo-$T$-algebras are as  important for us as  strict $T$-algebras but  in virtue of this theorem we should not worry about this difference.

\
 
 Now observe, that the terminal category $1$ is always a categorical $T$-algebra.
  
  \begin{defin} Let $A$ be  a categorical  $T$-algebra. An internal $T$-algebra $a$  in $A$ is given by a lax-morphism of $T$-algebras 
  $$a: 1 \rightarrow A .$$
  We have a notion of a natural transformation between internal $T$-algebras and so a category of internal 
  $T$-algebras $IAlg_T(A) .$  
  \end{defin}
  
Obviously, $IAlg_T(A)$ can be extended to   a $2$-functor $$IAlg_T: {\bf CAlg}_{T} \rightarrow Cat  .$$ 
 
\begin{theorem}\label{representInt1}
The $2$-functor  $IAlg_T$ is representable. The representing categorical $T$-algebra  $\H^T$ has a characteristic property that its simplicial nerve coincides with May's two sided bar construction $B_{\star}(T,T,1) $ i.e. with the cotriple resolution 
of the terminal $T$-algebra. \end{theorem}
\Proof 

Consider the following part of the cotriple resolution of the terminal    $T$-algebra in $Alg_T$ 

 {\unitlength=1mm

\begin{picture}(60,12)



\put(27,5){\makebox(0,0){\mbox{$T(1)$}}}

\put(39.5,6.5){\vector(-1,0){7}}
\put(39.5,3.5){\vector(-1,0){7}}
\put(32.5,5){\vector(1,0){7}}

\put(36,6.5){\shortstack{\mbox{\small $ $}}}

\put(45,5){\makebox(0,0){\mbox{$T^2(1)$}}}

\put(63,5){\makebox(0,0){\mbox{$T^3(1)$}}}

\put(57.5,6.5){\vector(-1,0){7}}
\put(57.5,5){\vector(-1,0){7}}
\put(57.5,3.5){\vector(-1,0){7}}



\end{picture}}

 Since $T$ is cartesian, the  object above is a truncated nerve of a categorical object $\H^T$ in $Alg_T .$ The Segal's conditions follow from the  naturality square for the multiplication of $T$ being a pullback.

Let us prove that $\H^T$ is a strict codescent object \cite{Lack,StD} of the  terminal  categorical $T$-algebra that is an appropriate weighted colimit of the following  diagram  ${\bf T}^{\star}(1) :$

 {\unitlength=1mm

\begin{picture}(60,12)



\put(27,5){\makebox(0,0){\mbox{${\bf T}(1)$}}}

\put(39.5,6.5){\vector(-1,0){7}}
\put(39.5,3.5){\vector(-1,0){7}}
\put(32.5,5){\vector(1,0){7}}
\put(36,6.5){\shortstack{\mbox{\small }}}

\put(45,5){\makebox(0,0){\mbox{${\bf T}^2(1)$}}}

\put(63,5){\makebox(0,0){\mbox{${\bf T}^3(1)$}}}

\put(57.5,6.5){\vector(-1,0){7}}
\put(57.5,5){\vector(-1,0){7}}
\put(57.5,3.5){\vector(-1,0){7}}




\end{picture}\label{codesc}}

Recall \cite{Lack,StD} that for a truncated cosimplicial category
$E^{\star}$

 {\unitlength=1mm

\begin{picture}(60,12)



\put(27,5){\makebox(0,0){\mbox{$E^0$}}}

\put(32.5,7.5){\vector(1,0){7}}
\put(32.5,2.5){\vector(1,0){7}}

\put(39.5,5){\vector(-1,0){7}}

\put(35,8.2){\shortstack{\mbox{\small $\scriptscriptstyle  d_0$}}}
\put(35,3.2){\shortstack{\mbox{\small $\scriptscriptstyle  d_1$}}}
\put(35,5.7){\shortstack{\mbox{\small $\scriptscriptstyle  s_0$}}}

\put(45.5,5){\makebox(0,0){\mbox{$E^1$}}}

\put(63,5){\makebox(0,0){\mbox{$E^2$}}}

\put(50.5,7.5){\vector(1,0){7}}
\put(50.5,5){\vector(1,0){7}}
\put(50.5,2.5){\vector(1,0){7}}

\put(52.5,8.2){\shortstack{\mbox{\small $\scriptscriptstyle  d_0$}}}
\put(52.5,3.2){\shortstack{\mbox{\small $\scriptscriptstyle  d_2$}}}
\put(52.5,5.7){\shortstack{\mbox{\small $\scriptscriptstyle  d_1$}}}



\end{picture}}

\noindent one can construct  the descent category $Desc(E^{\star})$ whose objects are pairs $(a,f)$ where $a$ is an object of $E^0$ and 
$f:d_0(a)\rightarrow d_1(a)$ is a morphism of $E^1$ satisfying the conditions that $s_0(f)$ is the identity morphism of $a$ and $d_1(f)$ is the composite of $d_2(a)$ and $d_0(f) .$ A morphism in 
$Desc(E^{\star})$ from $(a,f)$ to $(b,g)$ is a morphism $u:a\rightarrow b$ such that $d_0(u)\cdot f = g \cdot d_1(u) .$ 

Let $A$ be a categorical $T$-algebra and let 
$$E^{\star} = {\bf CAlg}_{ T}({\bf T}^{\star}(1), A).$$  
A direct verification demonstrates that   
$Desc(E^{\star})$ is isomorphic to the category ${\bf CAlg}_{ T}(\H^T, A) 
\footnote{ As observed in \cite{StD} this is a general fact about the descent category of a truncated cosimplicial category obtained as $Hom(Ner(X),A)$ where $Ner(X)$ is a nerve of a category considered as a discrete truncated simplicial  category.}.$ Therefore, $H^T$ is a codescent object of ${\bf T}^{\star}(1).$   

On the other hand  $${\bf CAlg}_{ T}({\bf T}^{k}(1), A)
\simeq {\bf Cat}(C)({\bf T}^{k-1}(1),A)$$ 
and the data for the objects and morphisms in $Desc(E^{\star})$ amount to the data for a lax-morphisms and their transformations
from $1$ to $A$ (see \cite{Lack} for   general and detailed consideration). 

Finally, the simplicial nerve of $H^T$ coincides with the bar construction due to the fact that $T$ is cartesian and, hence, all Segal's maps are isomorphisms

 \
 
 \Q
 
 \
 
 \begin{cor} The categorical $T$-algebra $\H^T$ has a terminal object given by its canonical internal $T$-algebra. In particular, it is contractible. 
 \end{cor}
  
\Proof  The terminality of the  internal algebra of $\H^T$  follows from the pullback of naturality for the unit of the monad $T .$

\

\Q

\

\Example Let $C = Set$ and $M$ be the free monoid monad. It is well known that $M$ is finitary and cartesian. The categorical pseudoalgebras of $M$ are equivalent to  monoidal categories. Then an internal $M$-algebra in a  monoidal category $V$  is just a monoid  in $V .$   
The category $\H^M$ is the category $\Delta_{alg} = \Omega_1 $ of all finite ordinals. 

\

\Example Let $C= Glob_n$ be the category of $n$-globular sets \cite{BatN} and let $D_n$ be the 
free $n$-category  monad on $Glob_n $ \cite{BatN}.  $D_n$ is cartesian and finitary \cite{StP}. The algebras of $D_n$ are $n$-categories, the categorical algebras are strict globular monoidal categories and pseudoalgebras are equivalent to globular monoidal categories \cite{BatN}. An internal $D_n$-algebra was called an $n$-globular monoid in a monoidal globular category. The category $\H^{D_n}$ is the monoidal globular category of trees (see section \ref{treesandmorphisms}):

 {\unitlength=1mm

\begin{picture}(60,12)
\put(25,5){\makebox(0,0){\mbox{$1= \Omega_0 $}}}

\put(39,5.5){\vector(-1,0){7}}

\put(39,4.5){\vector(-1,0){7}}

\put(36,6){\shortstack{\mbox{\small $ t$}}}
\put(36,2.5){\shortstack{\mbox{\small $ s$}}}

\put(43,5){\makebox(0,0){\mbox{$ \Omega_1$}}}

\put(53,5.5){\vector(-1,0){7}}

\put(53,4.5){\vector(-1,0){7}}

\put(50,6){\shortstack{\mbox{\small $ t$}}}
\put(50,2.5){\shortstack{\mbox{\small $ s$}}}

\put(54,3.5){\shortstack{\mbox{ $\ldots$}}}

\put(68,5.5){\vector(-1,0){7}}

\put(68,4.5){\vector(-1,0){7}}

\put(65,6){\shortstack{\mbox{\small $ t$}}}
\put(65,2.5){\shortstack{\mbox{\small $ s$}}}

\put(73,5){\makebox(0,0){\mbox{$\Omega_n  \ \ .$}}}

\end{picture}}

\

Now, suppose we have two  finitary  monads $S$ and $T$ on cocomplete  categories $C$ and $E$ respectively.  
Suppose also that there is a right adjoint $w: C \rightarrow E$ and a functor $d:Alg_S\rightarrow Alg_T$ making the following square commutative:

 \begin{center} {\unitlength=1mm
\begin{picture}(60,26)
\put(11,20){\makebox(0,0){\mbox{$Alg_S $}}}
\put(11,16){\vector(0,-1){10}}

\put(37,20){\makebox(0,0){\mbox{$Alg_T$}}}
\put(37,16){\vector(0,-1){10}}

\put(17,20){\vector(1,0){14}}
\put(22,22){\shortstack{\mbox{$\scriptstyle d $}}}

\put(68,10){\shortstack{\mbox{$ $}}}
\put(35,10){\shortstack{\mbox{$ $}}}
\put(11,3){\makebox(0,0){\mbox{$C$}}}
\put(37,3){\makebox(0,0){\mbox{$E$}}}
\put(17,4){\vector(1,0){14}}
\put(24,5.4){\shortstack{\mbox{$\scriptstyle w $}}}

\put(13,10){\shortstack{\mbox{$\scriptstyle U^S $}}}

\put(40,10){\shortstack{\mbox{$\scriptstyle U^T $}}}

\end{picture}}
\end{center}

\begin{pro}\label{adjsquare} The square above induces a commutative square of left adjoints. All together these adjunction can be included in a square which we will refer to as a commutative square of  adjunctions

 \begin{center} {\unitlength=1mm
\begin{picture}(60,26)
\put(10,20){\makebox(0,0){\mbox{$Alg_S $}}}
\put(9,6){\vector(0,1){10}}
\put(11,16){\vector(0,-1){10}}

\put(38,20){\makebox(0,0){\mbox{$Alg_T$}}}
\put(37,6){\vector(0,1){10}}
\put(39,16){\vector(0,-1){10}}

\put(17,21){\vector(1,0){14}}
\put(31,19){\vector(-1,0){14}}
\put(22,22){\shortstack{\mbox{$\scriptstyle d $}}}
\put(22,17){\shortstack{\mbox{$\scriptstyle p $}}}

\put(68,10){\shortstack{\mbox{$ $}}}
\put(35,10){\shortstack{\mbox{$ $}}}
\put(10,3){\makebox(0,0){\mbox{$C$}}}
\put(38,3){\makebox(0,0){\mbox{$E$}}}
\put(17,4.5){\vector(1,0){14}}
\put(31,2.2){\vector(-1,0){14}}
\put(24,5.4){\shortstack{\mbox{$\scriptstyle w $}}}
\put(24,0.1){\shortstack{\mbox{$\scriptstyle c $}}}

\put(4,10){\shortstack{\mbox{$\scriptstyle {\cal F}^S $}}}
\put(12,10){\shortstack{\mbox{$\scriptstyle U^S $}}}

\put(32,10){\shortstack{\mbox{$\scriptstyle {\cal F}^T $}}}
\put(40,10){\shortstack{${\scriptstyle U^T}\addtocounter{equation}{1}
 \ \ \ \ \ \ \ \ \ \ \ \ \ \ \ \  \ \ \ \ \ \ \ \ \  \ \ \  \ \  (\theequation) $}}

\end{picture}}
\end{center}\label{1sq1}

\end{pro}

\Proof This is the Adjoint Lifting Theorem 4.5.6 from \cite{Baurceux} and it also  follows from   Dubuc's adjoint triangle theorem \cite{ED} but for the sake of completeness we provide a proof below.

The problem here is to construct  a functor $p $ left  adjoint to $d .$ Immediately from the requirement of commutativity of the left adjoint square we  have $p \ \F^T \simeq \F^S c $  if $p$ exists.  We use this relation  as a definition of  $p$ on free algebras of the monad $T.$  Notice also that from our assumption of finitarity and cocompletness we get cocompleteness of the category $Alg_T .$

Let $X$ be an arbitrary algebra of  $T.$ Then $X$ is a canonical coequaliser in $Alg_T :$ 
$$ X \leftarrow T(X) \  \raisebox{-0.4mm}{\mbox{$\stackrel{\longleftarrow}{\scriptstyle \longleftarrow}$}} \ T^2(X).$$
The left adjoint $p$ must preserve coequalisers so we define   $p$ on $X $ as the following coequaliser in $Alg_S$:
$$p(X) \leftarrow \F^S c \ U^T (X) \  \raisebox{-0.4mm}{\mbox{$\stackrel{\longleftarrow}{\scriptstyle \longleftarrow}$}} \  \F^S c \ U^T \F^T U^T (X).$$
The first morphism in this coequaliser is induced by the $T$-algebra structure on $X .$ The second morphism is a component of the  natural transformation
\begin{equation}\label{modtrans}  \F^S c \longleftarrow \F^S c \ U^T \F^T \end{equation}
which we construct as follows. 

The existence of the functor $d$ making the first square commutative implies the existence of a natural transformation  
\begin{equation}\label{mapofmon} \Phi = w\  U^S \F^S c \longleftarrow  U^T\F^T \end{equation}
 which  actually can be completed to a map of   monads $\Phi \longleftarrow T .$  \ This gives us an adjoint natural transformation (or mate) $U^S\F^S c\leftarrow  c\ U^T\F^S .$ One more adjoint transformation gives us the transformation we  required.  

It is trivial to check that we have thus constructed a left adjoint to $d .$

\

\Q

\

\begin{corol} There is a canonical map of monads 
                          $$T\rightarrow U^T d \ F^S c \ .$$ 
  \end{corol}                        

If, in addition, $ T$ and $S$ are  cartesian monads, the above  results can be extended to the  categorical level. Abusing notation  once again, we will denote the categorical versions of the corresponding functors by the same letters if it does not lead to a confusion.

\begin{defin} Under the conditions of Proposition \ref{adjsquare},  an internal $T$-algebra inside a categorical  $S$-algebra $A$ is an internal $T$-algebra in $d(A) .$ 
\end{defin}

Let $A$ be a categorical $S$-algebra. The internal $T$-algebras in $A$ form a category $IAlg^S_T(A) .$ Moreover, we have a $2$-functor 
$$IAlg^S_T: {\bf CAlg}_S \rightarrow Cat .$$ 

Let us denote by $G$ the composite $\F^S c  .$ Then the transformation \ref{modtrans} equips  $G$ with the structure of a  module over the monad $T .$ We will  require  
this natural transformation to be cartesian which implies that the map of monads \ref{mapofmon} is also cartesian.

\begin{theorem}\label{representInt2}Let $S$ and $T$ be finitary cartesian monads and assume that the transformation \ref{modtrans} is cartesian. Then   the $2$-functor $IAlg^S_T$ is representable. The representing categorical $S$-algebra $h^T = p(H^T) $ has the characteristic property that its nerve coincides with May's two-sided bar-construction $B_{\star}(G,T,1) $ i.e. with the image under $p$ of the cotriple resolution of the terminal categorical $T$-algebra.

\end{theorem} 

\Proof  The proof is analogous to the proof of the theorem \ref{representInt1} but we use the cartesianess of $T$-action on $G$ to check that the Segal maps in the simplicial object

 {\unitlength=1mm

\begin{picture}(60,12)



\put(28,5){\makebox(0,0){\mbox{$G(1)$}}}

\put(39.5,5.5){\vector(-1,0){7}}
\put(39.5,4.5){\vector(-1,0){7}}

\put(36,6.5){\shortstack{\mbox{\small $ $}}}

\put(46,5){\makebox(0,0){\mbox{$GT(1)$}}}

\put(71,5){\makebox(0,0){\mbox{$GT^2(1) \ \ldots $}}}

\put(59.5,5.8){\vector(-1,0){7}}
\put(59.5,5){\vector(-1,0){7}}
\put(59.5,4.2){\vector(-1,0){7}}



\end{picture}}

\noindent are isomorphisms.

\

\Q

\

\Example  \label{internalobjects}  A trivial case of the adjunction square \ref{1sq1} is a map of monads $I\rightarrow T .$
The functor $d$ in this case is just the forgetful functor and $p = \F^T .$  So one can speak about an internal $I$-algebra inside a categorical $T$-algebra. Such an internal algebra amounts just to a morphism of internal categories 
$$1\rightarrow A .$$ 
 We will call them  {\it internal objects  of $A$ .}  
The corresponding representing $T$-algebra $h^T$ will be denoted by $H_d^T$ since this is  a discretisation of  the categorical $T$-algebra $H^T  $ and  is given by the constant simplicial object $T(1) .$  

\

Observe, that we can extend the functor $d$ from the square \ref{1sq1} to the lax-morphisms between categorical $S$-algebras.
For this we first observe that the map of monads \ref{mapofmon} induces a functor  between corresponding categories of algebras and their lax-morphism. Then we can construct a natural functor from $S$-algebras and their lax-morphisms to $\Phi$-algebras  and their lax-morphisms. We leave the details to the reader.  

Because $d(1) = 1 ,$ we can construct a $2$-natural transformation between $2$-functors which we denote by $\delta$ (but we  think of it as an internal version of the functor $d$): 
$$\delta: IAlg_S \rightarrow IAlg_T^S ,$$ 
which induces  a canonical map between representing objects: 
\begin{equation}\label{restrmap} \zeta: h^T \rightarrow H^S \ . \end{equation}

Another way to construct this map is the following. The algebra $H^S$ has a canonical internal $S$-algebra
$1\rightarrow H^S .$ If we apply $d$ to this lax morphism we will get 
$$1 = d(1) \rightarrow d(H^S) .$$
The last internal $T$-algebra can be represented by a map $H^T \rightarrow d(H^S)$ and by definition this gives a map  \ref{restrmap} . 

Summarizing we have the following
\begin{theorem}  The functor $\delta$ is naturally isomorphic to $\zeta^{\star} ,$ the restriction functor along $\zeta .$ The left adjoint to $\delta$ is isomorphic to the left Kan extension along $\zeta$  in the $2$-category of categorical $S$-algebras. \end{theorem}

\section{Symmetrisation of $n$-operads}\label{Intop} 

Let  $Coll_n(V)$ be the category of $(n-1)$-terminal $n$-collections \cite[section 6]{BatN} i.e. the category of $n$-globular functors $Tr^{(n)}\rightarrow \Sigma^n V.$ We can identify the objects of $Coll_n(V)$ with  families of objects of $V$ indexed by trees of height $n .$ The morphisms are levelwise morphisms. The category of $1$-collection $Coll_1(V)$  is, of course, the same as the category of nonsymmetric collections in the usual sense.

\begin{theorem}\label{cart}  If $V$ is a cocomplete symmetric monoidal category then the forgetful functor $R_n:O_n(V) \rightarrow Coll_n(V)$ is monadic with left adjoint $F_n .$ The free $n$-operad monad $\F_n$ on the category $Coll_n(Set)$ of $Set$ $n$-collections is finitary and cartesian. \end{theorem}

  \Proof 
We first give an inductive construction of the free $(n-1)$-terminal $n$-operad  on an $n$-collection in a cocomplete symmetric monoidal category $V .$

 Let us call an expression, given by an   $n$-tree $T ,$  an {\it admissible} expression of arity $T .$  We also have an admissible expression $e$ of arity $U_n .$   If $\sigma:T\rightarrow S$ is a morphism of  trees and  the admissible expressions   $x_S  , x_{T_1^{}},\ldots,  x_{T_k^{}}$ of arities $S, T_1, \ldots, T_k$ respectively  are already constructed then the expression 
$\mu_{\sigma}(x_S; x_{T_1^{}}, \ldots x_{T_k^{}})$ is also an admissible expression of arity $T .$
We also introduce an obvious equivalence relation on the set of admissible expressions generated by pairs of composable morphisms of  trees  and by two equivalences $T \sim \mu(T;e, \ldots, e) \sim \mu(e;T) $ generated by the identity morphism of $T$ and a unique morphism $T \rightarrow U_n .$ Notice however, that there are morphisms of trees all of  whose fibers  are equal to $U_n .$  We can form an admissible expression $\mu_{\sigma}(S;e, \ldots, e)$  corresponding to such a morphism but it is not equivalent to $S ,$ unless   $\sigma$ is equal to the identity.

Now if $C\in Coll_n(V)$
then, with every admissible expression $\tau$ of arity $T ,$ we can associate  by induction an object  $C(\tau) .$
We start from $C(T) = C_T,  C(e) = I$ and put 
$$C(\mu(x_S; x_{T_1^{}}, \ldots x_{T_k^{}})) = 
C(x_S)\otimes C(x_{T_1^{}})\otimes  \ldots \otimes C( x_{T_k^{}}) .$$
By the Mac Lane coherence theorem, this object  depends  on the equivalence class of an admissible expression only up to isomorphism.  So, we choose a representative of $C(\tau)$ for every equivalence class of admissible expressions. 

Now, the coproduct $\coprod_{\tau} C(\tau)$ over all equivalence classes of admissible expressions of arity $T$ gives us an  $n$-collection in $V .$ We also have a copy of the unit object $I$ of arity $U_n $ which corresponds to the admissible expression  $e .$

It is now a trivial exercise to check that in this way we indeed get a  free  $(n-1)-$ terminal $n$-operad  $F^n(C) $ on $C .$ 

It is also very obvious that the monad $\F_n = R_nF_n$ is finitary and cartesian if $V=Set$. Indeed, every admissible expression $\tau$ determines a non-planar tree decorated by $n$-trees (this tree actually has a canonical planar structure inherited from the planar structure of decorations). Such a decorated tree gives a collection $\alpha(\tau)$ which is empty in arities  which are not equal to any tree which is presented in the decoration of $\tau$    and equal to a $p$ element set $\{1,\ldots,p\}$ in arity $S$ if $S$ is presented in the decoration of $\tau$ exactly $p$ times.  Then 
$$\F_n (C) = \coprod_{\tau} Coll_n(\alpha(\tau), C)$$
so $\F_n$ is finitary and preserves pullbacks \cite{StP} . 

It is an easy exercise to check that the  multiplication and unit of $\F^n$ are cartesian natural transformations

\

\Q

\

We have a  functor 
$$ W_n: Coll_1(V) \rightarrow Coll_n(V)$$
defined on a nonsymmetric collection $A$  as follows: 
$$W_n(A)_T = A_{|T|} .$$
If $V$ has coproducts then $W_n$ has a left adjoint $C_n:$
\begin{equation}\label{adj22}C_n(B)_k = \coprod_{T\in Tr_n \ , \ |T| = k} B_T .\end{equation}

\begin{theorem}\label{cart1} If $V$ is a cocomplete symmetric monoidal category then the forgetful functor $R_{\infty}:SO_l(V) \rightarrow Coll_1(V)$ is monadic. 
The following square of right adjoints commutes

\begin{center} {\unitlength=1mm
\begin{picture}(60,26)
\put(10,20){\makebox(0,0){\mbox{$SO_l(V) $}}}
\put(9,6){\vector(0,1){10}}
\put(11,16){\vector(0,-1){10}}

\put(41,20){\makebox(0,0){\mbox{$O_n(V)$}}}
\put(38,6){\vector(0,1){10}}
\put(40,16){\vector(0,-1){10}}

\put(18,20){\vector(1,0){14}}
\put(22,22){\shortstack{\mbox{$\scriptstyle Des_n $}}}

\put(68,10){\shortstack{\mbox{$ $}}}
\put(35,10){\shortstack{\mbox{$ $}}}
\put(10,3){\makebox(0,0){\mbox{$Coll_1(V)$}}}
\put(41,3){\makebox(0,0){\mbox{$Coll_n(V)$}}}
\put(18,4.5){\vector(1,0){14}}
\put(32,2.2){\vector(-1,0){14}}
\put(24,5.4){\shortstack{\mbox{$\scriptstyle W_n $}}}
\put(24,0.1){\shortstack{\mbox{$\scriptstyle C_n $}}}

\put(3,10){\shortstack{\mbox{$\scriptstyle { F}_{\infty} $}}}
\put(13,10){\shortstack{\mbox{$\scriptstyle R_{\infty} $}}}

\put(33,10){\shortstack{\mbox{$\scriptstyle { F}_n $}}}
\put(42,10){\shortstack{${\scriptstyle R_n}$}} 

\end{picture}}
\end{center}

Therefore, by Proposition \ref{adjsquare}  this square can be completed to  the following commutative square of  adjoints 

\begin{center} {\unitlength=1mm
\begin{picture}(60,26)
\put(10,20){\makebox(0,0){\mbox{$SO_l(V) $}}}
\put(9,6){\vector(0,1){10}}
\put(11,16){\vector(0,-1){10}}

\put(41,20){\makebox(0,0){\mbox{$O_n(V)$}}}
\put(38,6){\vector(0,1){10}}
\put(40,16){\vector(0,-1){10}}

\put(18,21){\vector(1,0){14}}
\put(32,19){\vector(-1,0){14}}
\put(22,22){\shortstack{\mbox{$\scriptstyle Des_n $}}}
\put(22,16.7){\shortstack{\mbox{$\scriptstyle Sym_n $}}}

\put(68,10){\shortstack{\mbox{$ $}}}
\put(35,10){\shortstack{\mbox{$ $}}}
\put(10,3){\makebox(0,0){\mbox{$Coll_1(V)$}}}
\put(41,3){\makebox(0,0){\mbox{$Coll_n(V)$}}}
\put(18,4.5){\vector(1,0){14}}
\put(32,2.2){\vector(-1,0){14}}
\put(24,5.4){\shortstack{\mbox{$\scriptstyle W_n $}}}
\put(24,0.1){\shortstack{\mbox{$\scriptstyle C_n $}}}

\put(3,10){\shortstack{\mbox{$\scriptstyle { F}_{\infty} $}}}
\put(13,10){\shortstack{\mbox{$\scriptstyle R_{\infty} $}}}

\put(33,10){\shortstack{\mbox{$\scriptstyle { F}_n $}}}
\put(42,10){\shortstack{${\scriptstyle R_n} 
 \ \ \ \ \ \ \ \ \ \ \ \ \ \ \ \  \ \ \ \ \ \ \ \ \  \ \ \  \ \  
 $}}

\end{picture}}
\end{center}
The free   symmetric operad monad $\F_{\infty}= R_{\infty}F_{\infty}$ on the category of nonsymmetric $Set$-collections is finitary and cartesian and the canonical right $\F_n$-action on
$F_{\infty}C_n$ is cartesian.   

 \end{theorem}

\Proof The construction of the left adjoint $F_{\infty}(C)$ is classical: 
$$F_{\infty}(C)_n = \coprod_{\tau}C(\tau)$$
where $\tau$ run over the set of planar trees with $n$  marked leaves labelled by the natural numbers from $1$ to $n$. The object
$C(\tau)$ is the  tensor product of the $C_{|v|} ,$ where $v$ runs over all unmarked  vertices of $\tau$ and $|v|$ is the valency (number of input edges)  of $v .$ The symmetric groups act by permutation of labels and the substitution operation is grafting. The properties of $\F_{\infty}$ are obvious.

Now, if $V=Set,$ the composite $F_{\infty}C_n\F_n(C)$  is given by the set of labelled planar trees whose unmarked vertices are decorated by admissible expressions. The number of tips of the arity of the decoration should be equal to the valency of the vertex. As was observed before each admissible expression determines a canonical planar tree. So   an element of $F_{\infty}C_n\F_n(C)$ is given by the following data:
\begin{itemize}
\item a labelled planar tree $\tau$;
\item an assignment of a planar tree $\rho_v$ decorated by $n$-trees    to each internal vertices $v\in \tau$ such that the number of leaves of $\rho_v$ is equal to $|v| ;$  
\item an assignment of an element $c\in C_{T}$ for each $T$ from the decoration of $\rho_v .$ \end{itemize}

Then the action $F_{\infty}C_n\F_n(C) \rightarrow F_{\infty}C_n$ 
consists of gluing together  the planar trees $\rho_v$ according to the scheme provided by $\tau$  and forgetting the decorating $n$-trees. The labeling of leaves  and decorations by elements of $C$ remain in their places. This is obviously a cartesian transformation.

\Q

\

We finish this section by a theorem which will show that the functor $Sym_n$ from the commutative square
from  Theorem  \ref{adj22} is indeed a solution of the symmetrisation problem raised in Section \ref{nfoldsuspension}.

\begin{theorem}\label{sus} For an  $n$-operad $A$ in a closed symmetric monoidal category $V ,$ its symmetrisation $Sym_n(A)$ from  Theorem  \ref{adj22} is a solution of symmetrisation problem in the sense of Definition \ref{symproblem}.
\end{theorem}
\noindent {\bf Proof.} Indeed, for an object $x\in V$, an $A$-algebra structure is given by
a morphism of operads $$k: A\rightarrow End_n (x) . $$ By Proposition \ref{endomorphism}, $$End_n (x) \simeq
Des_n (End(x)),$$ 
and so $k$ determines, and is uniquely determined  by, a map of symmetric operads
$$Sym_n(A)\rightarrow End(x).$$

\Q

\section{Internal operads}\label{internaloperads} 

In virtue of   Theorem \ref{cart1}  we can develop a theory of internal $n$-operads inside categorical $n$-operads as well as consider internal $n$-operads and internal symmetric operads inside symmetric categorical operads.   
We would like to unpack our definition of internal operads and see what they really are on practice. 

Let   $A$ and $B$ be two $n$-operads in  $Cat$.  A {\it lax-morphism $$f:A\rightarrow B$$ }
consists of a
collection of functors
$$f_T:A_T\rightarrow B_T$$
together with natural transformations

{\unitlength=1mm

\begin{picture}(60,33)
\put(25,25){\makebox(0,0){\mbox{$A_S\times
A_{T_1}\times \ldots A_{T_k}$}}}
\put(25,21){\vector(0,-1){10}}
\put(53,26){\shortstack{\mbox{$ $}}}

\put(90,25){\makebox(0,0){\mbox{$B_S\times
B_{T_1}\times \ldots B_{T_k}$}}}
\put(90,21){\vector(0,-1){10}}

\put(43,25){\vector(1,0){28}}

\put(25,7){\makebox(0,0){\mbox{$A_T$}}}

\put(90,7){\makebox(0,0){\mbox{$B_T$}}}

\put(31,7){\vector(1,0){52}}

\put(55.8,12.4){\line(1,1){4}}
\put(55,13){\line(1,1){4}}

\put(55,12.2){\line(1,0){2}}
\put(55,12.2){\line(0,1){2}}
\put(54,16){\shortstack{\mbox{$\mu$}}}

\end{picture}}

\noindent for every $\sigma:T\rightarrow S$ and a morphism $\epsilon:
e_B\rightarrow e_A$, where $e_A,e_B$ are unit objects of $A$ and $B,$
respectively. They must satisfy the usual coherence conditions. If, however,  $\mu$ and $\epsilon$ are identities
the lax-morphism will be called  an operadic morphism (operadic functor).

In the  particular case $B =1 ,$ the terminal $Cat$-operad, a lax-morphism $a: 1\rightarrow A$ is an internal operad in $A .$ Explicitly this  gives the following

\begin{defin} Let $A$ be a categorical  $n$-operad  with
multiplication $m$ and unit object $e\in A_{U_n}$. An internal $n$-operad in $A$
consists of a collection of objects
$a_T
\in A_T\ , \ T\in Tr_n$, together  with a morphism
$$ \mu_{\sigma}: m_{}(a_S;a_{T_1},\ldots, a_{T_k}) \longrightarrow a_T $$
for every morphism of trees $\sigma:T\longrightarrow S$ and a morphism
$$\epsilon: e\longrightarrow a_{U_n}$$
which satisfy obvious conditions analogous to the conditions in the definition of
$n$-operad. 

A morphism $f:a\rightarrow b$ of internal $n$-operads is a collection of morphisms
$f_T: a_T\rightarrow b_T$ compatible with the operadic structures in the obvious sense. \end{defin}

\begin{defin} Let $A$ be a left symmetric operad in $Cat$. Then an internal $n$-operad in
$A$ is an internal $n$-operad in $Des_n(A)$. \end{defin}

 So an internal $n$-operad in a
symmetric categorical operad is given by a collection of objects $a_T \in A_{|T|} \ , \ T\in Tr_n ,$
together with a morphism   
$$ \mu_{\sigma}: m(a_S;a_{T_1},...,a_{T_k}) \longrightarrow \pi(\sigma) a_T $$
for every $\sigma:T\rightarrow S$ 
and $$\epsilon: e \longrightarrow a_{U_n} ,$$
which satisfy associativity and unitary conditions. 

For a notion of internal symmetric operad in a categorical symmetric operad  we have a choice of three different presentations of the category of symmetric operads. For technical reason it will be more convenient for us to use left categorical symmetric operads yet the $S$-version for internal symmetric operads. That is, we consider 
internal algebras for the following square of adjoints
\begin{center} {\unitlength=1mm
\begin{picture}(60,24)
\put(10,20){\makebox(0,0){\mbox{$SO_l(Set) $}}}
\put(9,6){\vector(0,1){10}}
\put(11,16){\vector(0,-1){10}}

\put(41,20){\makebox(0,0){\mbox{$O^s(Set)$}}}
\put(38,6){\vector(0,1){10}}
\put(40,16){\vector(0,-1){10}}

\put(18,21){\vector(1,0){14}}
\put(32,19){\vector(-1,0){14}}
\put(22,22){\shortstack{\mbox{$\scriptstyle  $}}}
\put(22,16.7){\shortstack{\mbox{$\scriptstyle  $}}}

\put(68,10){\shortstack{\mbox{$ $}}}
\put(35,10){\shortstack{\mbox{$ $}}}
\put(10,3){\makebox(0,0){\mbox{$Coll_1(Set)$}}}
\put(41,3){\makebox(0,0){\mbox{$Coll_1(Set)$}}}
\put(20,4.5){\vector(1,0){12}}
\put(32,2.2){\vector(-1,0){12}}
\put(22,5.4){\shortstack{\mbox{$\scriptstyle W = id $}}}
\put(22,0.1){\shortstack{\mbox{$\scriptstyle C = id $}}}

\put(3,10){\shortstack{\mbox{$\scriptstyle { F}_{\infty} $}}}
\put(13,10){\shortstack{\mbox{$\scriptstyle R_{\infty} $}}}

\put(33,10){\shortstack{\mbox{$\scriptstyle { F}^s $}}}
\put(42,10){\shortstack{${\scriptstyle R^s} 
$}}

\end{picture}}
\end{center}

\noindent where horizontal functors are the  isomorphisms of categories described in section \ref{symgroup} and the functors $F^s$ and $R^s$ are determined by commutativity of this square.
  The result of this mixture is the following definition:

\begin{defin} Let $A$ be a symmetric categorical operad. An internal symmetric operad
in $A$ consists of a collection of objects $a_{n}\in A_n \ , \ n\ge 0 ,$ 
together with a morphism    
$$ \mu_{\sigma}: m(a_{k};a_{n_1},...,a_{n_k}) \longrightarrow \pi(\sigma) a_{n} $$
for every $\sigma:[n]\rightarrow [k]$ in $\Omega^s$,
and $$\epsilon: e \longrightarrow a_{1},$$
which satisfy associativity and unitary conditions. \end{defin}

\

\Example Let $C$ be a category. We can consider the endomorphism operad $End(C)$ of $C$ in $Cat$.
An internal $1$-operad $a$ in $C$ is what we call a {\it multitensor } in $C$ \cite{BW}. This is a sequence of
functors
$$a_k:C^k \rightarrow C$$
satisfying the usual associativity and unitarity conditions. If  $a_k \ , \ k\ge 1 \ $, are isomorphisms then $a$ is
just a tensor product on $C$. Conversely, every tensor product on $C$ determines, in an obvious manner, a multitensor on
$C$.

It makes sense to consider categories enriched in a multitensor. In \cite{BW} we show that the
category of algebras of an arbitrary higher operad $A$ in $Span(C)$ is equivalent to the
category of categories enriched over an appropriate multitensor on the category of algebras of
another operad  $B(A)$ which is some sort of  delooping of $A$. 

\

\Example The internal symmetric operads in $End(C)$ were considered by J.McClure and J.Smith
\cite{letter} under the name of functor operads. These operads generalise symmetric monoidal structures on $C$ in
the same way as multitensors generalise monoidal structures.  

\

Let  $\bf CO_n$ be the $2$-category  whose   objects  are categorical $n$-operads,   morphisms are their operadic
morphisms, and the $2$-morphisms are operadic natural transformations. 
We have the $2$-functor
$$IO_n: {\bf CO_n}\rightarrow Cat$$ 
which assigns to an operad $A$ the category of internal $n$-operads in $A$.

Analogously, let $\bf SCO$ be the $2$-category  of left symmetric  categorical operads, their 
operadic functors, and operadic natural transformations. There is the $2$-functor 
$$IO^{sym}_n:{\bf SCO} \rightarrow Cat$$  
  which assigns to an operad $A$ the category of internal $n$-operads in $A$. For $n=\infty ,$ the functor
$IO^{sym}_{\infty}$ assigns  the category of internal symmetric operads in $A$.

\begin{theorem}\label{rep} \ \\
\vspace{-3mm}
\begin{itemize}
\item For every $1\le n < \infty$, there exists a categorical n-operad
$\H^n$ representing the $2$-functor  $IO_n: {\bf CO_n}\rightarrow Cat$; 

\item  there exists a categorical  symmetric operad
$\H^{\infty}$ representing the $2$-functor  $IO^{sym}_{\infty}: {\bf SCO} \rightarrow Cat$; 

\item For every $1\le n \le \infty$, there exists a symmetric $Cat$-operad
$\h^{n}$ representing the $2$-functor  $IO^{sym}_n:{\bf SCO} \rightarrow Cat$;

\item For a categorical symmetric operad $A$ if the left adjoint $sym_n$  to the internal desymmetrisation functor 
$$\delta_n: IO^{sym}_{\infty}(A) \rightarrow IO^{sym}_{n}(A)$$
exists then on an internal $n$-operad $a$ it is isomorphic to the left Kan extension of the representing operadic functor $\tilde{a}$ along 
the canonical morphism of categorical symmetric operads
$$\zeta: \h^n \rightarrow \H^{\infty} .$$

\end{itemize}
 \end{theorem}  

\Example  Let $V$ be a  symmetric strict monoidal category. Consider the following
symmetric 
 categorical operad $ V^{\bullet}$:
$$ V^{\bullet}_n = V,$$
the multiplication is given by iterated tensor product, the unit of
 $ V^{\bullet}$ is the unit object of
 $V$ and the action of the symmetric
groups is trivial. 

\begin{lem}\label{symoperad}  
There are the following isomorphisms of categories
$$IO^{sym}_{\infty}(V^{\bullet}) \rightarrow SO_l(V)$$
$$IO^{sym}_{n}(V^{\bullet}) \rightarrow O_n(V) .$$
The existence of one of the functors $sym_n$ or $Sym_n$ implies the existence of the other and moreover the following diagram commutes:
 
\begin{center} {\unitlength=1mm
\begin{picture}(60,26)
\put(10,20){\makebox(0,0){\mbox{$IO^{sym}_{\infty}(V^{\bullet}) $}}}

\put(10,16){\vector(0,-1){10}}

\put(41,20){\makebox(0,0){\mbox{$IO^{sym}_{n}(V^{\bullet})$}}}
\put(38,16){\vector(0,-1){10}}

\put(20,21){\vector(1,0){10}}
\put(30,19){\vector(-1,0){10}}
\put(23,22){\shortstack{\mbox{$\scriptstyle \delta_n $}}}
\put(22,16.7){\shortstack{\mbox{$\scriptstyle sym_n$}}}

\put(68,10){\shortstack{\mbox{$ $}}}
\put(35,10){\shortstack{\mbox{$ $}}}
\put(10,3){\makebox(0,0){\mbox{$SO_l(V)$}}}
\put(40,3){\makebox(0,0){\mbox{$O_n(V)$}}}
\put(18,4.5){\vector(1,0){14}}
\put(32,2.2){\vector(-1,0){14}}
\put(23,5.4){\shortstack{\mbox{$\scriptstyle Des_n $}}}
\put(23,0.1){\shortstack{\mbox{$\scriptstyle  Sym_n $}}}

\put(3,10){\shortstack{\mbox{$  $}}}
\put(13,10){\shortstack{\mbox{$  $}}}

\put(33,10){\shortstack{\mbox{$\scriptstyle  $}}}

\end{picture}}
\end{center}

\end{lem}

\Proof The proof is an easy exercise in definitions. 

\Q

\

\section{Combinatorial aspects of internal operads}\label{H^nh^n}

We will now show how to construct the categorical operads $\h^n$ and  $\H^n$ combinatorially using theorems \ref{representInt1}, \ref{representInt2}. 
 We concentrate first on the construction of $\h^n$.

Let $Tr_n$ be the result of application of the functor $C_n$ (see (\ref{adj22}))  to the collection $\F_n(1) .$
This collection is  the set of objects of $\Omega_n$ with the grading according to the number of tips. 
Now we can form a free symmetric operad $\F_{\infty}(Tr_n)$
 on this collection. The  elements of $\F_{\infty}(Tr_n)$
are the objects of $\h^{n}.$ 

 Now we want to define morphisms.
 We will do this by providing generators and relations.

Suppose we have a morphism $\sigma:T\rightarrow S$ in $\Omega_n$ and $T_1,..., T_k$
is its list of fibers. Then we will have a generator
$$ \gamma(\sigma):\mu(S;T_1,...,T_k) \rightarrow \pi(\sigma)T$$
where $\mu$ is the multiplication in $\F_{\infty}(Tr_n)$. By the equivariance requirement, we also have morphisms

 $$ \mu(\pi S;\xi_1T_1,...,\xi_k T_k) = \Gamma(\pi; \xi_1,...,\xi_k) \mu(S;T_1,...,T_k)
\longrightarrow $$
$$\stackrel{\Gamma(\pi;
\xi_1,...,\xi_k)\gamma(\sigma)}{-\!\!\! -\!\!\! -\!\!\! -\!\!\! -\!\!\!
-\!\!\! -\!\!\! -\!\!\! -\!\!\! -\!\!\!\longrightarrow} \ 
\pi\cdot\Gamma(\pi; \xi_1,...,\xi_k)(\sigma) T.$$

For every composite 
$$T\stackrel{\sigma}{\rightarrow} S \stackrel{\omega}{\rightarrow}
 R ,$$
we will have a relation given by the commutative diagram:

{\unitlength=1mm

\begin{picture}(200,47)(0,-3)

\put(20,35){\makebox(0,0){\small\mbox{$\mu(\mu( R;
S_{\bullet}); T_1^{\bullet},  ... 
T_i^{\bullet},  ... , T_k^{\bullet})   
$}}}
\put(20,31){\vector(0,-1){12}}

\put(94,31){\vector(0,-1){12}}

\put(84,35){\makebox(0,0){\small\mbox{$\mu(R;
\mu(S_{1}; T_1^{\bullet}) ,  ...
, \mu(S_{i};
 T_i^{\bullet}),  ... \mu( S_{k};
T_k^{\bullet}))   
$ }}}

\put(46,35){\makebox(0,0){\mbox{$= $}}}

\put(20,15){\makebox(0,0){\mbox{\small$\mu (\pi(\omega) S; 
T_1^{\bullet} ,  ...
, 
 T_i^{\bullet},  ...  T_k^{\bullet})
$}}}

\put(94,15){\makebox(0,0){\mbox{\small$\mu(R;\pi(\sigma_1)T_1,... 
, \pi(\sigma_k)T_k)
$}}}

\put(60,3){\makebox(0,0){\mbox{\small$\pi(\sigma)\cdot\Gamma(\pi(\omega);1,...,1)
T =\pi(\sigma\cdot\omega)\cdot\Gamma(1_k;
 \pi(\sigma_1),..., \pi(\sigma_k))
T  
$}}}

\put(30,11){\vector(4,-1){19}}

\put(90,11){\vector(-4,-1){19}}

\end{picture}}

\noindent We also have a generator
$$\epsilon: e \rightarrow U_n$$
and two  commutative diagrams:

{\unitlength=1mm
\begin{picture}(50,25)(28,2)

\put(92,20){\vector(-1,0){15}}

\put(60,17){\vector(0,-1){8}}

\put(60,20){\makebox(0,0){\mbox{\small$\mu(T; 
U_n, ... , U_n) 
$}}}

\put(104,20){\makebox(0,0){\mbox{\small$\mu(T; 
e,...,e) 
$}}}

\put(60,5){\makebox(0,0){\mbox{\small$T 
$}}}

\put(95,15){\vector(-4,-1){30}}

\put(85,9){\makebox(0,0){\mbox{\small$id
$}}}

\put(118,20){\makebox(0,0){\mbox{\small $   
$}}}

\end{picture}}

\noindent and

{\unitlength=1mm
\begin{picture}(50,25)(28,2)

\put(87,20){\vector(-1,0){15}}

\put(60,17){\vector(0,-1){8}}

\put(60,20){\makebox(0,0){\mbox{\small$\mu (U_n; 
T)
$}}}

\put(98,20){\makebox(0,0){\mbox{\small$\mu(e ;
T)
$}}}

\put(60,5){\makebox(0,0){\mbox{\small$T 
$}}}

\put(95,17){\vector(-3,-1){25}}

\put(84,11){\makebox(0,0){\mbox{\small$id
$}}}

\end{picture}}

\noindent as relations.

This operad contains an internal $n$-operad given by 
$a_T = T$.

We can construct $\H^n$ for all $n$ including $n=\infty$ in the same fashion.

To better understand the structure of $\h^n ,$ we can describe it in terms of decorated planar trees.

 An object
of $\h^n$ is a labelled planar tree with vertices decorated by trees from $Tr_n$ in the following sense:
to every vertex $v$  of valency $k$ we associate an $n$-tree with $k$-tips. The following picture
illustrates the concept for $n=2$. 

{\epsfxsize=200pt 
\makebox(260,180)[r]{\epsfbox{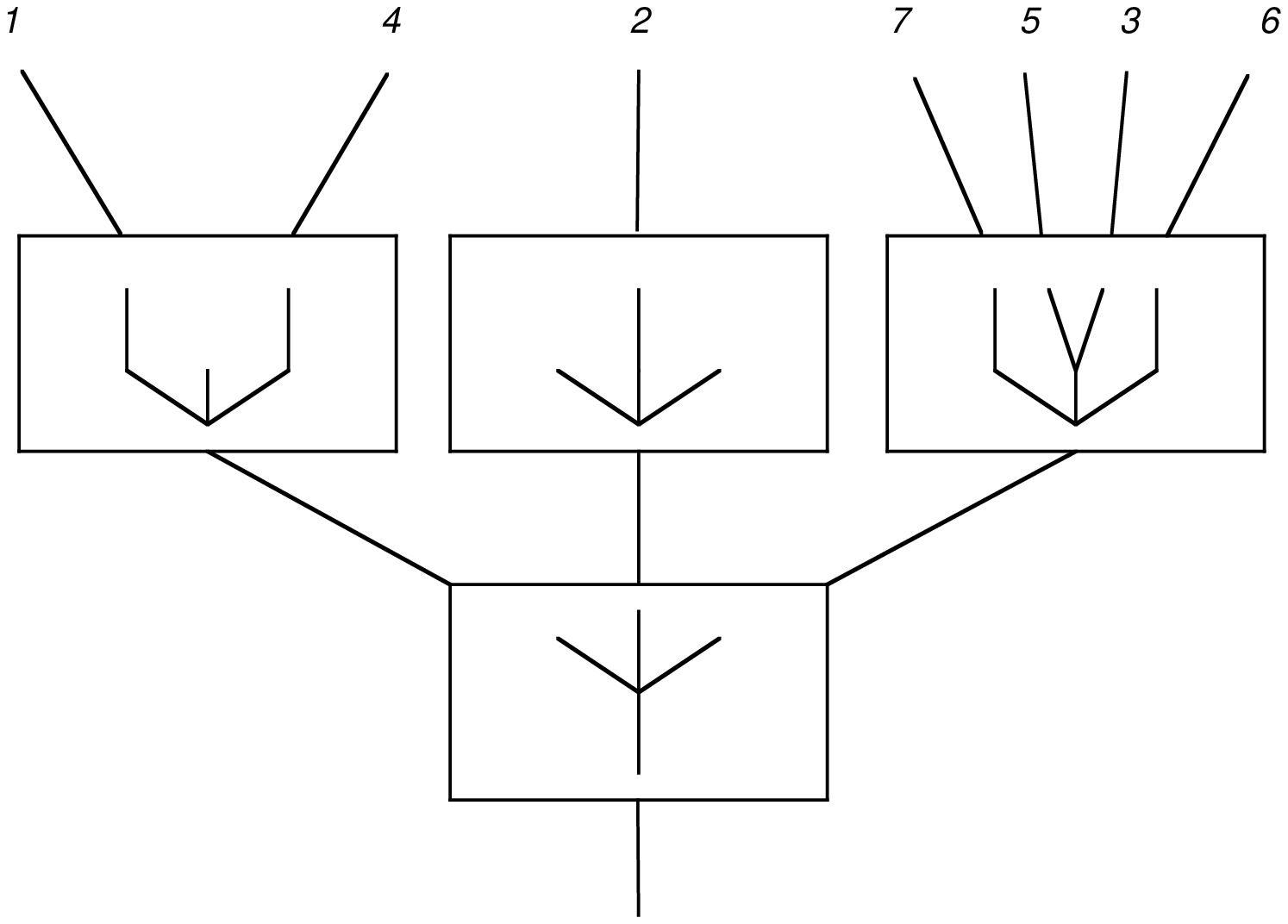}}}

So,  the objects of $\h^n$ are labelled planar trees with some extra internal structure. 
The morphisms are contractions or growing of internal edges, yet not all contractions are possible. It
depends on the extra internal structure. We can simultaneously contract the input edges of a vertex $v$ only 
if the corresponding $n$-trees in the vertices above $v$ can be pasted together in the 
$n$-category $Tr_n$
 according to the
globular pasting scheme determined by the tree at the vertex $v$. In the above example we see that the trees on 
the highest level
are fibers  over a map of trees:

 {\epsfxsize=140pt 
\makebox(220,60)[r]{\epsfbox{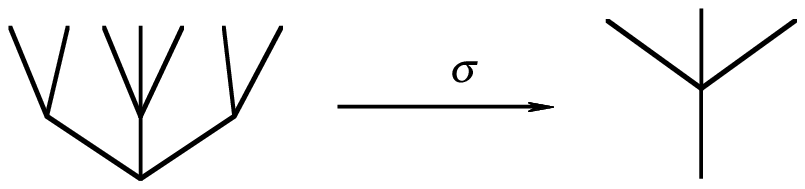}}}

So in $\h^2$ we have a morphism corresponding to the $\sigma$ from the object above to the object

 {\epsfxsize=80pt 
\makebox(200,130)[r]{\epsfbox{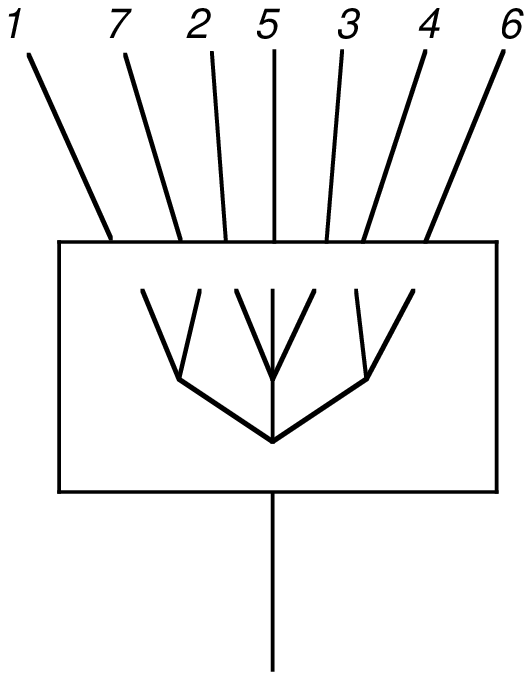}}}

The  case  $n=1$ is  well known.

Indeed, with $n=1$ all the decorations are meaningless. Yet, the morphisms in $\h^1$ correspond
only to order-preserving maps between ordinals. 
 
 Therefore the operad $\h^{1}$  coincides with the symmetrisation of the nonsymmetric operad $h$  described in
\cite{BatA}, which is, indeed, $\H^1$ in our present terminology.  For a discussion on it the reader may also look at  
\cite{DS}. The objects of $\H^{1}_{k} $  are bracketings of the strings 
consisting of several $0$'s and   symbols $1,\ldots,k$ in fixed order without repetition. Multiple
bracketing like
$(((\ldots)))$ and also empty bracketing $( \ )$ are allowable. The morphisms are throwing off $0$'s, removing and introducing
 a pair of brackets, and also a morphism $( \ ) \rightarrow (1)$. The symmetric groups act by permuting the
symbols 
$1,\ldots,k$. The operad multiplication is given by replacing  one of the symbols by a corresponding expression. 

It is
clear that 
$\pi_{0}(\h^1_k)=\Sigma_k$ and all higher homotopy groups vanish. In other words $\h^1$ is an
$A_{\infty}$-operad. The algebras of $\h^1$ in $Cat$ are categories equipped with an $n$-fold tensor product 
satisfying some obvious associativity and unitarity conditions. For example, instead of a single associativity
isomorphism
 we will
have  two, perhaps  noninvertible, morphisms from two different combinations of binary products
to the triple tensor product; that is, a cospan 
$$(a\otimes b)\otimes c \longrightarrow a\otimes b\otimes c \longleftarrow a\otimes (b \otimes c).$$ 
Instead of the pentagon, we will have a barycentric subdivision of it and so on. Such categories were called
lax-monoidal in
\cite{DS}.
 
The operad  $\H^{\infty}$ is also classical. All the decorations again collapse to a point. But  morphisms are
more complicated and correspond to the maps of finite sets. So  we can give the following 
 description of  the
operad
$\H^{\infty}$. A typical object of
$\H^{\infty}_n$  is a planar tree  equipped with an
 injective function (labelling) from $[n]$ to the set of vertices of this  tree. The symmetric group
acts by permuting the labels. The morphisms are generated by 
 contraction of an internal edge, growing of an internal edge, and dropping unlabelled leaves,
with usual relations of associativity and unitarity.
 We also will have an isomorphism $T\rightarrow \pi T$ for every
permutation $\pi \in
S_n$. This isomorphism should satisfy obvious equivariancy conditions. Again the $\H^{\infty}$-algebras in $Cat$ are symmetric
lax-monoidal categories in the terminology of \cite{DS1}.

The internal operad is given by the trees  

{\epsfxsize=90pt 

\makebox(200,90)[r]{\epsfbox{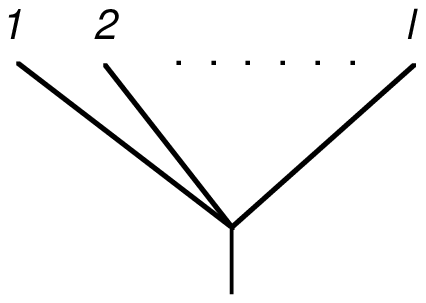}}}

\noindent and this is a terminal object in $\H^{\infty}$. Hence, the nerve of
$\H^{\infty}$  is an $E_{\infty}$-operad.     

\

\Remark The trees formalism from
\cite{GK}[Section 1.2] is actually
 a special case of our theorem \ref{rep} with $n= \infty$ and $ \ A =
 V^{\bullet}$ for a symmetric monoidal category $V$.  

\

To clarify  the structure of $\h^n$ for $2\le n < \infty$ we provide a part of
 the picture of $\h^2_2$:

{\epsfxsize=300pt 

\makebox(320,250)[r]{\epsfbox{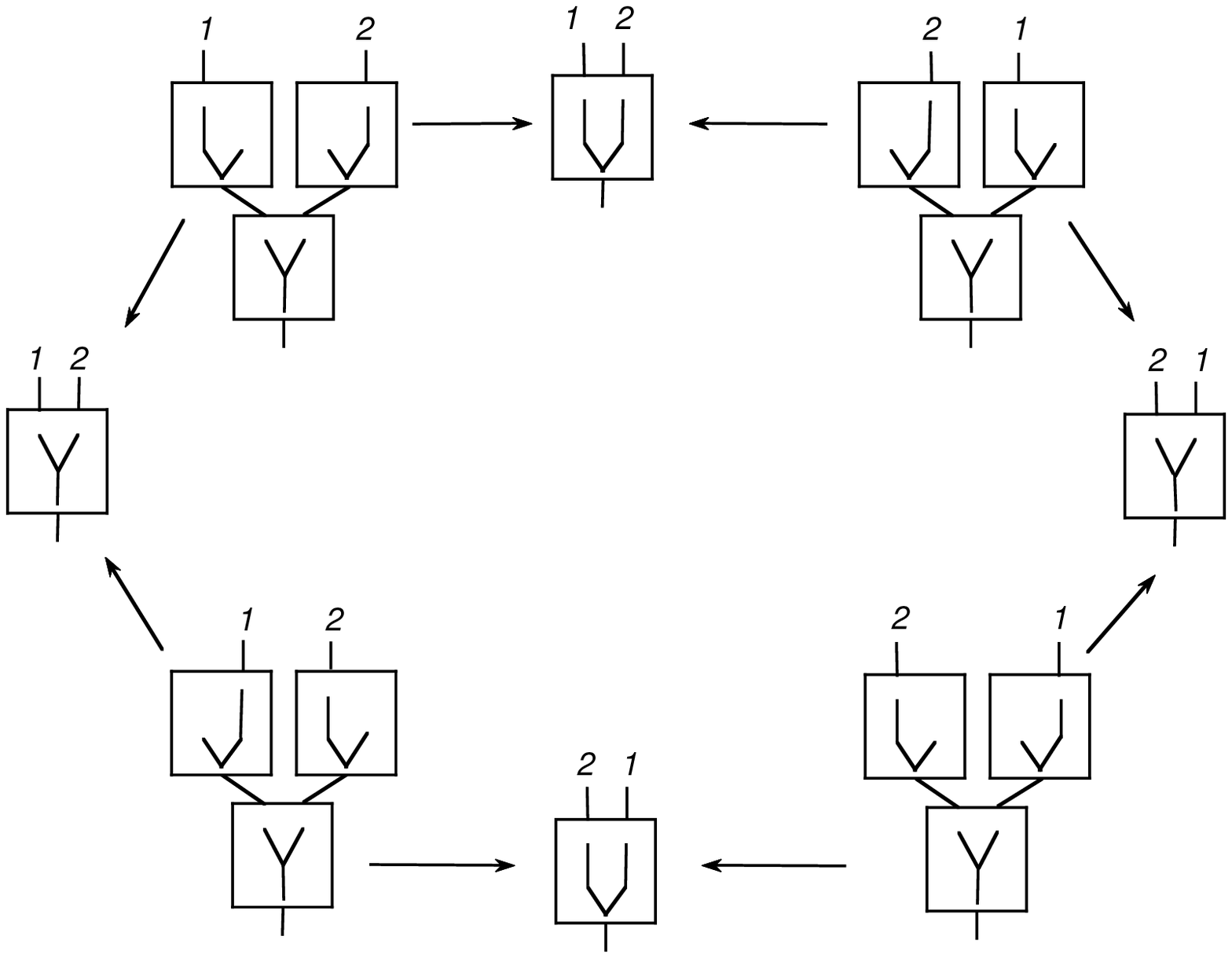}}}

  The reader  can find an analogy
with a diagram from the construction  of the braiding in Proposition 5.3 from \cite{JS}. 
The reader may also look at a similar picture
 for the category   $\M^2_2$  where $\M^2$ is a $Cat$-operad constructed in
\cite{F}. We also recommend  the reader  look at the picture of $\M^3_2$ in \cite{F}, which
looks like a two dimensional sphere, and   try to construct a similar picture in $\h^3_2$. 
 Of course, these are not  
 accidental coincidences as we will show in the next paper \cite{BFM}.

\section{Example: iterated monoidal category operad}\label{IMon}

First of all we briefly review the construction of the iterated monoidal category operad $\ML^n$ introduced in \cite{F}. 

The objects  of $\ML^n_k$ are all finite expressions generated by the symbols $1,\ldots,k$ and 
$n$ associative operations $\ten{1},\ldots,\ten{n}$ in which each generating symbol occurs exactly once.
 There is a natural left action 
of symmetric group on $\ML^n_k$ and an   operation of substitution which provides an
operadic structure on the objects of $\ML^n$. 

Now we can describe the 
 morphisms in $\ML^n$. They are generated by 
the middle interchange laws
\begin{equation}\label{interchange}\eta^{ij}:(1\ten{i}2)\ten{j}(3\ten{i}4) \rightarrow 
(1\ten{j}3)\ten{i}(2\ten{j}4) , \ j< i,\end{equation}
 substitutions and permutations, and must satisfy the coherence conditions specified 
in the first section of \cite{F}. It was shown in \cite{F} that the  operad $\ML^n$ is a poset
operad. The algebras of $\ML^n$ in $Cat$ are iterated
$n$-monoidal categories, i.e. categories with $n$  strict
monoidal structures which are related by interchange morphisms
(not necessary isomorphisms) satisfying some natural
coherence conditions. They are also monoids in the category of iterated $(n-1)$-monoidal categories and they {\it lax} monoidal functors.

We also would like to introduce another categorical symmetric operad $\M^n$ which is constructed in the same way as $\ML^n$ but 
we use operations $\ten{0},\ldots,\ten{n-1}$ and we reverse the direction of the interchange law (\ref{interchange}.

This is  the   picture of $\M^2_2$ .

{\unitlength=1mm

\begin{picture}(60,33)

\put(55,27){\makebox(0,0){\mbox{$ 1\ten{0}2 $}}}
\put(55.5,4){\makebox(0,0){\mbox{$ 2\ten{0}1$}}}

\put(40,15){\makebox(0,0){\mbox{$ 1\ten{1}2 $}}}
\put(70,15){\makebox(0,0){\mbox{$ 2\ten{1}1 $}}}

\put(45,18){\vector(1,1){6}}
\put(45,12){\vector(1,-1){6}}

\put(65,18){\vector(-1,1){6}}
\put(65,12){\vector(-1,-1){6}}

\put(44.5,21.5){\shortstack{\mbox{$\scriptstyle
\eta^{01}
      $}}}

\put(63,21.5){\shortstack{\mbox{$\scriptstyle
\eta^{01}
      $}}}

\put(42.5,8){\shortstack{\mbox{$\scriptstyle
\eta^{01}
      $}}}

\put(63,8){\shortstack{\mbox{$\scriptstyle
\eta^{01}
      $}}}

\end{picture}}

There is an obvious  isomorphism of operads $\ML^n$ and $(\M^n)^{op} .$ So the algebras of $\M^n$ are  monoids in the category of iterated $(n-1)$-monoidal categories and their {\it oplax} monoidal functors. 
  Of course, these two operads have the same  homotopy type. 
We consider here the operad $\M^n$ simply because it is better adapted to our agreement about directions of middle interchange cells and numeration of operations, which makes our proof easier to follow. 

\begin{theorem} \label{M} The categorical symmetric  operad ${\ML}^n$ contains both an internal $n$-operad and internal $n$-cooperad. The same is true for the operad $\M^n .$ \end{theorem}

\

\Remark We did not discuss the notion of internal $n$-cooperad but it can be easily obtained from the definition of internal $n$-operad by inverting the structure cells.  

\

\noindent {\bf Proof.} We will give a proof that $\M^n$ contains an internal $n$-operad. The other statements of the theorem follow. It is sufficient to change the numeration and reverse the direction of morphisms in an appropriate way. 

We have to assign an object $a_T\in {\M}^n_k$ to every $n$-tree $T$
with $|T|=k$.  We will do it by induction.  We put $a_T=1$ for all trees with
$|T|=0,1.$ In particular, $a_{U_n}= 1$.  Now, suppose we have already constructed
$a_T$ for  all trees which are $(n-k)$-fold suspensions.
  Suppose  a tree $T$ is an $(n-k-1)$-fold suspension. 
Take a canonical decomposition $$T= T_1\otimes_{n-k-1} T_2\otimes_{n-k-1} ...
\otimes_{n-k-1} T_r$$ 
 Then we put
$$a_T = m(1\diamond_{\scriptscriptstyle n-k-1}2\diamond_{\scriptscriptstyle n-k-1} ... \diamond_{\scriptscriptstyle n-k-1} r
; a_{T_1}, a_{T_2}, ... ,
 a_{T_r}),$$ 
where $m$ is the multiplication in ${\M}^n$.

\

\Example 
To give an idea how the operad multiplication in $a$ looks  we present the following $2$-dimensional example.
 \vskip15pt \hskip40pt \includegraphics[width=160pt]{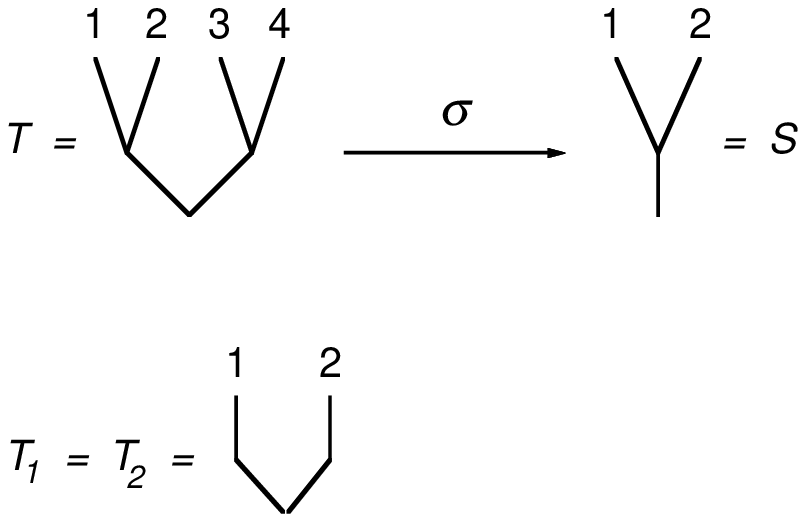}\vskip15pt

In this picture the map of trees is given by
$$\sigma(1)=1 , \sigma(2)=2, \sigma(3)=1 , \sigma(4)=2 ,$$
$$\pi(\sigma) = (1324).$$

Then $$a_T = (1\diamond_{\scriptscriptstyle 1} 2)\diamond_{\scriptscriptstyle 0} (3\diamond_{\scriptscriptstyle 1}  4). $$
$$m(a_S;a_{T_1},a_{T_2}) = m(1\diamond_{\scriptscriptstyle 1} 2; 1\diamond_{\scriptscriptstyle 0}
2, 1\diamond_{\scriptscriptstyle 0} 2) = (1\diamond_{\scriptscriptstyle 0} 2)\diamond_{\scriptscriptstyle 1}
(3\diamond_{\scriptscriptstyle 0} 4),$$ 
and the operadic multiplication  $\mu_{\sigma}$ is given by the
middle interchange morphism
$$\eta_{1,2,3,4}:(1\diamond_{\scriptscriptstyle 0} 2)\diamond_{\scriptscriptstyle 1} (3\diamond_{\scriptscriptstyle 0}
4) \longrightarrow (1\diamond_{\scriptscriptstyle 1} 3)\diamond_{\scriptscriptstyle 0}
(2\diamond_{\scriptscriptstyle 1} 4).$$

Before we construct the multiplication in general we have to formulate the following lemma  whose proof is obtained by
an obvious induction.
\begin{lem} \label{pass} Let the $n$-tree $T$ be
$$T=
[k_n]\stackrel{\rho_{n-1}}{\longrightarrow}[k_{n-1}]\stackrel{\rho_{n-2}}{\longrightarrow}...
\stackrel{\rho_0}{\longrightarrow} [1]
$$  
then an element
$u
\diamond_{\scriptscriptstyle i} v$ is in
$a_T$ in the sense of {\rm \/ \cite{F}}
 if and only if $u<v$ and  $$ \rho_{n-1}\cdot ... \cdot \rho_{i}(u) = \rho_{n-1}\cdot
...
\cdot \rho_{i}(v)$$ but 
$$ \rho_{n-1}\cdot ... \cdot \rho_{i+1}(u) \ne \rho_{n-1}\cdot ...
\cdot \rho_{i+1}(v).$$
\end{lem}

Now we want to construct the multiplication $m_{\sigma}$ in the special case
where $$\sigma: T\longrightarrow M_k^2 . $$
So we have to construct a morphism
$$m(1\diamond_{\scriptscriptstyle k} 2 ; a_{T_1}, a_{T_2})\longrightarrow \pi(\sigma) a_T.$$
According to \cite{F} we have to check that  $u\diamond_{\scriptscriptstyle i} v$ in $m(1\diamond_{\scriptscriptstyle k} 2 ;
a_{T_1}, a_{T_2})$ implies  either $u\diamond_j v$ in $\pi(\sigma)a_T$ for $j\le i$ or 
$v\diamond_{\scriptscriptstyle j} u$ in
$\pi(\sigma)a_T$ for
$j< i$.

Recall that $m(1\diamond_{\scriptscriptstyle k} 2 ;
a_{T_1}, a_{T_2}) = a_{T_1}\diamond_{\scriptscriptstyle k} \bar{a_{T_2}} $ where $\bar{a_{T_2}}$ is
the same expression as $a_{T_2}$ but all numbers are shifted on $|T_1|$.  Let $\xi_i:
T_i\rightarrow T, i=1,2$ be inclusions of $T_i$ as $i$-th fiber.  

 Now, suppose 
$u\diamond_{\scriptscriptstyle i} v$ is in
$a_{T_1}.$ By our lemma it means that $u<v$ and  $$ \rho_{n-1}\cdot ... \cdot \rho_{i}(u) = \rho_{n-1}\cdot ...
\cdot \rho_{i}(v)$$ but 
$$ \rho_{n-1}\cdot ... \cdot \rho_{i+1}(u) \ne \rho_{n-1}\cdot ...
\cdot \rho_{i+1}(v)$$ in $T_1$. Hence, we have
$\xi_1(u)<\xi_1(v)$ and  $$ \rho_{n-1}\cdot ... \cdot \rho_{i}(\xi_1(u)) = \rho_{n-1}\cdot ...
\cdot \rho_{i}(\xi_1(v))$$ but 
$$ \rho_{n-1}\cdot ... \cdot \rho_{i+1}(\xi_1(u)) \ne \rho_{n-1}\cdot ...
\cdot \rho_{i+1}(\xi_(v))$$ in $T$. But $\pi(\sigma)(\xi_1(w))=w$ by definition of
$\pi(\sigma)$. Therefore,
$\pi(\sigma)^{-1}u < \pi(\sigma)^{-1}v$ and  $$ \rho_{n-1}\cdot ... \cdot \rho_{i}(\pi(\sigma)^{-1}u) =
\rho_{n-1}\cdot ...
\cdot \rho_{i}(\pi(\sigma)^{-1}v)$$ but 
$$ \rho_{n-1}\cdot ... \cdot \rho_{i+1}(\pi(\sigma)^{-1}u) \ne \rho_{n-1}\cdot ...
\cdot \rho_{i+1}(\pi(\sigma)^{-1}v)$$ in $T$. By our lemma it follows that $u\diamond_i v $ is in $a_T$.

The same argument applies if $u\diamond_{\scriptscriptstyle i} v$ is in $\bar{a_{T_2}}$ but all numbers must be shifted on
$|T_1|$.

Now suppose $u$ is in $a_{T_1}$ but $v$ is in $\bar{a_{T_2}}$. This means that $u\diamond_{\scriptscriptstyle k} v$ is
in 
$a_{T_1}\diamond_{\scriptscriptstyle k} \bar{a_{T_2}}.$  

 We have two possibilities. The first is 
$$ \rho_{n-1}\cdot ... \cdot \rho_{k}(u) = \rho_{n-1}\cdot ...
\cdot \rho_{k}(v)$$
where the first composite is in $T_1$ and the second is in $T_2$. This means that 
$$ \rho_{n-1}\cdot ... \cdot \rho_{k}(\xi_1(u)) = \rho_{n-1}\cdot ...
\cdot \rho_{k}(\xi_2(v))$$
already in $T$.  But $\sigma$ is a morphism of trees, hence, preserves order on fibers of $\rho_k$ and we
have $\xi_1(u)<\xi_2(v)$, hence, again $\pi(\sigma)^{-1}u < \pi(\sigma)^{-1}v$ and 
$$ \rho_{n-1}\cdot ... \cdot \rho_{k}(\pi(\sigma)^{-1}u) =
\rho_{n-1}\cdot ...
\cdot \rho_{k}(\pi(\sigma)^{-1}v)$$
and therefore $u\diamond_{\scriptscriptstyle k} v$ is in $\pi(\sigma)a_T$.

The last possibility is  
$$ \rho_{n-1}\cdot \ldots \cdot \rho_{l}(u) = \rho_{n-1}\cdot\ldots
\cdot \rho_{l}(v)$$
for some $l<k$ but $$ \rho_{n-1}\cdot \ldots \cdot \rho_{l+1}(u) \ne \rho_{n-1}\cdot \ldots
\cdot \rho_{l+1}(v)$$ again in $T_1$ and $T_2$ respectively. Then 
$$ \rho_{n-1}\cdot \ldots \cdot \rho_{l}(\xi_1(u)) = \rho_{n-1}\cdot \ldots
\cdot \rho_{l}(\xi_2(v))$$
for some $l<k ;$ but $$ \rho_{n-1}\cdot \ldots \cdot \rho_{l+1}(\xi_1(u)) \ne \rho_{n-1}\cdot \ldots
\cdot \rho_{l+1}(\xi_2(v))$$
already in $T$. By the usual argument it follows that either $u\diamond_{\scriptscriptstyle l} v$ or
$v\diamond_{\scriptscriptstyle l} u$ is in
$\pi(\sigma)a_T$, and  that finishes the proof of the special case. 

Now, suppose we have constructed $m_{\sigma}$ for all $\sigma$ whose codomain is
$S = M_k^j$ and where  $j\le m$. 
 Then, for $S = M_k^{m+1},$   $$a_S = m(1\diamond_{\scriptscriptstyle k} 2;a_{S_1},a_{S_2})$$ and an easy
inductive argument can be
applied.

In general, let $\sigma:T\rightarrow S$ be a morphism of trees. If $S=U_n$ then we put
$\mu_{\sigma}= id$. Now suppose we already have constructed $\mu_{\sigma}$ for all $\sigma$
with codomain being an $(n-k)$-fold suspension. Let $S$ be an $(n-k-1)$-fold suspension. Then
the canonical decomposition of $S$  gives us  
$$\omega: S \rightarrow M^j_{n-k-1}$$
 with
$(n-k)$-fold suspensions $S_i, \ 1\le i \le r $, as fibers. We have 
$$m(a_S;a_{T_1},\ldots,a_{T_k})= m(m(1\diamond_{\scriptscriptstyle n-k-1} ... \diamond_{\scriptscriptstyle n-k-1} r
; a_{S_1}, 
\ldots ,
 a_{S_r}),a_{T_1},\ldots ,a_{T_k}) =$$ $$= m(1\diamond_{\scriptscriptstyle n-k-1} ... \diamond_{\scriptscriptstyle
n-k-1} r;m(a_{S_1};a_{T_1^1},\ldots ,a_{T_1^{m_1}}),\ldots,
m(a_{S_r};a_{T_r^{1}},\ldots,a_{T_r^{m_r}})).$$  By the inductive hypothesis we already have
$m_{\sigma_i}$ for the fibers of $\sigma$. So we have a
 morphism
$$m(1,m_{\sigma_1},\ldots,m_{\sigma_r}): m(1\diamond_{\scriptscriptstyle n-k-1} 
... \diamond_{\scriptscriptstyle n-k-1}
r;m(a_{S_1};a_{T_1^1},\ldots ,a_{T_1^{m_1}}),\ldots $$ $$\ldots,
m(a_{S_r};a_{T_r^{1}},\ldots, a_{T_r^{m_r}}))
\longrightarrow m(1\diamond_{\scriptscriptstyle n-k-1} ... \diamond_{\scriptscriptstyle n-k-1}
r;\pi(\sigma_1)a_{T'_1},\ldots ,\pi(\sigma_r)a_{T'_r}),$$
where $T'_1,\ldots ,T'_r$ are fibers of $\sigma\cdot\omega$. 
But
$$m(1\diamond_{\scriptscriptstyle n-k-1} ... \diamond_{\scriptscriptstyle n-k-1}
r;\pi(\sigma_1)a_{T'_1},\ldots,\pi(\sigma_r)a_{T'_r}) = $$ $$ =
\Gamma(1,\pi{\sigma_1},\ldots,\pi{\sigma_r}) m(1\diamond_{\scriptscriptstyle n-k-1} ... \diamond_{\scriptscriptstyle
n-k-1} r;a_{T'_1},\ldots,a_{T'_r}).$$
Now  we  already have the morphism
$$m_{\sigma\cdot\omega}:m(1\diamond_{\scriptscriptstyle n-k-1} ... \diamond_{\scriptscriptstyle n-k-1}
r;a_{T'_1},\ldots,a_{T'_r}) \rightarrow \pi(\sigma\cdot\omega)a_T .$$
So we have
$$\Gamma(1,\pi{\sigma_1},\ldots,\pi{\sigma_r})m_{\sigma\cdot\omega}:
\Gamma(1,\pi{\sigma_1},\ldots,\pi{\sigma_r})
m(1\diamond_{\scriptscriptstyle n-k-1} ... \diamond_{\scriptscriptstyle n-k-1}
r;a_{T'_1},\ldots,a_{T'_r})\rightarrow$$ $$\longrightarrow
\Gamma(1,\pi{\sigma_1},\ldots,\pi{\sigma_r})\pi(\sigma\cdot\omega)a_T.$$ By Lemma
\ref{pisigma}, $$\Gamma(1,\pi{\sigma_1},\ldots,\pi{\sigma_r})\pi(\sigma\cdot\omega) =
\Gamma(\pi(\omega),1,\ldots,1)\pi(\sigma).$$
But $\omega$ is order preserving, hence, the last permutation is $\pi(\sigma)$.
So the composite 
  $$m(1,m_{\sigma_1},\ldots,m_{\sigma_r})\cdot
\Gamma(1,\pi{\sigma_1},\ldots,\pi{\sigma_r})m_{\sigma\cdot\omega}$$
gives us the required morphism
$$\mu_{\sigma}: m(a_S;a_{T_1},...,a_{T_k}) \longrightarrow \pi(\sigma)a_T.$$
Associativity and unitarity of this multiplication are trivial because $\M^n$ is a
poset operad.

\Q

\

 In \cite{F}, a morphism of categorical  operads
$$\ML^n \longrightarrow {\cal K}^{(n)} $$
is constructed. Here ${\cal K}^{(n)}$ is the $n$-th filtration of Berger's complete graph operad \cite{Berger}, which plays a
central role in his theory of cellular operads. So we have 
\begin{cor} \label{K} ${\cal K}^{(n)}$  contains
 an internal $n$-operad and an internal $n$-cooperad. \end{cor}

\section{\label{freeoperad}Free internal operads }

In this section we apply the techniques described in   Theorem \ref{rep} to get some  formulas which will be of use in the final section as well as in \cite{BFM}.

\begin{defin} We call a categorical  $n$-operad {\it cocomplete} if each category $A_T$ is cocomplete and multiplication 
in $A$ preserves colimits in each variable. We give a similar definition of cocompleteness  for symmetric operads.\end{defin}

An internal object in an $n$-operad $A$ (see example on page  \pageref{internalobjects}) will be called an {\it internal $n$-collection in $A$ .} So we have a category of internal $n$-collections $IColl_n(A) $ 
and the corresponding categorical operad $\H^n_d$ which represents this $2$-functor.

\

\Example Let $A = Des_n( V^{\bullet})$ for a symmetric monoidal category $V$. Then the category  $IColl_n(A)$ is isomorphic to the category $Coll_n(V)$ of
$n$-collections in $V .$ 

\

Given an internal $n$-collection $x$ in $A$ we
will denote by
$\tilde{x}:
\H^n_d\rightarrow A$  the corresponding  operadic functor.

\begin{theorem} \label{freeop} Let $A$ be a cocomplete categorical  $n$-operad.

The free internal  n-operad on an $n$-collection $x$ is given by the 
formula
$${\cal F}_n(x)_T = \coprod_{\scriptstyle W\in {\scriptstyle {\bf H}^n_T}} \tilde{x}(W).$$
More generally, the $k$-th iteration of ${\cal F}$ is given by  the formula
$${\cal F}_n^k(x)_T = \coprod_{W_1\stackrel{{\scriptstyle f}_{1}}{\leftarrow\!-}W_2\stackrel{\scriptstyle f_2}{\leftarrow\!-}\
\ldots\stackrel{\scriptstyle f_{\scriptscriptstyle k-1}}{\leftarrow\!-}\ W_k}
\tilde{x}(W_k),$$
where $f_1,\ldots,f_{k-1}$ are morphisms in $\H^n_T$. 

\end{theorem}

\noindent {\bf Proof.}
  The  left Kan extension in the $2$-category 
of categorical $n$-collections
of $\tilde{x}$ along the inclusion $i: \H^n_d \rightarrow \H^n$ is given by the following 
formula 
\begin{equation}\label{ext}  Lan_i(\tilde{x})(W) = \coprod_{\scriptstyle W\leftarrow W'} \tilde{x}(W').\end{equation} 

  We are going to prove that it is also a left Kan extension in ${\bf CO}_n$.
We have thus to show  that the functor $Lan_i$ is operadic.

 Indeed, let $\sigma: T \rightarrow S$ be a morphism of trees and let
$W_S \in \H^n_{S},W_1 \in \H^n_{T_1},\ldots, W_k \in \H^n_{T_k}.$  Then 
$$\mu_A(Lan_i(\tilde{x})(W_S);Lan_i(\tilde{x})(W_1),\ldots ,Lan_i(\tilde{x})(W_k)) \simeq $$
$$\simeq \coprod_{\scriptstyle W_S\leftarrow W'_S,W_1\leftarrow W'_1,\ldots,W_k\leftarrow W'_k}\tilde{x}(W'_S)\times
\tilde{x}(W'_1)\times
\ldots
\tilde{x}(W'_k) = \ $$ 
$$= \coprod_{\scriptstyle \mu(W_S;W_1,\ldots,W_k) \leftarrow \mu(W'_S;W'_1,\ldots,W'_k)} 
\tilde{x}(\mu(W'_S;W'_1,\ldots,W'_k)) \simeq \ \ \ \ \ 
   $$
$$\simeq \coprod_{\scriptstyle \mu(W_S;W_1,\ldots,W_k) \leftarrow W'}\tilde{x}(W') =
Lan_i(\tilde{x})(\mu(W_S;W_1,\ldots,W_k)),$$  since $\tilde{x}$ is an  operadic functor and by the inductive construction of
objects and morphisms in $\H^n$.  Analogously one can prove that the counit of this adjunction is operadic. It is
straightforward now to show that this is really an operadic adjunction. 

Let us denote by ${\bf Lan}$  the monad generated by the adjunction $Lan_i\dashv i^{\star}$. Then  from the
formula (\ref{ext}) we have the following formula for the iteration of this monad: 
\begin{equation}\label{ext1}{\bf Lan}^k(\tilde{x})(W) = \coprod_{W\stackrel{{\scriptstyle
f}_{0}}{\leftarrow\!-}W_1\stackrel{{\scriptstyle f}_{1}}{\leftarrow\!-}W_2\stackrel{\scriptstyle f_2}{\leftarrow\!-}\
\ldots\stackrel{\scriptstyle f_{\scriptscriptstyle k-1}}{\leftarrow\!-}\ W_k}
\tilde{x}(W_k),\end{equation}
where $f_0,\ldots,f_{k-1}$ are morphisms in $\H^n$.

To obtain the formula for the free operad it is enough to evaluate the formula (\ref{ext}) at $T$. 
Since $T$ is the terminal object in $\H^n_T $ we get the formula as in the statement of the theorem. Analogously one obtains the formula
for the iterated free operad monad. 

\

\Q

\

\Remark We will encounter a similar situation with the calculation of a Kan extension in Theorem \ref{formula}.
 
\

The analogous result holds in the case of symmetric operads.

\begin{theorem} Let $A$ be a cocomplete symmetric categorical operad.

The free internal symmetric operad on a nonsymmetric internal collection $x$ is given by the 
formula
$${\cal F}_{\infty}(x)_m = \coprod_{\scriptstyle W\in {\scriptstyle {\bf H}^{\infty}_m}} \tilde{x}(W).$$
More generally, the $k$-th iteration of ${\cal F}_{\infty}$ is given by  the formula
$${\cal F}_{\infty}^k(x)_m = \coprod_{W_1\stackrel{{\scriptstyle f}_{1}}{\leftarrow\!-}W_2\stackrel{\scriptstyle
f_2}{\leftarrow\!-}\
\ldots\stackrel{\scriptstyle f_{\scriptscriptstyle k-1}}{\leftarrow\!-}\ W_k}
\tilde{x}(W_k),$$
where $f_1,\ldots,f_{k-1}$ are morphisms in $\H^{\infty}_m$. 

\end{theorem}

 \section{Colimit formula for symmetrisation}\label{suspension}
 
 Now we return to the study of the canonical operadic functor
 $$\zeta: \h^{n} \rightarrow \H^{\infty}$$

\begin{lem} For $n\ge 2$ the functor $\zeta$ is final (in the sense of \/ {\rm \cite{ML}} ). \end{lem}

\noindent {\bf Proof.} The functor $\zeta$ is surjective on objects by construction. Hence, it will be sufficient
to prove that,  for any morphism $f:a\rightarrow b$ in $\H^{\infty}$ and any objects $a',b'\in \h^n$ such that
$\zeta(a')=a$ and $\zeta(b')=b,$ there exist a chain of morphisms in $\h^n$
$$b'\leftarrow x_1\rightarrow \ldots \leftarrow x_{i+1} \stackrel{f'}{\rightarrow} x_{i} \leftarrow \ldots
  \leftarrow x_m
\rightarrow a'$$ with the following properties:

- there exists an $0\le i \le m $ such that $\zeta(f') = f ;$ 

- the image  under $\zeta$ of any other arrow
is either a retraction or its right inverse;  

- the image under $\zeta$ of a composite of the appropriate morphisms or their inverses gives an identity $\zeta
x_i\rightarrow a$;

-  the image under $\zeta$ of a composite of the appropriate morphisms or their inverses gives an identity $\zeta
x_{i+1}\rightarrow b$.

If these all are  the case then   the following commutative diagram
provides a path between  any two objects in the comma-category of
$\zeta$ under the object $a$ from $\H^{\infty}$.

{\unitlength=1mm

\begin{picture}(60,33)(-9,0)

\put(45,25){\makebox(0,0){\mbox{$a$}}}
\put(45,21){\vector(0,-1){10}}
\put(55,26){\shortstack{\mbox{$g $}}}

\put(70,25){\makebox(0,0){\mbox{$\zeta(c)$}}}
\put(41,25){\vector(-1,0){16}}
\put(33,26){\shortstack{\mbox{$f $}}}
\put(20,25){\makebox(0,0){\mbox{$\zeta(b)$}}}
\put(49,25){\vector(1,0){16}}

\put(41,23){\vector(-3,-1){16}}
\put(41.5,21){\vector(-1,-1){10}}

\put(21,16){\makebox(0,0){\mbox{$\zeta(x_i)$}}}
\put(30,20.5){\shortstack{\mbox{\small $f $}}}

\put(29,8){\makebox(0,0){\mbox{$\zeta(x_{i+1})$}}}
\put(32,15.5){\shortstack{\mbox{\small $id $}}}
\put(20,18.5){\vector(0,1){4}}
\put(15.5,19){\shortstack{\mbox{\small $id $}}}

\put(27,10.5){\vector(-2,1){5}}
\put(26,12.5){\shortstack{\mbox{\small $f $}}}

\put(40,8){\vector(-1,0){4.5}}
\put(37,9){\shortstack{\mbox{\small $id $}}}

\put(41,15){\shortstack{\mbox{$id $}}}

\put(45,8){\makebox(0,0){\mbox{$\zeta(a') $}}}

\put(49,23){\vector(3,-1){16}}
\put(48.5,21){\vector(1,-1){10}}

\put(70,16){\makebox(0,0){\mbox{$\zeta(y_j)$}}}
\put(58.5,20.5){\shortstack{\mbox{\small $g $}}}

\put(64,8){\makebox(0,0){\mbox{$\zeta(y_{j+1})$}}}
\put(55,15.5){\shortstack{\mbox{\small $id $}}}
\put(70,18.5){\vector(0,1){4}}
\put(71.2,19){\shortstack{\mbox{\small $id $}}}

\put(62,10.5){\vector(2,1){5}}
\put(61,12.5){\shortstack{\mbox{\small $g $}}}

\put(50,8){\vector(1,0){5.5}}
\put(51,9){\shortstack{\mbox{\small $id $}}}

\put(41,15){\shortstack{\mbox{$id $}}}

\end{picture}}

It is clear that it will be enough to show that the above property
holds for a generating morphism 
$$f: \mu([k];[n_1],\ldots,[n_k]) \longrightarrow \pi(\sigma)[m]$$
which corresponds to the morphism of ordinals
$\sigma:[m]\rightarrow [k]$.

Let the trees $T$, $S$ and $T_1,\ldots , T_k$ be such that
$$\zeta(T)= [m] \ ,$$
$$\zeta(S)= [k] \ ,$$
$$\zeta(T_i)= [m_i], \ 1\le i\le k .$$
Then $$\zeta(\mu(S,T_1,\ldots,T_k))  =
\mu([k],[m_1],\ldots,[m_k]) .$$
Let  $$T' = {M_0^m} \  \  \ ,  \ \ \ S' =
M_{n-1}^k \mbox{\ \ \ \ and \ \  \  }
 T'_i = {M_0^{m_i}} .$$

Then $\sigma$ determines a unique morphism $\sigma':T'\rightarrow S'$ in
$\Omega_n$ with $\sigma'_n = \sigma$.  This morphism  gives 
the following morphism in $\h^n :$
$$f': \mu(S';T'_{1},\ldots,T'_{k}) \rightarrow \pi(\sigma)T'$$
with $\zeta(f')=f$.

There is also a unique
morphisms $S\rightarrow S'$ with $\xi_n = id$, which gives a morphism 
$$\xi: \mu(S';S'_1,\ldots, S'_k)\rightarrow S$$
in $\h^n$. Every $S'_i$ has a unique tip. Hence, we  have a morphism 
$$\psi: \mu(S';S'_1,\ldots,S'_k)\rightarrow \mu(S';U_n,\ldots,U_n) \rightarrow S'$$ in $\h^n .$
Now $$\zeta(\xi)=\zeta(\psi):\mu([k];[1],\ldots,[1])\rightarrow [k]$$ is a retraction in $\H^{\infty}.$
So we get a chain of  morphisms 
$$\mu(S';T'_{1},\ldots,T'_{k})\leftarrow \mu(\mu(S';S'_1,\ldots,S'_k);T'_{1},\ldots,T'_{k}) \rightarrow
\mu(S;T'_{1},\ldots,T'_{k})$$
in $\h^n .$

We continue by choosing a unique morphism $\phi:T'\rightarrow T$ with $\phi_n = id$ and construct the other side of 
the chain analogously. Finally observe, that we have morphisms $\sigma'_i:T_i'\rightarrow T_i$ with 
$(\sigma'_i)_n = id$ which allow us to complete the construction. 

\

\Q

\

Recall that  $$\zeta^{\ast}:IO_{\infty}^{sym}(A) \rightarrow IO_n^{sym}(A).$$
means  the restriction functor along  $\zeta .$ 
By Theorem \ref{rep}, $\zeta^{\ast}$ is isomorphic to the functor of internal desymmetrisation $\delta_n .$

\begin{theorem}\label{formula} Let $A$ be a cocomplete categorical symmetric operad then  a left adjoint $sym_n$ to $\zeta^{\ast}$ exists,  and on  
   an internal
$n$-operad $a\in IO_n^{sym}(A)$, is given by the formula 
$$(sym_n(a))_k \ \simeq \ \mbox{\rm co}\!\lim\limits_{\h^{n}_k} \tilde{a}_k$$
where $\tilde{a}_k:\h^{n}_k \rightarrow A_k$ is the operadic functor representing the  operad $a$.
\end{theorem} 

\noindent {\bf Proof.} The case $n=1$ is well known. For an internal  symmetric operad $x$ the internal $1$-operad $\zeta^{\star}(x)$
has the same underlying collection as $x$ and the same multiplication for the orderpreserving maps of ordinals.
So the left adjoint to $\zeta^{\star}$ on object $a$ is given by 
$$(sym_1(a))_n = \coprod_{\Sigma_n} a_n$$
which is the same as the colimit of $\tilde{a}$ over $\h^1$ (see the description of $\h^1$ in Section \ref{H^nh^n}).

 Let $x:\h^n\rightarrow A$ be an operadic functor, $n\ge 2$. If we forget about the operadic structures on
$\h^{n},\H^{\infty}$ and
$A ,$ we can take a left Kan extension 
{\unitlength=1mm

\begin{picture}(60,33)

\put(45,25){\makebox(0,0){\mbox{$\h^n$}}}
\put(45,21){\vector(0,-1){10}}
\put(55,26){\shortstack{\mbox{$x $}}}

\put(69,25){\makebox(0,0){\mbox{$A $}}}

\put(49,25){\vector(1,0){16}}

\put(41,15){\shortstack{\mbox{$\zeta $}}}

\put(45,7){\makebox(0,0){\mbox{$\ \H^{\infty}$}}}

\put(49,11){\vector(3,2){16}}

\put(51.6,16.4){\line(1,1){4}}
\put(51,17){\line(1,1){4}}

\put(51,16.4){\line(1,0){2}}
\put(51,16.4){\line(0,1){2}}
\put(50.5,19.8){\shortstack{\mbox{$\scriptstyle\phi$}}}
\put(59,15){\shortstack{\mbox{$L = Lan_{\zeta}(x)$}}}

\end{picture}}

 \noindent  of \ $x$ \ along $\zeta$
in the
$2$-category of symmetric
$Cat$-collections. Since multiplication $m$ in $A$ preserves colimits in each variable, the
following diagram is a left Kan 
extension.

{\unitlength=1mm

\begin{picture}(60,33)(13,0)

\put(25,25){\makebox(0,0){\mbox{$\h^n_k\times\h^n_{n_1}\times \ldots \times \h^n_{n_k}$}}}
\put(25,21){\vector(0,-1){10}}
\put(104,26){\shortstack{\mbox{$m $}}}

\put(86,25){\makebox(0,0){\mbox{$A_k\times A_{n_1}\times \ldots \times A_{n_k} $}}}
\put(118,25){\makebox(0,0){\mbox{$A_{n_1+\ldots + n_k}$}}}
\put(76,13){\makebox(0,0){\mbox{$L_k\times L_{n_1}\times \ldots \times L_{n_k} $}}}

\put(42,25){\vector(1,0){25}}
\put(103,25){\vector(1,0){5}}


\put(25,7){\makebox(0,0){\mbox{$\H^{\infty}_k\times\H^{\infty}_{n_1}\times \ldots \times
\H^{\infty}_{n_k}$}}}

\put(45,11){\vector(3,1){28}}

\put(53.6,16.4){\line(1,1){4}}
\put(53,17){\line(1,1){4}}

\put(53,16.4){\line(1,0){2}}
\put(53,16.4){\line(0,1){2}}
\put(29,18.8){\shortstack{\mbox{$\scriptstyle\phi_k\times \phi_{n_1}\times \ldots \times
\phi_{n_k}$}}}




\end{picture}}

On the other hand, since $\zeta$ and $x$ are strict operadic functors  we have a natural
transformation 

{\unitlength=1mm

\begin{picture}(60,43)(13,-10)

\put(25,25){\makebox(0,0){\mbox{$\h^n_k\times\h^n_{n_1}\times \ldots \times \h^n_{n_k}$}}}
\put(25,21){\vector(0,-1){20}}
\put(104,26){\shortstack{\mbox{$m $}}}

\put(86,25){\makebox(0,0){\mbox{$A_k\times A_{n_1}\times \ldots \times A_{n_k} $}}}
\put(118,25){\makebox(0,0){\mbox{$A_{n_1+\ldots + n_k}$}}}

\put(42,25){\vector(1,0){25}}
\put(103,25){\vector(1,0){5}}


\put(25,-3){\makebox(0,0){\mbox{$\H^{\infty}_k\times\H^{\infty}_{n_1}\times \ldots \times
\H^{\infty}_{n_k}$}}}
\put(78,-3){\makebox(0,0){\mbox{$\H^{\infty}_{n_1+\ldots + n_k}$}}}  \put(108,-4){\makebox(0,0){\mbox{$,$}}}
\put(46,-3){\vector(1,0){19}}
\put(55,-2){\shortstack{\mbox{$\mu $}}}
\put(55,18){\shortstack{\mbox{$\mu $}}}

\put(85,1){\vector(3,2){30}}

\put(77,13){\makebox(0,0){\mbox{$\h^{n}_{n_1+\ldots + n_k}$}}}


\put(101,8){\shortstack{\mbox{$L$}}}

\put(92.6,9.4){\line(1,1){4}}
\put(92,10){\line(1,1){4}}

\put(92,9.4){\line(1,0){2}}
\put(92,9.4){\line(0,1){2}}
\put(91.5,13){\shortstack{\mbox{$\scriptstyle\phi$}}}

\put(40,20.5){\vector(4,-1){26}}
\put(72,10){\vector(0,-1){9}}
\put(80,14.5){\vector(4,1){27}}

\end{picture}}

\noindent and by the universal property of Kan extension we have a natural transformation
$$\rho: m(L_k; L_{n_1}, \ldots ,L_{n_k})\rightarrow L_n(\mu) $$
which determines a structure of lax-operadic functor on $L$. Moreover, $\phi$ becomes an
operadic natural transformation.   

Now the sequence of objects $L(p) = L([p]), p\ge  0 ,$ has a structure of an internal symmetric operad 
in
$A$. For a map of ordinals
$\sigma:[p]\rightarrow [k],$ let us define an internal multiplication $\lambda_{\sigma}$ by the composite
$$m(L(k);L(p_1),...,L(p_k)) \stackrel{\rho}{\longrightarrow} L(\mu([k],[p_1],...,[p_k]))
\longrightarrow L(\pi(\sigma)[p]) = \pi(\sigma)L(p).$$ 
Let us denote this operad by ${\cal L}(x)$.

 The calculation of $L(p)$ can be performed by the classical formula for 
pointwise left Kan extension \cite{ML}. 
It is therefore \ $\mbox{\rm co}\!\!\!\!\! \lim\limits_{f\in\zeta/[p]} \delta $ \ , where \ $\delta(f) =
x(S)$ \ for an object \ 
$f:\zeta(S)\rightarrow [p]$ \  of the comma category $\zeta/[p]$.  But according to the remark
after  Theorem \ref{rep}, \ $[p]$ \ is a terminal object of
$\H^{\infty}_p$ \ and therefore  $$\mbox{\rm co}\!\!\!\!\! \lim\limits_{ f\in\zeta/[p] } \delta  \
\simeq
\
\mbox{\rm co}\!\lim\limits_{\h^{n}_p} \ x_p \ .$$

It remains to prove that the internal operad ${\cal L}(\tilde{a})$ is $sym_n(a)$. Indeed, for a given  operadic morphism 
${\cal L}(\tilde{a})\rightarrow b$ the composite $$\tilde{a}\stackrel{\phi}{\rightarrow}\zeta^{\ast}{\cal L}(\tilde{a})\rightarrow
\zeta^{\ast}\tilde{b}$$ is operadic since $\phi$ is operadic. But $\zeta$ is final and, therefore,   the counit
of the adjunction $\zeta^{\ast}\dashv Lan_{\zeta}$ is an isomorphism. So for a given operadic morphism $\tilde{a}\rightarrow
\zeta^{\ast}\tilde{b}$ \  of internal
$n$-operads the morphism 
$${\cal L}(\tilde{a}) \rightarrow   {\cal L}(\zeta^{\ast}\tilde{b}) \simeq \tilde{b} $$
is operadic, as well. So the proof of the theorem is completed. 

\

\Q

\begin{cor}Let  
$A$ be an  
$n$-operad in   a cocomplete   symmetric monoidal category and $V,$  then 
$$(Sym_n(A))_k \ \simeq \ \mbox{\rm co}\!\lim\limits_{\h^{n}_k} \tilde{A}_k$$
where $\tilde{A}_k:\h^{n}_k \rightarrow V^{\bullet}$ is the  operadic functor representing the  operad $A$.
\end{cor}

\

\begin{theorem} The isomorphism $$\h^n \longrightarrow  Sym_n(\H^n)$$
induces
 a  canonical isomorphism 
 $$N(\h^n)\longrightarrow Sym_n(N(\H^n)).$$ \end{theorem}

\noindent {\bf Proof.} We have to calculate the result of the application of $Sym_n$ to the simplicial $Set$ $n$-operad ${\cal
F}_n^{\star}(1) = B({\cal F}_n,{\cal F}_n,1)$.  

We have  the following commutative diagram 

{\unitlength=0.9mm

\begin{picture}(60,35)(10,0)

\put(31,25){\makebox(0,0){\mbox{$ {\bf CO}_n(\H^n_d,Des_n({S}et^{\bullet}))$}}}
\put(31,21){\vector(0,-1){12}}

\put(86,25){\makebox(0,0){\mbox{${\bf CO}_n(\H^n,Des_n({S}et^{\bullet}))$}}}
\put(84,21){\vector(0,-1){12}}

\put(65,25){\vector(-1,0){14}}
\put(56,26){\shortstack{\small \mbox{$i^{\star}$}}}

\put(31,5){\makebox(0,0){\mbox{${\bf SCO}(\h^n_d,{S}et^{\bullet})$}}}

\put(86,5){\makebox(0,0){\mbox{${\bf SCO}(\h^n,{S}et^{\bullet})$}}}

\put(67,5){\vector(-1,0){16}}

\put(101,8){\shortstack{\small \mbox{$sev$}}}
\put(103,19){\shortstack{\small \mbox{$ev$}}}

\put(57,6){\shortstack{\mbox{\small $j^{\star} $}}}

\put(116,15){\makebox(0,0){\mbox{$O_n(Set)$}}}
\put(97,21){\vector(2,-1){10}}
\put(97,8){\vector(2,1){10}}

\end{picture}} 

\noindent where the vertical morphisms are canonical isomorphisms and the horizontal morphisms are the corresponding
restriction functors. In this diagram  $$sev: {\bf SCO}(\h^n,{\cal S}et)
\longrightarrow O_n(Set)$$ is the isomorphism  which gives an $n$-operad by evaluating an operadic functor on the generating
objects of $\h^n$ and $ev$ the corresponding evaluation functor for $n$-operads. 

Observe that the  simplicial operad ${\cal
F}_n^{\star}(1)$ is the result of application of
the functor
$ev$ to the  simplicial operadic functor  ${\bf Lan}^{\star}(\tilde{1}) = B({\bf Lan},{\bf Lan},\tilde{1})$  from Theorem
\ref{freeop}.

This diagram above shows that   $$ev({\bf Lan}^{\star}(\tilde{1})) = sev({\bf lan}^{\star}(\tilde{1})),$$ 
where the monad ${\bf lan}$ is the monad generated by the adjunction $j^{\star}\vdash Lan_j$, which is an operadic adjunction
by an argument analogous to the proof of the Theorem \ref{freeop}. We can also  prove
the following analogue for the formula (\ref{ext1}) :
\begin{equation}\label{ext2}{\bf lan}^k(\tilde{x})(W) = \coprod_{W\stackrel{{\scriptstyle
f}_{0}}{\leftarrow\!-}W_1\stackrel{{\scriptstyle f}_{1}}{\leftarrow\!-}W_2\stackrel{\scriptstyle f_2}{\leftarrow\!-}\
\ldots\stackrel{\scriptstyle f_{\scriptscriptstyle k-1}}{\leftarrow\!-}\ W_k}
\tilde{x}(W_k),\end{equation}
where $f_0,\ldots,f_{k-1}$ are morphisms in $\h^n$.   

Applying the functor $Lan_{\zeta}$ to the  operadic functor ${\bf lan}^k(\tilde{x})$
we have 
\begin{equation}\label{ext3}Lan_{\zeta}({\bf lan}^k(\tilde{x})) = Lan_{\zeta}(Lan_j(j^{\star}{\bf lan}^{k-1}(\tilde{x})))
\simeq
 Lan_{\zeta\cdot j}(j^{\star}{\bf lan}^{k-1}(\tilde{x})).\end{equation}
The left Kan extension $Lan_{\zeta\cdot j}(\tilde{x})$ is given by the formula
\begin{equation}\label{ext4}Lan_{\zeta\cdot j}(\tilde{x})(V) = \coprod_{\ \ \scriptstyle V\stackrel{f}{\longleftarrow} \
\zeta(W)} \tilde{x}(W),
\end{equation}
 where $f$ runs over the morphisms of $\H^{\infty}$. So, combining formulas (\ref{ext2}),
(\ref{ext3}) and (\ref{ext4}) we get
\begin{equation}\label{ext5}
Lan_{\zeta}({\bf lan}^k(\tilde{x}))(V) \simeq \coprod_{\ \ \scriptstyle V\stackrel{f}{\longleftarrow} \
\zeta(W)}
\coprod_{W\stackrel{{\scriptstyle
f}_{0}}{\leftarrow\!-}W_1\stackrel{{\scriptstyle f}_{1}}{\leftarrow\!-}W_2\stackrel{\scriptstyle f_2}{\leftarrow\!-}\
\ldots\stackrel{\scriptstyle f_{\scriptscriptstyle k-2}}{\leftarrow\!-}\ W_{k-1}}
\!\!\!\!\!\tilde{x}(W_{k-1})
\end{equation}
Now to calculate the $p$-th space of $Sym_n(N(\H^n)) =
Sym_n(sev({\bf lan}^{\star}(\tilde{1})))$ we have to put $x = 1$ and evaluate   (\ref{ext5}) at $V=[p]$. Since $[p]$ is a
terminal object  we have
$$Sym_n(N(\H^n)_p)^k \simeq \coprod_{W_0\stackrel{{\scriptstyle
f}_{0}}{\leftarrow\!-}W_1\stackrel{{\scriptstyle f}_{1}}{\leftarrow\!-}W_2\stackrel{\scriptstyle f_2}{\leftarrow\!-}\
\ldots\stackrel{\scriptstyle f_{\scriptscriptstyle k-2}}{\leftarrow\!-}\ W_{k-1}}\!\!\!\!\!\!\!\!\!\!\!\!\!\!\!
1 \ \ \ \ \ \ \ \ \ \ \ = N(\h^n_p)^k.$$
It is not hard to see that these isomorphisms agree with face and degeneracy operators.

\

\Q

\end{document}